\documentclass[12pt, a4paper, parskip=half, abstracton]{scrartcl}

\usepackage{array}
\usepackage{marginnote}
\usepackage{xcolor}\usepackage{hyperref}
\usepackage{amscd,amssymb,amsfonts,amsmath,latexsym,amsthm}
\usepackage[all,cmtip]{xy}
\usepackage{mathrsfs}
\usepackage{graphicx}
\usepackage{bm} 
 
\usepackage{etoolbox}

\def\hB{\hspace*{\fill}$\qed$}

\usepackage[nottoc]{tocbibind} 

\usepackage{defs_pp1}
\usepackage{slashed}
\usepackage[utf8]{inputenc}
\usepackage{microtype}
\usepackage[english]{babel}
\usepackage{mathtools}
\usepackage{bm}
\usepackage{esvect}
\usepackage[nameinlink, capitalise]{cleveref}

\title{$K$-theory of crossed products  via homotopy theory}
\author{
Ulrich Bunke\thanks{Fakult{\"a}t f{\"u}r Mathematik,
Universit{\"a}t Regensburg,
93040 Regensburg,
ulrich.bunke@mathematik.uni-regensburg.de} 
}

\numberwithin{equation}{section}
\newtheorem{theorem}{Theorem}[section] 
\newtheorem{prop}[theorem]{Proposition}
\newtheorem{lem}[theorem]{Lemma}

\newtheorem{ddd}[theorem]{Definition}
\newtheorem{kor}[theorem]{Corollary}

\theoremstyle{remark}
\theoremstyle{definition}

\newtheorem{ex}[theorem]{Example}
\newtheorem{rem}[theorem]{Remark}

\newcommand{\se}{\mathrm{se}}
\newcommand{\Nm}{\mathrm{Nm}}
 
\newcommand{\ee}{\mathrm{e}}
 \newcommand{\group}{\mathrm{group}}
\newcommand{\UCT}{\mathcal{B}}
\newcommand{\tR}{\mathrm{R}}
\newcommand{\All}{ \mathcal{A}\mathrm{ll}}
\newcommand{\Cyc}{\mathcal{C}\mathrm{yc}}

\newcommand{\fin}{\mathrm{fin}}

\newcommand{\EE}{\mathrm{E}}
\newcommand{\sepa}{\mathrm{sep}}
\newcommand{\LCH}{\mathbf{LCH}}

\newcommand{\kkG}{\mathrm{kk}^{G}}
\newcommand{\KKG}{\mathrm{KK}^{G}}

\newcommand{\nCalg}{C^{*}\mathbf{Alg}^{\mathrm{nu}}}
\newcommand{\can}{\mathrm{can}}

\newcommand{\Cof}{\mathrm{Cof}}
\newcommand{\Res}{\mathrm{Res}}

\newcommand{\Orb}{\mathbf{Orb}}

\newcommand{\Ob}{\mathrm{Ob}}

\newcommand{\Cofib}{\mathrm{Cofib}}

\newcommand{\Fib}{{\mathrm{Fib}}}

\newcommand{\incl}{\mathrm{incl}}

\newcommand{\cP}{\mathcal{P}}

\newcommand{\cK}{\mathcal{K}}

\newcommand{\CAlg}{{\mathbf{CAlg}}}

\newcommand{\bG}{{\mathbf{G}}}

\newcommand{\PSh}{{\mathbf{PSh}}}

\newcommand{\const}{{\mathtt{const}}}

\newcommand{\triv}{\mathrm{triv}}

\newcommand{\cO}{{\mathcal{O}}}

\newcommand{\Spc}{\mathbf{Spc}}

\newcommand{\Calg}{{\mathbf{C}^{\ast}\mathbf{Alg}}}

\renewcommand{\Pr}{\mathbf{Pr}}

\newcommand{\op}{\mathrm{op}}

\newcommand{\kk}{\mathrm{kk}}
\newcommand{\KK}{\mathrm{KK}}

\newcommand{\nCcat}{C^{*}\mathbf{Cat}^{\mathrm{nu}}}

\renewcommand{\Spc}{\mathbf{Spc}}

\newcommand{\exa}{\mathrm{ex}}

\newcommand{\K}{\mathrm{K}}
\renewcommand{\LCH}{\mathrm{LCH}}
 \newcommand{\hfin}{\mathrm{hfin}}

 \newcommand{\prp}{\mathcal{P}\mathrm{rp}}
\newcommand{\st}{\mathrm{st}}

\begin{document}

\maketitle\begin{abstract} In this paper, we analyse for a $G$-$C^{*}$-algebra $A$ to which extent one can calculate
the $K$-theory of  the reduced crossed product $K(A\rtimes_{r}G)$
from the $K$-theory spectrum $K(A)$   with the induced $G$-action. 
 We also consider some cases where one  allows to use
the $K$-theories of crossed products for some proper subgroups of $G$.
 Our central goal is to demonstrate the usefulness of a homotopy theoretic  approach.
We mainly concentrate on finite groups.
	\end{abstract}\tableofcontents
\setcounter{tocdepth}{5}

\section{Introduction}

Given a $C^{*}$-algebra $A$ with a continuous action    of a second countable locally compact group $G$,
one can form a new $C^{*}$-algebra called the reduced crossed product $A\rtimes_{r} G$.  
It is then a natural question  whether one can calculate the
$K$-theory  spectrum $K(A\rtimes_{r}G)$    in terms of 
the $K$-theory spectrum $K(A)$ 
 with the induced action of $G$. This will of course be only possible in very simple situations, but in more general cases one may hope to calculate $K(A\rtimes_{r}G)$ in terms of
 the spectra $K(A\rtimes_{r}H)$ for some proper subgroups $H$ of $G$. This  would lead to  an inductive approach.
\begin{ex}
The prototypical example is the case of the group 
 $G=\Z$. There is a canonical assembly map
$$\colim_{B\Z} K(A)\stackrel{}{\to} K(A\rtimes_{r}\Z)$$ which turns out to be an equivalence.
This  fact is a consequence of the validity  of the Baum-Connes conjecture with coefficients for the amenable group $\Z$ (see also \cref{wergjiwoeferfw}). However,  it can also be deduced from the more classical 
  Pimsner-Voiculescu long exact sequence \cite{pv} of $K$-theory groups
$$\dots\to K_{*}(A)\xrightarrow{\id-1_{*}} K_{*}(A) \to K_{*}(A\rtimes_{r}\Z)\to K_{*-1}(A)\ ,$$
where $1_{*}$ is the induced action of the generator $1$ of $\Z$ on $K_{*}(A)$. 
\hB
\end{ex}

\begin{ex}\label{wergjiwoeferfw}
For general locally compact groups $G$,   the Baum-Connes conjecture with coefficients (see e.g. \cite{MR4300553} for a recent survey) predicts for any $G$-$C^{*}$-algebra  $A$ an equivalence 
\begin{equation}\label{vqporkpoqcvxdwqe}\colim_{G_{\mathcal{C}\mathrm{pt}}\Orb} K^{G}_{A} \xrightarrow{\simeq}  K(A\rtimes_{r} G)\ .
\end{equation} 
Here $G_{\mathcal{C}\mathrm{pt}}\Orb$ is the full subcategory of the  orbit $\infty$-category of $G$ on orbits 
 with compact stabilizers, and $K_{A}^{G}$ is a  spectrum valued functor on $G\Orb$	 associated to $A$. It sends the orbit    $G/L$ to the spectrum $K_{A}^{G}(G/L)\simeq \KK^{G}(C_{0}(G/L),A)$. For  open (and closed, of course) subgroups $L$ the latter can be identified with  
 $  K(A\rtimes_{r}L)$, using the adjunctions \eqref{dfbvdfvdfvdfvsdfv} and \eqref{bsdfvkspdfvqreve}.  
 
 In the case of a discrete group $G$, the orbit category $G\Orb$ is ordinary and the construction of this functor  is  due to \cite{davis_lueck} (with corrections by \cite{joachimcat}), see also
 \cite{kranz} or \cite[Sec. 16]{bel-paschke}. It   was shown by \cite{kranz} (see also \cite[Sec. 16]{bel-paschke})  that the Baum-Connes-Higson assembly map  featuring the classical formulation of the Baum-Connes conjecture is equivalent to the Davis-Lück assembly  map appearing in \eqref{vqporkpoqcvxdwqe}.
 
 It is known that the Baum-Connes conjecture with coefficients holds for amenable groups, hyperbolic groups and many others, but it is false in general \cite{zbMATH01773367}.

 %
If the Baum-Connes conjecture with coefficients holds for $G$,  then it reduces the problem of calculating $K(A\rtimes_{r}G)$ 
to the calculation of the spectra $K(A\rtimes_{r}  L)$ for all compact subgroups   $L$ of $G$ and the homotopy theoretic problem of the 
 calculation of a colimit of a spectrum-valued functor over $G_{\mathcal{C}\mathrm{pt}}\Orb$.
 \hB
 \end{ex}

As explained in  \cref{wergjiwoeferfw}, the Baum-Connes conjecture takes the $K$-theories of crossed products by compact groups as an input in order to calculate the $K$-theory of crossed products by non-compact groups.
In contrast, in the present paper we are interested in the calculation of the $K$-theory of crossed products by compact groups themselves using homotopy theoretic methods. We will mainly concentrate on finite groups, i.e., the case of discrete compact groups.

We will first consider cases where one can calculate the spectrum $K(A\rtimes_{r}G)$ from the spectrum $
 K(A)$ with its induced $G$-action. Later we study the more complicated situation 
 where one also must take crossed products by subgroups of $G$ into account.

\begin{ex} \label{wtkohpwerfrefwref}
Assume that $G$ is  discrete and that
 the $G$-$C^{*}$-algebra $A$ is induced from the trivial subgroup, i.e., $A\cong C_{0}(G,B)$ with the $G$-action by left translations for some $C^{*}$-algebra $B$. Then
we have an equivalence \begin{equation}\label{vssodjvposvadvdscadca}\colim_{BG} K(A)\simeq K(A\rtimes_{r}G)\ .
\end{equation}
In order to see this, one can calculate both sides explicitly.
For the right-hand side,  we  employ Green's imprimitivity theorem:
$$K(A\rtimes_{r}G)\stackrel{\eqref{bwefljrfopwerfrerfwerfwreferf}}{\simeq} K(B)\ .$$
 For the left-hand side, we observe that 
$A\cong \bigoplus_{G}B$  in $\Fun(BG,\nCalg)$ with the $G$-action by left translations on the index set.  Since the $K$-theory functor preserves sums,  we get 
$K(A)\simeq \bigoplus_{G}K(B)$ in $\Fun(BG,\Mod_{KU}(\Sp))$ and can conclude that 
$$\colim_{BG} K(A)\simeq  K(B)\ .$$ \hB
\end{ex}
If $G$ is finite and  $A$ is a  $G$-$C^{*}$-algebra, then the Rohklin property     (in the form introduced by \cite{Izumi_2004}, \cref{rktohoprthrgtgegetrg}) says that $A$ is induced from the trivial subgroup in an asymptotic sense.
 The following proposition can thus be considered as an extension of  \cref{wtkohpwerfrefwref}.
  \begin{prop}[\cref{gkerpogkpwergwerfwerf1}]\label{gkerpogkpwergwerfwerf}
If $G$ is finite and $A$ is a $G$-$C^{*}$-algebra with  the Rokhlin property, then the assembly map is an equivalence
$$\colim_{BG} K(A)\stackrel{\simeq}{\to}  K(A\rtimes_{r}G)\ .$$
\end{prop}
 A variety of examples of  $G$-$C^{*}$-algebras with the Rokhlin property is given in \cite{Izumi_2004}, see also 
   \cref{weotkgperfgrefrefwrefwf} below.

\begin{rem}
As indicated in the examples above, we consider a homotopy theoretic formula
for the spectrum $K(A\rtimes_{r}G)$ as an eligible solution.  Our result reduces the
 calculation of the $K$-theory groups $K_{*}(A\rtimes_{r}G)$ to  homotopy theoretic calculations
 which still might be a complicated task involving  spectral sequences and similar; see for example, the open problem mentioned in \cref{jerotgrgertgertgfb}. But there is are situations  that will be discussed in   detail in \cref{oierjgoiegwergwerg} below where this calculation is straightforward.
See also \cref{iorjfoqrggegerfrefwref} and \cref{iorjfoqrggegerfrefwref1} for further concrete cases.
   \hB
\end{rem}

A natural framework for the $G$-equivariant homotopy theory of $C^{*}$-algebras  is the stable $\infty$-category $\KK^{G}$ representing equivariant Kasparov $KK$-theory. We will  review the basic features of this category 
 in \cref{wtkopgwfrefrefw}. In particular, it  has a symmetric monoidal structure whose tensor unit will be denoted by $\beins_{\KK^{G}}$.   The commutative ring spectrum \begin{equation}\label{lekjkqrlfqwedqewdqde}R(G):=\map_{\KK^{G}}(\beins_{\KK^{G}},\beins_{\KK^{G}})
\end{equation}  is the    homotopy theoretic version of the representation ring of $G$. As  an object in presentable stable $\infty$-categories $\Pr^{L}_{\mathrm{st}}$, the category $\KK^{G}$   becomes a commutative algebra in modules over the commutative algebra $\Mod_{R(G)}(\Sp)$ in $\Pr^{L}_{\mathrm{st}}$,  the presentably symmetric monoidal category of modules over the commutative ring spectrum $R(G)$. This is a neat way to say that $\KK^{G}$  has a highly coherent  enrichment over $\Mod_{R(G)}(\Sp)$. In particular, 
 the mapping spectra in $\KK^{G}$ become $R(G)$-modules.   
 The $\Mod_{R(G)}(\Sp)$-module structure of $\KK^{G}$ provides 
 a bifunctor \begin{equation}\label{vsdfvdsvsvfdvsfdvsfdvre}-\otimes-:\Mod_{R(G)}(\Sp)\times \KK^{G}\to \KK^{G}
\end{equation} 
 which preserves colimits in each argument and is characterized by the property that $R(G)\otimes -$ is equivalent to identity. In particular, any element $\sigma$ in $\pi_{0}R(G)$ gives rise to an object-preserving endofunctor of $\KK^{G}$,
 i.e., it acts as an endomorphism on any object in a way which is compatible with morphisms.
 
 The equivariant $K$-theory functor for $C^{*}$-algebras  factorizes over $\KK^{G}$:
 $$K^{G}:G\nCalg\xrightarrow{\kk^{G}} \KK^{G}\xrightarrow{\K^{G}} \Mod_{R(G)}(\Sp)\ ,$$
 where \begin{equation}\label{hertoijiojgretbth}\K^{G}:=\map_{\KK^{G}} (\beins_{\KK^{G}},-)\ .
\end{equation}
 The $R(G)$-module structure on $K$-theory plays an important role in this note since we want to consider
 localizations  and completions at elements in $\pi_{0}R(G)$.
 In the case of the trivial group, we will omit the superscript $G$
 in the notations for $K$ and $KK$-theory functors.

 The homotopy theoretic point of view leads to an immediate generalization of \cref{wtkohpwerfrefwref}.
 We shall see  that the map inducing the equivalence \eqref{vssodjvposvadvdscadca} comes from a natural transformation  \eqref{frewfwrewoifjowrfwrefwerf}  of colimit preserving functors    from
 $\KK^{G}$ to $\Mod_{KU}(\Sp)$.   The equivalence can therefore be extended from induced $C^{*}$-algebras to   $G$-$C^{*}$-algebras   
  whose $\KK^{G}$-class belongs  to the localizing subcategory  $\langle \Ind^{G}(\KK) \rangle $  generated by
 the image of the induction functor $\Ind^{G}:\KK\to \KK^{G}$.
 Even more generally, we could take the $\KK^{G}$-classes themselves as parameters. 
 Ideas of this kind  in the context of the triangulated equivariant $KK$-category were first  used in \cite{MR2193334} in connection with the Baum-Connes conjecture.
  
  For $A$ in $\KK^{G}$, we use the notation  $\widehat \Res^{G}(A)$ for the object  in $\Fun(BG,\KK)$ representing  the  underlying  $\KK$-class of $A$ with the induced $G$-action, see \eqref{qfqwefqwedqwdq} for the precise definition. Note that for a $G$-$C^{*}$-algebra, we keep writing $ K(A)$  for the $K$-theory spectrum with the induced $G$-action.

  \begin{prop}[\cref{jeirgowregwerfrefw1}]\label{jeirgowregwerfrefw}
 Assume that $G$ is discrete. If $A$ is in $\langle \Ind^{G}(\KK) \rangle$, then we have an equivalence
\begin{equation}\label{vferopvjopfdewdedq}\colim_{BG} \K(\widehat \Res^{G}(A)) \simeq \K(A\rtimes_{r}G)\ .
\end{equation}  
 \end{prop}

 \begin{rem} For a $G$-$C^{*}$-algebra $A$,  \cref{jeirgowregwerfrefw} implies the equivalence \eqref{vssodjvposvadvdscadca}
  under the assumption that the object $\kk^{G}(A)$ can be obtained as a colimit  of a diagram of
  objects in the image of the induction functor. The case considered in  \cref{wtkohpwerfrefwref}
  corresponds to the simple situation where $\kk^{G}(A)$ itself is in this image. 
 
  The statement of \cref{jeirgowregwerfrefw}  in the form of an isomorphism  
  \begin{equation}\label{fqweqdewdedeq}\colim_{W\subseteq EG}\KK^{G}_{*}(C_{0}(W) ,A)\stackrel{\simeq}{\to} K_{*}(A\rtimes_{r}G)
\end{equation} 
of $K$-theory groups
  for $G$-$C^{*}$-algebras $A$ such that $\kk^{G}(A)\in \langle \Ind^{G}(A)\rangle$ is a well-known fact. Here the filtered
 colimit runs over the finite $G$-CW-subcomplexes of a $G$-CW model for the classifying space $EG$ of $G$.
   But note that the identification of the left-hand sides of \eqref{fqweqdewdedeq} and \eqref{vssodjvposvadvdscadca} is a non-trivial matter and similar to the identifications of the domains of the Baum-Connes and the Davis-Lück assembly maps  in \cite{kranz}, or using Paschke duality in \cite{bel-paschke}.   We refer to 
     \cref{wejkogpefrefewrfewrfewrfrewfw} for further references.

  Even more generally, using the notation introduced in \cref{qewfoiihqowedewqdeqwd} and the method from the proof of \cref{gkerpogkpwergwerfwerf}, one can show that  the equivalence \eqref{vssodjvposvadvdscadca}
 holds  for $A$ in $\KK^{G}$ under the assumption that $c^{G}(A)$ (see \eqref{wrfqwdewdewfqef} for $c^{G}$) belongs to  closure of $ \langle \Ind^{G}(\EE_{\oplus})\rangle$ under universal  phantom retracts. 
  \hB
 \end{rem}

    The following proposition characterizes a class of  examples where the assumption of \cref{jeirgowregwerfrefw} is satisfied.
  If $n$ is an integer, then it induces an endomorphism $n\id_{A}$  of any object $A$ of a pre-additive (so in particular of a stable) $\infty$-category. If this endomorphism is an equivalence, then we say that $n$ acts as an equivalence on $A$. Let $\rho $ in $\pi_{0}R(G)$ be the class of the left-regular representation of $G$ on $L^{2}(G)$ in the representation ring. Then $\rho$ acts as an endomorphism of $A$ in $\KK^{G}$ and we could consider the condition that this endomorphism is an equivalence.

Let $A$ be an object of $\KK^{G}$.
\begin{prop}[\cref{qirjfofdewdewdqewd1}]\label{qirjfofdewdewdqewd} \footnote{The author thanks S. Nishikawa for a discussion about examples helped to locate a mistake in the statement of this proposition in an earlier version of this paper.}
We assume that $G$ is finite. 
 The following assertions are equivalent.
\begin{enumerate}
\item \label{qirjfofdewdewdqewdf} $\rho$ acts as an equivalence on $A$.
\item We have  $A\in  \langle \Ind^{G}(\KK) \rangle$ and $|G|$ acts an equivalence on $A$.
\end{enumerate}
%
%
\end{prop}


\begin{rem}
%
%

Under the equivalent conditions of 
  \cref{qirjfofdewdewdqewd}, we have the equivalence\begin{equation*}\label{vferopvjopfdewdedqfff}\colim_{BG} \K(\widehat \Res^{G}(A)) \simeq \K(A\rtimes_{r}G)\ .
\end{equation*}  In addition, we know
   that $\K(\widehat \Res^{G}(A))$ is projective by \cref{ekorgpefreqwdqedqed} so that we also have  an isomorphism of $K$-theory groups $$ \colim_{BG} \K_{*}(\widehat \Res^{G}(A))\simeq \K_{*}(A\rtimes_{r}G)\ .$$ \hB
 \end{rem}

%
%


In general, the  spectrum  $K(A\rtimes_{r}G)$  is not determined by the spectrum $K(A)$ with its induced $G$-action alone. We now start to
consider  situations where for the calculation of $K(A\rtimes_{r}G)$  we need to take the contributions  of crossed products for  some proper  subgroups of $G$ into account. To this end, 
in \cref{wtgijowergrwefwrefwef} we   construct a functor $$\K^{G}_{-}:\KK^{G}\times G\Orb^{\comp}\to \Mod_{R(G)}(\Sp)\ , \quad (A,S)\mapsto \K^{G}_{A}(S)\ .$$
It is an analogue of the functor $K^{G}_{A}$ appearing in \cref{wergjiwoeferfw}. In particular, for a compact group $G$ and a  cofinite subgroup $L$,  by \cref{qerjigoiewrgwerfrefw} it also has the values $\K^{G}_{A}(G/L)\simeq \K(A\rtimes_{r}L)$. This equivalence equips the $KU$-modules $\K(A\rtimes_{r}L)$ with $R(G)$-module structures.

  If $G$ is finite, then by \cref{rtkohpertgregertgetg}  the $KU$-module  $\K^{G}_{A}(G)$ with the $G$-action induced by right-multiplication on $G$ is equivalent to   $\K(\widehat \Res^{G}(A))$, i.e., the $KU$-module $\K(A)$  with $G$-action induced by functoriality. The functor $\K^{G}_{A}$ thus  interpolates between the  spectrum $\K(\widehat \Res^{G}(A))\simeq \K^{G}_{A}(G)$ with $G$-action and the spectrum $\K(A\rtimes_{r}G)\simeq \K^{G}_{A}(*)$.

We now restrict to the case of finite groups $G$.    Let $G_{\prp}\Orb$ denote the full subcategory of the orbit category of $G$  on the orbits whose stabilizer is a proper subgroup of $G$.
 We then analyse  the cofibre sequence   \begin{equation}\label{vwvssfdvsfvs}\colim_{G_{\prp}\Orb} \K^{G}_{-} \to \K^{G}_{-}(*) \to \Cof^{G}(-)
\end{equation}
  of functors from $\KK^{G}$ to $\Mod_{R(G)}(\Sp)$, where  
 $ \Cof^{G}$ is defined as the cofibre. The  evaluation of
 the first map  at $A$ in $\KK^{G}$ is  called the assembly map for the functor $\K^{G}_{A}$ and the family $\prp$ of proper subgroups.  Since  $\K^{G}_{A}(*)\simeq \K(A\rtimes_{r} G)$,  this sequence suggests
 an inductive approach to the calculation of $\K(A\rtimes_{r} G)$ in terms of the cofibre and
 a colimit of $K$-theories of crossed products $\K^{G}_{A}(G/L)\simeq \K(A\rtimes_{r}L)$ by the proper subgroups $L$ of $G$.   
 We are in particular interested in situations where this cofibre is trivial.
 
 Here is a very simple instance where this is the case. We let $\langle \Ind_{\prp}^{G} \rangle$ denote the localizing subcategory of $\KK^{G}$ generated by the collection of $\Ind_{L}^{G}(\KK^{L})$ for all proper subgroups $L$ of $G$.

  \begin{prop}[\cref{werjoigpwegerfwref1}]\label{werjoigpwegerfwref} If $G$ is finite and  $A$ is in $\langle \Ind_{\prp}^{G} \rangle$, then  the assembly map is an equivalence
   \begin{equation}\label{vsdfvsdcscdsf}\colim_{G_{\prp}\Orb} \K^{G}_{A}\xrightarrow{\simeq} \K(A\rtimes_{r}G)\ .
\end{equation}
  \end{prop}
 
 \begin{rem} If $A=\Ind_{L}^{G}(B)$ for some $B$ in $\KK^{L}$ and proper subgroup $L$ of $G$, then 
  the right-hand side of \eqref{vsdfvsdcscdsf} can of course be calculated by Green's imprimitivity theorem: $$\K(B\rtimes_{r}L)\stackrel{\eqref{bwefljrfopwerfrerfwerfwreferf}}{\simeq}  \K(A\rtimes_{r}G)$$
 and the spectrum on the left is one of the values of the functor contributing to the colimit in  \eqref{vsdfvsdcscdsf}. So the \cref{werjoigpwegerfwref} is not very useful for the calculation of the right-hand side of \eqref{vsdfvsdcscdsf} in terms of the colimit, but it might be used  in the reverse direction to extract some information about the structure of the colimit on the left-hand side.
\hB\end{rem}

 We continue to assume that    $G$ is finite.
 We  let 
  $V$ be the orthogonal complement  of the trivial representation in the finite-dimensional $G$-representation $L^{2}(G)$.
We then define the element $\xi:=\Lambda_{-1}(V)$ in $\pi_{0}R(G)$.
In \cref{weijgowiergijrefrewf9}, we further define the endofunctors $(-)[\xi^{-1}]$ and $S_{\xi}:=\Fib(\id\to (-)[\xi^{-1}])$ of any  module category in $\Pr^{L}_{\mathrm{st}}$ over
$\Mod_{R(G)}(\Sp)$. 
 \begin{theorem}[\cref{wekopgwerfrewfw1}]\label{wekopgwerfrewfw}
 We assume that $G$ is finite.
The cofibre sequence \eqref{vwvssfdvsfvs} is equivalent to the cofibre sequence
\begin{equation}\label{sfvsdfvsfdcssc}S_{\xi}(\K^{G}_{-}(*))\to \K^{G}_{-}(*)\to  \K^{G}_{-}(*)[\xi^{-1}]\end{equation}  of functors from $\KK^{G}$ to $\Mod_{R(G)}(\Sp)$.
\end{theorem}
 The main insight given by \cref{wekopgwerfrewfw} is that $\xi$ acts an equivalence on the values of the cofibre in \eqref{vwvssfdvsfvs}. If $\xi=0$, then we could conclude that $\Cof^{G}\simeq 0$.
 In \cref{kopthethegtrg1}, we will observe that  if $G$ is finite and not cyclic, then $\xi=0$.



 \begin{kor}\label{gjkwegokwerpferfwerf}
 If $G$ is finite and not cyclic, then the assembly map is an equivalence. 
 $$\colim_{G_{\prp}\Orb}\K^{G}_{A}\stackrel{\simeq}{\to} \K(A\rtimes_{r}G)$$
 for every $A$ in $\KK^{G}$.
  \end{kor}

 Let $G_{\Cyc}\Orb$ denote the full subcategory of the orbit category of $G$ on orbits whose stabilizer is a cyclic subgroup. By an inductive procedure,  using \cref{gjkwegokwerpferfwerf} in the induction step, we can then deduce:
\begin{prop}[\cref{geropfrfwrefrefrwfwrf1}]\label{geropfrfwrefrefrwfwrf} If $G$ is finite, then
we have an equivalence
$$\colim_{G_{\Cyc}\Orb}\K^{G}_{A}\xrightarrow{\simeq} \K(A\rtimes_{r}G)$$ for every $A$ in $\KK^{G}$.
\end{prop}

This result reduces the calculation of  $\K(A\rtimes_{r}G)$ to the calculation of  the $K$-theory of crossed products by cyclic groups
and  the  colimit  of a spectrum-valued functor over $G_{\Cyc}\Orb$.
For $A=\beins_{\KK^{G}}$, this result is similar in spirit to the classical Artin- and Brower induction theorems.


\begin{rem}
\cref{geropfrfwrefrefrwfwrf} can also be deduced from the results of \cite{Mathew_2019}. To this end one observes that 
$\K^{G}_{A}$ (considered as a $\Sp$-valued functor) refines to spectral Mackey functor\footnote{The extension of $\K^{G}_{A}$ to a Mackey functor will be discussed in a subsequent paper.}.   This refinement has furthermore a lax symmetric monoidal refinement in the variable $A$  and therefore becomes a module over $\K^{G}_{\beins_{\KK^{G}}}$. The latter spectral Mackey functor  can\footnote{In order to make this precise we would need a precise definition of what the usual usual genuine equivariant topological $K$-theory spectrum $KU^{G}$ as a commutative algebra in spectral Mackey functors is. The reference \cite{Mathew_2019}  is vague at this point.} be identified with the usual genuine equivariant topological $K$-theory spectrum $KU^{G}$.   The assertion of \cref{geropfrfwrefrefrwfwrf}  for $ \K^{G}_{\beins_{\KK^{G}}}$ is then the same as  \cite[Prop. 5.7]{Mathew_2019}. It is easy to see that it implies \cref{geropfrfwrefrefrwfwrf} also for all modules over $\K^{G}_{\beins_{\KK^{G}}}$.
 \hB
\end{rem}

 In the following, we consider situations where the cofibre in the sequence \eqref{vwvssfdvsfvs}, or equivalenty the third term in \eqref{sfvsdfvsfdcssc}, may be non-trivial but splits off as a summand. Let $p $ be a prime, $k$ be in $\nat\setminus \{0\}$,  and let $A$ be  in $\KK^{C_{p^{k}}}$. 
 
 \begin{prop}[\cref{gijeoggerwfreffvvsdfvsfdvsfv1}]\label{gijeoggerwfreffvvsdfvsfdvsfv}
 If $p$ acts as an equivalence in $A$, then the  sequence   
 $$S_{\xi}(\K^{C_{p^{k}}}_{A}(*))\to \K^{C_{p^{k}}}_{A}(*)\to  \K^{C_{p^{k}}}_{A}(*)[\xi^{-1}]$$ splits naturally.
 \end{prop}   
 Under the assumptions of \cref{gijeoggerwfreffvvsdfvsfdvsfv},
 we thus have  an equivalence 
 \begin{equation}\label{fewqiofhiwodewdedqwedewd}
S_{\xi}(\K^{C_{p^{k}}}_{A}(*)) \oplus   K(A\rtimes_{r}C_{p^{k}})[\xi^{-1}]\simeq K(A\rtimes_{r} C_{p^{k}})
\end{equation}
If $k=1$, then $$S_{\xi}(\K^{C_{p}}_{A}(*)) \simeq \colim_{BC_{p}}  \K(\widehat \Res^{C_{p}}(A)))\ ,$$ and the
  homotopy groups  of this colimit can be calculated using
  \cref{ekorgpefreqwdqedqed}: \begin{equation}\label{efqewdqewdqededqed} \pi_{*}(S_{\xi}(\K^{C_{p}}_{A}(*)))\cong 
\pi_{*}(\colim_{BC_{p}} \K(\widehat \Res^{C_{p}}(A)))\cong \colim_{BC_{p}} \K_{*}(\widehat \Res^{C_{p}}(A))\ .
 \end{equation}
 For $k>1$ the calculation of this summand is more complicated and involves a colimit over $C_{p^{k},\prp}\Orb$.
 
 In general, the homotopy groups of the second summand $K(A\rtimes_{r}C_{p^{k}})[\xi^{-1}]$ in \eqref{fewqiofhiwodewdedqwedewd}
 are still complicated to calculate.  In the following we provide an interesting class of examples
 where this calculation is possible.

 \begin{ex} \label{iorjfoqrggegerfrefwref}This example is added in order to demonstrate that the calculation of the $K$-theory of crossed products using the decomposition  \eqref{fewqiofhiwodewdedqwedewd} confirms known results. We consider the group $G:=C_{p}$ .
The simplest object of $\KK^{C_{p}}$ on which $p$ acts as an equivalence is of course
 $A=\beins_{\KK^{C_{p}}}[p^{-1}]$.  Then  $\K(\widehat \Res^{C_{p}}(A))\simeq \triv^{C_{p}}(KU[p^{-1}])$  and the 
 first summand in   \eqref{fewqiofhiwodewdedqwedewd}  is given by $\colim_{BC_{p}}\triv^{C_{p}}(KU[p^{-1}])\simeq KU[p^{-1}]$.
  The second summand is given as a $KU$-module  by
 $\K(\beins_{\KK^{C_{p}}}[p^{-1}]\rtimes_{r}C_{p})[\xi^{-1}]\simeq R(C_{p})[p^{-1}][\xi^{-1}]\simeq  KU[p^{-1}]^{\oplus (p-1)}$.  Consequently, we have
 \begin{equation}\label{fvsdfvdfvfdvsdfv}\K(\beins_{\KK^{C_{p}}}[p^{-1}]\rtimes_{r}C_{p})\simeq   KU[p^{-1}]^{\oplus p} \ .
\end{equation} 
 This shows that the calculation based on \eqref{fewqiofhiwodewdedqwedewd} is compatible with the well-known general fact that for any finite group $G$,  we have an equivalence of $KU$-modules  $$\K(\beins_{\KK^{G}}\rtimes_{r}G)\simeq R(G)\simeq KU^{\bigoplus |\hat G|} \ . $$
This equivalence persists after inverting $p$ and yields the equivalence \eqref{fvsdfvdfvfdvsdfv}  in the case $G=C_{p}$.
  \hB 
 \end{ex}

  \begin{ex} \label{iorjfoqrggegerfrefwref1}
  In this example we demonstrate the unitility of the results of the present paper for the explicit calculation of the $K$-theory of a crossed product  of a rather complicated  $A$ in $\KK^{G}$.
If we knew how $A$  is built from objects induced from smaller subgroups, then we could reduce the calculation of $\K(A\rtimes_{r}G)$ to the calculation of crossed products for smaller subgroups. In contrast, multiplicative induction provides another source of objects of $\KK^{G}$ which turn out to be much more difficult to understand.

Assume that $G$ is a finite group and $H$ is a subgroup of $G$.
 On the level of $C^{*}$-algebras, the multiplicative induction sends a $H$-$C^{*}$-algebra $B$ to the $G$-$C^{*}$-algebra given by the tensor product   $B^{\otimes G/H}:=\bigotimes_{G/H}B$ with an appropriate $G$-action \eqref{dafasdfqffadfs}.  Thereby the choices of the minimal or maximal tensor products leads to two different versions of   multiplicative induction.

In order to work on the level of homotopy theory, in  \cref{iorgoergfgsfg} we show that the multiplicative induction
descends to a  symmetric monoidal functor $$(-)^{\otimes G/H} : \KK^{H}\to \KK^{G}$$
which is not exact in general.

 If $I\to A\to B$ is a fibre sequence in $\KK$ and $Z$ is a finite  $G$-set, then the relation of 
$A^{\otimes Z}$ in $\KK^{G}$ with  multiplicative inductions  of $I$ and $B$ is given by a tower that 
will be discussed in detail in \cref{krwfopqwewfewdffqwf}. The associated spectral sequence is employed, e.g., in the proof of \cref{wergojowpergrefgrwefwref} below.

We consider the group $G:=C_{p^{k}}$ for a prime $p$ and $k\in \nat\setminus\{0\}$. We further choose an integer $\ell$ such that $p\not| \ell$ and   consider $B:=\beins_{\KK^{H}}/\ell$ in $\KK^{H}$, where $H$ is any subgroup of $C_{p^{k}}$ for the moment. The assumption on $\ell$ implies that $p$ acts as an equivalence on $B$. Then by \cref{egrjiweorgerferwfrewf} it 
  also acts as an equivalence on  $A:=B^{\otimes C_{p^{k}}/H}$. Therefore   
 \cref{gijeoggerwfreffvvsdfvsfdvsfv} applies to $A $.
 
 If one likes to think in concrete terms, note that
 $$A\simeq \kk^{C_{p^{k}}}(\underbrace{\cO_{\ell+1} \otimes \dots \otimes \cO_{\ell+1}}_{C_{p^{k}}/H})$$ with the $C_{p^{k}}$ action by permutations of tensor factors, and where $\cO_{\ell+1}$ is the Cuntz algebra.

  The following calculates the second summand in the decomposition \eqref{fewqiofhiwodewdedqwedewd} as a $KU$-module.
\begin{lem}[Lemma \ref{oejgoerpwgrefdssfvfsfsv1}]\label{oejgoerpwgrefdssfvfsfsv}
We have $$\K(A\rtimes_{r}C_{p^{k}})[\xi^{-1}]\simeq R(C_{p^{k}})[\xi^{-1}]/\ell\simeq (KU/\ell)^{\oplus (p-1)p^{k-1}}\ .$$ 
\end{lem}
Note that the result is independent of $H$.
The crucial fact enabling us to do this calculation is that the composition of the  cofibre functor
$\Cof^{C_{p^{k}}}$  from \eqref{vwvssfdvsfvs} with the multiplicative induction functor preserves all colimits by \cref{ewgojwoepgfrefwfwerfwrf}, so is in particular exact


  In order to calculate the first summand in the decomposition \eqref{fewqiofhiwodewdedqwedewd}, we restrict to the case $k=1$, the trivial subgroup $H$, and use \eqref{efqewdqewdqededqed}.
  We therefore must calculate the groups with $C_{p}$-action  $\K_{*}(\widehat \Res^{C_{p}}(A))$  and form coinvariants.
  
  \begin{lem}[\cref{wergojowpergrefgrwefwref1}]\label{wergojowpergrefgrwefwref}
  We have $$\colim_{BC_{p}} K(\widehat \Res^{C_{p}}(A))\simeq (KU/\ell)^{n_{0}(p)}\oplus (\Sigma KU/\ell)^{n_{1}(p)}$$ with
  $$n_{i}(p):=\left\{\begin{array}{cc} \frac{2^{p-2} +  \frac{(p+1)(p-1)}{2}  }{p}&i=0\\ \frac{2^{p-2}  - \frac{(p-1)^{2}}{2}   }{p}& i=1 \end{array} \right.$$ for $p>2$ and
  $$n_{i}(2):=\left\{\begin{array}{cc} 1&i=0\\ 0&i=1  \end{array} \right.\ .$$
  \end{lem}

Combining these two results with \eqref{fewqiofhiwodewdedqwedewd} we get:  
\begin{kor} For a prime $p$, an integer $\ell$ with $p\not|\ell$, and $A:=(\beins_{\KK}/\ell)^{\otimes C_{p}}$ in $\KK^{C_{p}}$ we have 
 $$\K(A\rtimes_{r}C_{p})\simeq (KU/\ell)^{n_{0}(p)+(p-1)}\oplus  (\Sigma KU/\ell)^{n_{1}(p)}\ .$$
\end{kor}

We have the following table:   
$$\begin{array}{|c||c|c|c|c|c|c|}  \hline   p&2&3&5&7&11&\dots\\ \hline n_{0}(p)+p-1&2&4&8&14&62&\dots\\ \hline  n_{1}(p)&0&0&0&2&42 &\dots\\\hline  \end{array}$$
The case $p=2$ is due to \cite{Izumi}, and the case $p=3$ has previously been obtained using different methods by Kranz and Nishikawa.\footnote{These calculations are part of an ongoing joint project with J. Kranz and S. Nishikawa 
with the main goal to settle the case where $p=\ell$, see \cref{jerotgrgertgertgfb}.} 
         \hB
 \end{ex}

%
%
%
%

We now come back to the question how one can get rid of the  cofibre in the sequence \eqref{vwvssfdvsfvs}, or equivalenty the third term in \eqref{sfvsdfvsfdcssc}, by either adding additional assumptions on $A$, or by
applying auxiliary operations.

By   \cref{erjigowergfregfw}.\ref{werkjgowerferwfqewdqewdacf},  we have a Bousfield localization
$$L_{\xi}:\Mod_{R(G)}(\Sp)\leftrightarrows L_{\xi}\Mod_{R(G)}(\Sp):\incl$$
whose kernel consists precisely of the $\xi$-acyclic objects in $\Mod_{R(G)}(\Sp)$, i.e., the objects  $C$ satisfying $S_{\xi}(C) \simeq 0$. Recall that $\K(A\rtimes_{r}G)$ has a canonical $R(G)$-module structure.
If we combine \cref{wekopgwerfrewfw} with  \cref{erjigowergfregfw}.\ref{wekgopwegrfwerfrf} we get:
\begin{kor}\label{wjogpwergerfrefwerfw}
We assume that $G$ is finite. We have an equivalence
$$L_{\xi} (\colim_{G_{\prp}\Orb} \K^{G}_{A})\xrightarrow{\simeq} L_{\xi} (\K(A\rtimes_{r}G))$$
for every $A$ in $\KK^{G}$.
\end{kor}
Thus, up to $\xi$-completion, we can calculate $\K(A\rtimes_{r}G)$ in terms of the $K$-theory of  crossed products of $A$ by proper subgroups of $G$ and a colimit of a spectrum-valued functor over
$G_{\prp}\Orb$.

If one is interested in $p$-torsion effects for a prime $p$, then it is natural to consider the $p$-completion $L_{p}$. 
 We again consider   the group $G=C_{p^{k}}$ for $k$ in $\nat\setminus \{0\}$.  For 
  $A$  in  $\KK^{C_{p^{k}}}$ 
  we  can  indeed  calculate
the $K$-theory spectrum $\K(A\rtimes_{r} C_{p^{k}})$ after $p$-completion in terms of the spectrum $\K(\widehat \Res^{C_{p^{k}}}(A))$, i.e., the underlying $K$-theory spectrum of $A$ with the induced $C_{p^{k}}$-action. 
 \begin{theorem}[\cref{erijgowergrefwrefr1}]\label{erijgowergrefwrefr} For any $A$ in $\KK^{C_{p^{k}}}$ 
 the assembly map
\begin{equation}\label{wervwoerijvewoivwev} \colim_{BC_{p^{k}}} \K(\widehat \Res^{C_{p^{k}}}(A))\to   \K(A\rtimes_{r}C_{p^{k}})
\end{equation}   
 becomes    an   equivalence after $p$-completion.
\end{theorem}

\begin{ex}
If $\K(A)\simeq 0$, then \cref{erijgowergrefwrefr}  implies that the $p$-completion of  $\K(A\rtimes_{r} C_{p})$ vanishes. By \cref{erjigowergfregfw}.\ref{wekgopwegrfwerfrf},  this is equivalent to the assertion that $p$ acts on $\K(A\rtimes_{r} C_{p})$ as an equivalence, or that  the groups $\K_{*}(A\rtimes_{r} C_{p})$ are uniquely   $p$-divisible. This fact has been observed in  \cite[Lemma 4.4]{Izumi_2004} with a different proof.
In    \cite{Izumi_2004},  \cite{MR3572256} the authors describe   many examples of actions of $C_{p}$ on a $C^{*}$-algebra $A$ 
with $K(A)\simeq 0$  such that $K_{*}(A\rtimes_{r} C_{p})$ are non-trivial uniquely $p$-divisible groups.
This   shows that the map in  \eqref{wervwoerijvewoivwev} is in general   not   an equivalence before $p$-completion.
\hB
  \end{ex}

\begin{ex}\label{jerotgrgertgertgfb}

By  \cref{iorgoergfgsfg}, for any $B$ in $\KK$, we can form the multiplicative induction $A:=B^{\otimes C_{p}}$ in $\KK^{C_{p}}$. By \cref{erijgowergrefwrefr}, the assembly map induces an equivalence 
\begin{equation}\label{efdewdewdewdq}L_{p}(\colim_{BC_{p}} \K(\widehat \Res^{C_{p}}(A))) \to L_{p}(\K(A\rtimes_{r}C_{p}))\ .\end{equation}   If $B$ is in the bootstrap class $\UCT$, then we have an equivalence $$\K(\widehat \Res^{C_{p}}(A))\simeq \K(B)^{\otimes C_{p}}$$ in $\Fun(BC_{p},\Mod_{KU}(\Sp))$,
where the right-hand side has the induced action by permutation of tensor factors. 
So the theorem reduces the calculation   of the $p$-completion of  $K(A\rtimes_{r}C_{p}) $  to the homotopy theoretic problem of calculating the $p$-completion of $\colim_{BC_{p}}  \K(B)^{\otimes C_{p}}$. In general, to calculate the homotopy groups of this colimit is still  a complicated task.

The following is the counterpart of the example discussed in \cref{iorjfoqrggegerfrefwref1}, but now in the case that $p|\ell$.
We consider $B:=\beins_{\KK}/\ell$ which belongs to the bootstrap class $\UCT$. A concrete realization of $A$ is then given in terms of the Cuntz algebra by 
$$A\simeq \kk^{C_{p}}(\underbrace{\cO_{\ell+1}\otimes \dots\otimes \cO_{\ell+1}}_{C_{p}})\ .$$
In this case the    homotopy groups of $\K( A \rtimes_{r}C_{p}) $ are known for $p=2$ by \cite{Izumi} 
$$\K_{*}( (\beins_{\KK}/2)^{\otimes C_{2}} \rtimes_{r}C_{2}) \simeq \left\{\begin{array}{cc} \Z/2\ell \Z\oplus \Z/\frac{\ell}{2}\Z &*=0\\0 &*=1  \end{array} \right.$$
and
$p=3$ by Kranz and Nishikawa (yet unpublished).
 As far as we know,  the case $p\ge 5$ is still an open problem. 
 
 One could try to use the spectral sequence associated to the tower obtained by applying
 $\K(-\rtimes_{r}C_{p})$ to the tower \eqref{ewqfqwoijoiqwedewdq} for the fibre sequence
 $\beins_{\KK}\to \beins_{\KK}/p\to \Sigma \beins_{\KK}$ in place of \eqref{fqewifhqiofsffeadf}. But in contrast to the case $p\not| \ell$ considered in \cref{iorjfoqrggegerfrefwref1},
this spectral sequence may have non-trivial higher differentials and is therefore  complicated to analyse.
    \hB
\end{ex}

 
 We finally discuss conditions on $A$ in $\KK^{C_{p^{k}}}$r implying  that  the assembly map \eqref{wervwoerijvewoivwev} is an equivalence even before $p$-completion. We concentrate on the case $k=1$.
In view of \cref{wekopgwerfrewfw}, the optimal condition is of course that $\K^{C_{p}}(A)[\xi^{-1}]\simeq 0$.  But  since this condition already involves the object $\K(A\rtimes_{r}C_{p})$ that we want to calculate  it does not seem to be useful in practice in general, but see \cref{ewfoijqwofdqewdqewdqewdq}.
 In order to get a more explicit condition we let $\bG(p)$ denote the localizing subcategory of $\KK^{C_{p}}$ generated by objects of the following form:
\begin{enumerate}
\item $\Res_{C_{p}}(A')$ for some $A'$ in $\KK$ such that $\K_{*}(A')$  is annihilated by $p^{n}$ for some $n$ in $\nat$.
\item  $\Ind^{C_{p}}(A'')$ for some $A''$ in $\KK$.
\end{enumerate}

\begin{prop}[\cref{wetghkopwergwerfw1}]\label{wetghkopwergwerfw}
If $A$ is in $\bG(p)$, then
the assembly map is an equivalence
\begin{equation}\label{efdewdewdewdq1nn} \colim_{BC_{p}}  \K(A)  \stackrel{\simeq}{\to}   \K( A \rtimes_{r}C_{p}) 
\end{equation} 
\end{prop}

We consider $B:=\beins_{\KK}/p^{n}$ and $A:=B^{\otimes C_{p}}$. It is not clear that $A$ belongs to $\bG(p)$. Nevertheless we have the following  result stating that here is no need of $p$-completion in order to calculate this crossed product.
\begin{prop}[\cref{ewfoijqwofdqewdqewdqewdq1}]\label{ewfoijqwofdqewdqewdqewdq}
 $$\colim_{BC_{p}} \K( \widehat \Res^{C_{p}}((\beins_{\KK}/p^{n})^{\otimes C_{p}}))\simeq \K((\beins_{\KK}/p^{n})^{\otimes C_{p}}\rtimes_{r}C_{p})\ .$$
\end{prop}

 
%

 {\em Acknowledgement: The author was supported by the SFB 1085 (Higher Invariants) funded by the Deutsche Forschungsgemeinschaft (DFG). The author thanks  B. Dünzinger, J. Kranz, M. Land, S.  Nishikawa, U. Pennig and Ch. Winges   for discussions about parts of this project.}



%
%
%

 \section{Some homotopy theory}

 \subsection{Equivariant KK-theory}\label{wtkopgwfrefrefw}

Equivariant group-valued $KK$-theory was introduced in \cite{kasparovinvent}. Following \cite{higsondiss}, \cite{Thomsen-equivariant}, \cite{Meyer:aa}, it should be considered as an additive category receiving
a functor from $G$-$C^{*}$-algebras that can be characterized by a universal property. By \cite{MR2193334}, the additive equivariant $KK$-category has a canonical structure of a triangulated category.
In order to perform homotopy theoretic constructions, e.g., homotopy colimits over group actions on objects in this category, we will work with the $\infty$-categorical enhancement of equivariant $KK$-theory. It is given by a functor from $G$-$C^{*}$-algebras to 
a stable $\infty$-category that is again characterized a universal property. 
The classical picture can be recovered by going over to the homotopy category.
 In this subsection we recall the main features of this $\infty$-categorical version of
 equivariant $KK$-theory. In the case of a discrete group $G$ (which is the only relevant case for the applications in the present paper), a reference   is \cite{KKG}.
\footnote{The general case of  second countable locally compact groups was the topic of a course taught in Regensburg 2023 whose notes are in preparation.}

 Let $G$ be a second countable locally compact group. We then consider the category of $G$-$C^{*}$-algebras $G\nCalg$.
 Its objects are $C^{*}$-algebras with an action of  $G$ by automorphisms such that, for all $a$ in $A$,
the  maps $G\to A$, $g\mapsto ga$, are
 continuous.  The morphisms of $G\nCalg$ are equivariant homomorphisms.

In the following, we will use the language of $\infty$-categories as developed in  \cite{htt}
\cite{HA}.  There exists a presentable stable  $\infty$-category
$\KK^{G}$ and a functor 
\begin{equation}\label{vwervhewiercecsffc}\kk^{G}:G\nCalg\to \KK^{G}
\end{equation} 
that
   is characterized by the following universal property:  For any cocomplete stable $\infty$-category $\cC$
  precomposition by $\kk^{G}$ induces an equivalence 
\begin{equation}\label{sdfbdfbvsfdvfdvsfvsfdvsfdvsf}\Fun^{\colim}(\KK^{G},\cC)\xrightarrow{\simeq} \Fun^{h,K_{G},\mathrm{se},\mathrm{sfin}}(G\nCalg,\cC)
\end{equation} 
of functor categories.
 Here the superscripts indicate full subcategories of functors satisfying corresponding conditions. The superscript $\colim$ stands for colimit-preserving functors. The superscripts
 $h$, $K_{G}$, $\mathrm{se}$ and $\mathrm{sfin}$ stand  for the conditions of  homotopy invariance, $K_{G}$-stability, semi-exactness and $s$-finitariness that will be explained in the following. Consider a functor $F$ on $G\nCalg$:
 \begin{enumerate}
 \item  It is homotopy invariant if it sends the canonical homomorphisms $A\to C([0,1],A)$, $a\mapsto \const_{a}$, to equivalences for all $A$ in $G\nCalg$.
 \item It is $K_{G}$-stable if it sends $K_{G}$-equivalences to equivalences. Here, a morphism $f:A\to B$ is a $K_{G}$-equivalence if $f\otimes \id_{G}:A\otimes K_{G}\to B\otimes K_{G}$ is a homotopy equivalence, where $K_{G}$ is the algebra of compact operators on  the $G$-Hilbert space  $L^{2}(G)\otimes L^{2}(\nat)$ with $G$-action  given by the  left-regular representation on $L^{2}(G)$.   \item It is semiexact if it sends semisplit exact sequences $0\to A\to B\to C\to 0$ in $G\nCalg$ to fibre sequences. The condition of being semisplit means that the map $B\to C$ admits  an equivariant completely positive contractive right inverse.
 \item \label{irwegjoergrffwre} It is $s$-finitary if the canonical morphisms $\colim_{A'\subseteq_{\sepa}A } F(A')\to F(A)$
 are  equivalences for all $A$ in $ G\nCalg$. Here, $A'\subseteq_{\sepa}A$ is the poset of $G$-invariant separable subalgebras
 of $A$. \end{enumerate}
 The existence of the functor \eqref{vwervhewiercecsffc} is shown by  providing a construction
 whose main idea is to force the desired universal properties using Dwyer-Kan localizations and $\Ind$-completion.

 The  minimal tensor product equips 
 the category $G\nCalg$ with a symmetric monoidal structure $\otimes$. We will denote the tensor unit by $\beins_{\KK^{G}}$. Using the symmetric monoidal version of Dwyer-Kan localizations \cite{hinich} and $\Ind$-completion,
 the category $\KK^{G}$ acquires a presentably symmetric monoidal structure
 and  functor $\kk^{G}$ has a symmetric monoidal refinement characterized by the following universal property: For every presentably symmetric monoidal $\infty$-category $\cC$, the precomposition by $\kk^{G}$ induces an equivalence   \begin{equation}\label{wfrrwefwrfrefvfdsv}
\Fun_{\otimes/\mathrm{lax}}^{\colim}(\KK^{G},\cC)\xrightarrow{\simeq} \Fun_{\otimes/\mathrm{lax}}^{h,K_{G},\mathrm{se},\mathrm{sfin}}(G\nCalg,\cC)\ ,
 \end{equation}
 where the subscripts $\otimes$ or $\mathrm{lax}$ stand for symmetric monoidal or lax symmetric monoidal functors, respectively. 
 There is an analogous statement for the maximal tensor product on $G\nCalg$.

 By the universal property of $\kk^{G}$, we get  factorizations
\begin{equation}\label{vwfsfjkvopqrefddcs}\xymatrix{G\nCalg\ar[r]^{L_{h}} \ar[dr]_{\kk^{G}}& G\nCalg_{h}\ar[d]^{\kk_{h}^{G}}\ar[r]^{L_{K_{G}}}&L_{K_{G}}\nCalg_{h}\ar[dl]^{\kk^{G}_{h,K_{G}}}\\&\KK^{G}&}\ ,
\end{equation} 
where $L_{h}$ and $L_{K_{G}}$ are the Dwyer-Kan localizations inverting homotopy equivalences and $K_{G}$-equivalences. 

\begin{rem}
Note that in the $KK$-theory functors considered in the present paper do not preserve countable coproducts. 
One could force this property by a further localization; see \cite[Rem. 3.4]{KKG}. \hB
\end{rem}

A main feature of equivariant $KK$-theory is the presence of a variety of group-change functors and various adjunctions between them. On the level of triangulated categories, they are  reviewed in \cite{MR2193334}. 
In the following we list their refinements to the $\infty$-categorical level. Thereby, we will provide sufficient data to fix the functors and adjunctions essentially uniquely.

If $H\to G$  is a homomorphism, then we have a restriction functor \begin{equation}\label{vsoijo3rfeffsvfvs}\Res^{G}_{H}:G\nCalg \to H\nCalg
\end{equation} 
that restricts the $G$-action to a $H$-action along this homomorpism.

Furthermore, for  a  closed  subgroup $H$ of $G$ have an  induction functor \begin{equation}\label{asdcadscaacd}  \Ind_{H}^{G}:H\nCalg\to G\nCalg\ .
\end{equation}    It sends $A$ on $H\nCalg$ to  the closure  in  $C_{b}(G,A)$ (the bounded  continuous functions on $G$ with values in $A$ with the supremum norm)  of  the subalgebra of functions $f$
 which are $H$-invariant in the sense that $f(gh)=hf(g)$  for all $h$ in $H$, and whose support projects to a compact subset of $G/H$. The group $G$ acts on this subalgebra  by the left-regular representation
 $(g',f)\mapsto (g\mapsto f(g^{\prime,-1}g))$.
 
 Both functors descend to equivariant $KK$-theory, i.e., we have canonical commutative squares
$$\xymatrix{G\nCalg\ar[r]^{\Res^{G}_{H}}\ar[d]^{\kk^{G}} & H\nCalg\ar[d]^{\kk^{H}} \\ \KK^{G}\ar@{..>}[r]^{\Res_{H}^{G}} &\KK^{H} } \qquad \xymatrix{H\nCalg\ar[r]^{\Ind^{G}_{H}}\ar[d]^{\kk^{H}} &G\nCalg \ar[d]^{\kk^{G}} \\\KK^{H} \ar@{..>}[r]^{\Ind_{H}^{G}} & \KK^{G}} \  .$$


We have a symmetric monoidal   functor 
$$\widehat \Res^{G}: G\nCalg \to \Fun(BG,\nCalg_{h})$$
that interprets   $G$-$C^{*}$-algebra   as an object of $G\nCalg_{h}$   with the induced $G$-action. 
This functor also descends to $KK$-theory in the sense that we have a canonical commutative square \begin{equation}\label{qfqwefqwedqwdq}\xymatrix{G\nCalg\ar[r]^-{\widehat \Res^{G}}\ar[d]^{\kk^{G}} &  \Fun(BG,\nCalg_{h})\ar[d]^{\kk_{h}} \\\KK^{G} \ar@{..>}[r]^-{\widehat \Res^{G}}&\Fun(BG, \KK) } \ .
\end{equation} 

Let  $UM(A)$ denote the unitary group of the multiplier algebra of $A$ equipped with the strict topology. 
\begin{lem} \label{wreojgkpwegwrefwrefwref}If
  the $G$-action on $A$ is induced by conjugation with a   continuous homomorphism $\rho: G\to UM(A)$,
 then we have an equivalence \begin{equation}\label{fvdsvervsdfvfdv}\kk^{G}(A)\simeq \Res_{G}\Res^{G}(\kk^{G}(A))
\end{equation} in $\KK^{G}$.  \end{lem}
\begin{proof}The equivalence in \eqref{fvdsvervsdfvfdv} is induced by the cocycle isomorphism
$$(\id_{A},\rho):\Res_{G}(\Res^{G}(A))\to A\ .$$ This is not simply a morphism in $G\nCalg$.
However,   cocycle morphisms can naturally be interpreted as morphisms in $L_{K_{G}}\nCalg_{h}$ and therefore also give rise to morphisms in $\KK^{G}$ via the factorization \eqref{vwfsfjkvopqrefddcs}.
 \end{proof}
 
 \begin{rem}
 Under the assumptions of \cref{wreojgkpwegwrefwrefwref},
we can conclude that the $G$-action on $\kk(A)$ in $\KK$ that is induced by functoriality
 is trivial 
  since this $G$-object can be written formally as $$\widehat \Res^{G}( \kk^{G}(A))\simeq \widehat \Res^{G}\Res_{G}\Res^{G}(\kk^{G}(A))\ ,$$ and the right-hand side obviously has the trivial $G$-action.
 This refines the classical fact that the action of a group on a $C^{*}$-algebra by inner automorphisms induces
 the trivial action on $K$-theory groups. \hB
\end{rem}
 
If $H$ is closed and open in $G$, then we have an adjunction
\begin{equation}\label{dfbvdfvdfvdfvsdfv}  \Ind^{G}_{H}:\KK^{H}\leftrightarrows \KK^{G}:\Res^{G}_{H}\end{equation}  
whose unit $A\mapsto \Res^{G}_{H}(\Ind^{G}_{H}(A))$ is already defined on the level of algebras and   sends $a$  in $A$ to the function
$$g \mapsto \left\{\begin{array}{cc} g^{-1}a&g\in H\\ 0&g\not\in H  \end{array} \right.\ .$$ 
This function is indeed continuous by the openess assumption on $H$.
If $G/H$ is in addition finite, then we also have an adjunction
 \begin{equation}\label{dfbvdfvdfvdfvsdfv1}  \Res^{G}_{H}:\KK^{G}\leftrightarrows \KK^{H}:\Ind^{G}_{H}\end{equation}  
 whose unit $A\to \Ind_{H}^{G}(\Res^{G}_{H}(A))$ and counit   
 $\Res^{G}_{H}(\Ind^{G}_{H}(A))\to A$ are already defined on the level of algebras and are given by 
the maps $a\mapsto (g\mapsto g^{-1}a)$ and the  
 evaluation $f\mapsto f(e)$, respectively.
 
 For $A$ in $\KK^{H}$ and $B$ in $\KK^{G}$, we   have a bi-natural projection formula
 \begin{equation}\label{dfbvdfvdfvdfvsdfv2} \Ind_{H}^{G}(A)\otimes B\simeq \Ind_{H}^{G}(A\otimes \Res^{G}_{H}(B))    \end{equation} 
 that follows from  the analoguous isomorphism 
 $f\otimes b\mapsto  (g\mapsto f(g)\otimes g^{-1} b)$ on the level of $C^{*}$-algebras.

 For $A$ in $G\nCalg$, we can form the convolution $*$-algebra $C_{c}(G,A)$. It is a pre-$C^{*}$-algebra  in the sense that all elements have a finite maximal norm. The completion of $C_{c}(G,A)$ is the maximal crossed product $A\rtimes_{\max}G$. If $A\to B(H)$ is a  faithful representation of $A$ on some Hilbert space $H$, then
the representation of $A\rtimes_{\max}G$
 on $L^{2}(G,H)$  induces the reduced norm on $A\rtimes_{\max}G$. The completion of $A\rtimes_{\max}G$ with respect to the reduced norm is called  the reduced crossed product $A\rtimes_{r}G$. 
  
The crossed product functors
$$-\rtimes_{r}G, -\rtimes_{\max}G:G\nCalg\to \nCalg$$ both descends to equivariant $KK$-theory, i.e., we have canonical commutative squares
$$\xymatrix{G\nCalg\ar[r]^{-\rtimes_{r}G}\ar[d]^{\kk^{G}}&  \nCalg\ar[d]^{\kk} \\ \KK^{G}\ar@{..>}[r]^{-\rtimes_{r}G} & \KK} \ , \qquad \xymatrix{G\nCalg\ar[r]^{-\rtimes_{\max}G}\ar[d]^{\kk^{G}}&  \nCalg\ar[d]^{\kk} \\ \KK^{G}\ar@{..>}[r]^{-\rtimes_{\max}G} & \KK}  \ .$$
 For a closed subgroup $H$ of $G$, we have  equivalences of functors (Green's imprimitivity theorem)
\begin{equation}\label{bwefljrfopwerfrerfwerfwreferf}-\rtimes_{r}H\simeq \Ind_{H}^{G}(-)\rtimes_{r}G:\KK^{H}\to \KK\ , \quad -\rtimes_{\max}H\simeq \Ind_{H}^{G}(-)\rtimes_{\max}G:\KK^{H}\to \KK
\end{equation}

If $G$ is discrete, then this equivalence is induced by the composition of  $-\rtimes_{?}H$ applied to the unit of the $(\Ind^{G}_{H},\Res^{G}_{H})$-adjunction and the  transformation $\Res^{G}_{H}(-)\rtimes_{?}H\to -\rtimes_{?}G$ induced on the level of convolution algebras by the extension by zero $C_{c}(H,\Res^{G}_{H}(-))\to C_{c}(G,-)$ (it is here where we use discreteness).  In the general  case, the unit is induced by a   Morita bimodule and therefore starts to exist in $L_{K_{G}}G\nCalg_{h}$. 

If $G$ is compact, or more generally, amenable, then the canonical  transformation
$$-\rtimes_{\max}G\to -\rtimes_{r}G$$ is an equivalence.

  An
 instance of Green's imprimitivity theorem for $H=\{e\}$  is the equivalence  $ \beins_{\KK}\simeq \kk^{G}(C_{0}(G))\rtimes_{r}G$, 
 where we use that  $\kk^{G} (C_{0}(G))\simeq \Ind^{G}(\beins_{\KK})$. 
Here, we consider $C_{0}(G)$ in $G\nCalg$ using the $G$-action by left translations.  
The object $C_{0}(G)$ in $G\nCalg$ has an additional action of $G$ induced by right translations.
 We thus get a refined object $\kk^{G} (C_{0}(G))$ in $\Fun(BG,\KK^{G})$ and therefore
$\kk^{G} (C_{0}(G))\rtimes_{r}G$ in $\Fun(BG,\KK)$. 
Let $\triv^{G}:\KK\to \Fun(BG,\KK)$ be the functor which  equips an object with the trivial $G$-action.
\begin{lem} \label{wtkorhopweferwfwregw} The Green's imprimitivity equivalence refines to an equivalence 
  \begin{equation}\label{qewfiojioefqwdwedewqdedqd}\triv^{G}(\beins_{\KK})\simeq \kk^{G}(C_{0}(G))\rtimes_{r}G
\end{equation}
in $\Fun(BG,\KK)$. \end{lem}
 \begin{proof}
We consider the  Hilbert space $L^{2}(G)$. The algebra $C_{0}(G)$ acts
by multiplication operators. The action is covariant for the action of  $G$ on
the Hilbert space and on the algebra by left translations. By the universal property of the reduced crossed product, we therefore get a homomorphism
\begin{equation}\label{fwewdqededqd}C_{0}(G)\rtimes_{r}G\to K(L^{2}(G))
\end{equation}  which turns out to be an isomorphism. 

We let $\rho$ denote the action on $G$ on $L^{2}(G)$ by right-multiplication. 
We interpret $\rho$ as a continuous  homomorphism $G\to UM(K(L^{2}(G)))$.
We then observe that, under the isomorphism \eqref{fwewdqededqd},  the $G$-action on  
$C_{0}(G)\rtimes_{r}G$ induced by right-translations on $C_{0}(G)$ goes to the action on $ K(L^{2}(G))$ induced by conjugation with $\rho$.
We view $K(L^{2}(G))$ with this conjugation action as an object of $G\nCalg$. 
We then have
$$\kk^{G}(K(L^{2}(G))) \stackrel{\eqref{fvdsvervsdfvfdv}}{\simeq} \Res_{G}(\Res^{G}(\kk^{G}(K(L^{2}(G)))))\ .$$
By stability, the right-hand side is equivalent to $\Res_{G}(\beins_{\KK})$. 
This implies that the object $\widehat \Res^{G} \kk^{G}(K(L^{2}(G)))$ in $\Fun(BG,\KK)$ is equivalent to  the object obtained by equipping $\beins_{\KK}$ with the trivial $G$-action.
%
%
%
%
  \end{proof}

 For $B$ in $\KK$, we  have an equivalence of functors (also natural in $B$)
\begin{equation}\label{fwerfeiuhieuworf}B\otimes (-\rtimes_{r}G) \simeq (\Res_{G}(B)\otimes -)\rtimes_{r}G :\KK^{G}\to \KK 
\end{equation} 
 that comes from a similar natural isomorphism 
   on the level of $C^{*}$-algebras. 
  Here, it is relevant that we work with the minimal tensor products.
  There is an analogous result for the maximal crossed produce involving maximal tensor products.
  

 If $G$ is compact, then we have an adjunction (the Green-Julg theorem)
\begin{equation}\label{bsdfvkspdfvqreve}\Res_{G}:\KK\leftrightarrows \KK^{G}:-\rtimes_{r}G
\end{equation} 
whose unit $A\to \Res_{G}(A)\rtimes_{r}G$  is already defined on the level of algebras  and given by $a\mapsto \const_{a}$.

Finally, if $G$ is discrete, then we have an adjunction (the dual Green-Julg theorem)
\begin{equation}\label{verojvpevfrefwfref}-\rtimes_{\max}G:\KK^{G}\leftrightarrows  \KK:\Res_{G}
\end{equation} 
whose unit $A\mapsto \Res_{G}(A\rtimes_{\max}G)$ is given by a morphism which already exists in  $L_{K_{G}}G\nCalg_{h}$.  It is given by the cocycle morphism $A\ni a\mapsto  a\delta_{e} \in C_{c}(G,A)$.

For both adjunctions,   the counits  also exist in  $L_{K_{G}}G\nCalg_{h}$, and the adjunctions
are descended from this intermediate localization, see \eqref{vwfsfjkvopqrefddcs}

Let $\beins:=\kk(\C)$ denote the tensor unit of $\KK$. 
The commutative algebra $\map_{\KK}(\beins,\beins)$ in $\CAlg(\Sp)$ is the  complex $K$-theory spectrum from homotopy theory that is   usually denoted by $KU$. As a consequence, the category
$\KK$ is a $\Mod_{KU}(\Sp)$-module in the category $\Pr^{L}_{\mathrm{st}}$ of stable presentable $\infty$-categories. In particular, its mapping spectra are $KU$-modules   in a canonical way.

The symmetric monoidal restriction functor $\Res_{G}:\KK
\to \KK^{G}$  turns $\KK^{G}$ into a $\KK$- and hence into a $\Mod_{KU}(\Sp)$-module and   we have 
  $\beins_{\KK^{G}}  \simeq \Res_{G}(\beins)$.   
  \begin{ddd}
We define   the representation ring  $$R(G):=\map_{\KK^{G}}(\beins_{\KK^{G}} ,\beins_{\KK^{G}} )$$ of $G$  in $\CAlg(\Mod_{KU}(\Sp))$. \end{ddd} 
The category $\KK^{G}$ is   naturally a module over $\Mod_{R(G)}(\Sp)$ in $\Pr^{L}_{\mathrm{st}}$, and
its mapping spectra are naturally $R(G)$-modules.

\begin{rem}\label{wrjtohipwgwregfwerfw}
 If $G$ is compact, then  $\pi_{0}R(G)$ can be identified with the classical representation ring of $G$.  In particluar, for any   finite-dimensional representation $V$ of $G$  we get an element $[V]$ in $\pi_{0}R(G)$.
 \hB
\end{rem}
\begin{ex}\label{wekogpwefgrefwerfwefew}
 If $H$ is open in $G$ and $G/H$ is finite, then the functor $\Res^{G}_{H}:\KK^{G}\to \KK^{H}$ has  a left and a right adjoint that  are both given by the functor $\Ind_{H}^{G}:\KK^{H}\to \KK^{G}$ and the units described above.
 The identity of $\Ind_{H}^{G}$ identifies the left and the right adjoint and refines
 $ \Res^{G}_{H}$ to an iso-normed functor $(\Res^{G}_{H},\Nm^{G}_{H})$ in the sense of  \cite[Def. 2.1.3]{zbMATH07526076}. Since $\Res^{G}_{H}$ is symmetric monoidal and we have the projection formula
 \eqref{dfbvdfvdfvdfvsdfv2}, the normed functor is actually $\otimes$-normed in the sense of 
 \cite[Def. 2.3.1]{zbMATH07526076}.
 
  By  \cite[Def. 2.1.11]{zbMATH07526076}, for every two objects $A,B$ in $\KK^{G}$, we have an integration
 $$\int_{\Res^{G}_{H}}:\KK^{H}(\Res^{G}_{H}(A),\Res^{G}_{H}(B))\to \KK^{G}(A,B)$$ that sends
$f:\Res^{G}_{H}(A)\to \Res^{G}_{H}(B)$ to the composition
\begin{equation}\label{feqwdewdadc}A\to \Ind^{G}_{H}(\Res^{G}_{H}(A)) \stackrel{\Nm^{G,-1}_{H}}{\simeq}\Ind^{G}_{H}(\Res^{G}_{H}(A))
\xrightarrow{\Ind^{G}_{H}(f)} \Ind^{G}_{H}(\Res^{G}_{H}(B))\to B\ ,
\end{equation} 
where the first map is the unit of the $(\Res^{G}_{H},\Ind^{G}_{H})$-adjunction, and the last map 
is the counit of the  $(\Ind^{G}_{H},\Res^{G}_{H})$-adjunction. 
We want to understand this map explicitly in the case where $f$ comes from a $C^{*}$-algebra homomorphism.

So changing notation, let $A,B$ be in $G\nCalg$ and let $f:\Res^{G}_{H}(A)\to \Res^{G}_{H}(B)$ be an $H$-equivariant  homomorphism.
Then we want to calculate $\int_{\Res^{G}_{H}}\kk^{H}(f)$.
The first three maps in the composition \eqref{feqwdewdadc} are defined already on the level 
of $C^{*}$-algebras and send $f$ to the  function 
 $g \mapsto  f(g^{-1}a) $.
 The counit of the $(\Ind^{G}_{H},\Res^{G}_{H})$-adjunction
 only exists  in the localization $L_{K_{G}}G\nCalg_{h}$ explained in  \eqref{vwfsfjkvopqrefddcs}.  
  It can be described as follows.  A special case of the 
 projection formula \eqref{dfbvdfvdfvdfvsdfv2} provides an isomorphism 
$$\Ind_{H}^{G}(\Res^{G}_{H}(B))\cong  C(G/H)\otimes B$$ in $G\nCalg$.
We furthermore embed $C(G/H)\otimes B$ into $K(L^{2}(G/H))\otimes B$,
where the first factor acts by multiplication operators.
The counit is then represented by the  sequence of homomorphisms of $G$-$C^{*}$-algebras
\begin{equation}\label{sadvadsacaewcc}\Ind_{H}^{G}(\Res^{G}_{H}(B))\cong  C(G/H)\otimes B\to K(L^{2}(G/H))\otimes B\xleftarrow{\simeq} B\ ,
\end{equation} 
where the right map is a corner embedding sending $b$ in $B$ to the operator $p\otimes b$, where
$p$ in $ K(L^{2}(G/H))$ is the projection to the constant functions on $G/H$.
This map is a $K_{G}$-equivalence and therefore invertible in $L_{K_{G}}G\nCalg_{h}$. 
The integral of $f$ is then represented (up to the last equivalence in \eqref{sadvadsacaewcc}) by the  homomorphism of $G$-$C^{*}$-algebras
$A\to  K(L^{2}(G/H))\otimes B$ given by 
$$A\ni a\mapsto \sum_{gH\in G/H}  p_{gH}\otimes  gf(g^{-1}a) \in K(L^{2}(G/H))\otimes B\ ,$$
where $p_{gH}$ in $ K(L^{2}(G/H))$ is the projection onto the one-dimensional subspace of functions supported on $gH$.

If $f$ is $G$-invariant, then we read off that 
$$\int_{\Res_{H}^{G}} \Res_{H}^{G}(f):A\xrightarrow{a\mapsto 1_{K(L^{2}(G/H))}\otimes f(a)} K(L^{2}(G/H))\otimes B\stackrel{\simeq}{\leftarrow} B \ .$$
The zig-zag of homomorphisms $$\C\xrightarrow{1\mapsto 1_{K(L^{2}(G/H))} }K(L^{2}(G/H))\stackrel{\simeq}{\leftarrow} \C$$
represents the class $\rho_{G/H}$ in the representation ring $\pi_{0} R(G)$ given by the representation of $G$ on $L^{2}(G/H)$.  So $\int_{\Res_{H}^{G}} \Res_{H}^{G}(f)\simeq \rho_{G/H} f$.
In particular, using the notation of  \cite[Def. 2.13]{zbMATH07526076}, we have 
 $$\int_{\Res^{G}_{H}} \id_{\Res^{G}_{H}}(A)\simeq |\Res^{G}_{H}|_{A}\simeq \rho_{G/H}\ .$$
\hB
 \end{ex}

For certain applications, the equivariant $KK$-theory functor  \eqref{vwervhewiercecsffc} has two undesirable problems: 
It does not  send all infinite (say countable) sums of $C^{*}$-algebras to sums in $\KK$, and it does not send every exact sequence
of $C^{*}$-algebras to a fibre sequence. 
One can improve on both by modifying the construction of $\KK^{G}$. The equivariant $E$-theory functor
 \begin{equation}\label{ffqewfewffqewf}\ee^{G}:G\nCalg\to \EE^{G}\end{equation}  is the universal homotopy invariant, $K_{G}$-stable, exact, s-finitary and countable sum preserving functor to a cocomplete stable $\infty$-category. 
Here exact and  countable sum preserving mean
%
 that the functor $\ee^{G}$ sends all (as opposed to all semiexact in the case of $\KK^{G}$) exact sequences of $G$-$C^{*}$-algebras to fibre sequences and that $$ \bigoplus_{i\in I} \ee^{G}(A_{i})\to \ee^{G} (\bigoplus_{i\in I}A_{i})$$ is an equivalence for every
countable  family $(A_{i})_{i\in I}$ in $G\nCalg$.
 Thus if $\cC$
is a cocomplete stable $\infty$-category, then pull-back along 
$\ee^{G}$ induces an equivalence
\begin{equation}\label{bvwervopkpsdfvdfvsfdv}\Fun^{\colim}(\EE^{G},\cC)\xrightarrow{\simeq} \Fun^{h,K_{G},\mathrm{ex},\oplus,\mathrm{sfin}}(G\nCalg,\cC)\ .\end{equation}
The functor $\ee^{G}$ has a symmetric monoidal  enhancement for the maximal tensor product in $G\nCalg$ such that $\EE^{G}$ has a presentably symmetric monoidal structure and  for every presentably symmetric monoidal    stable $\infty$-category $\cC$ the pull-back along $\ee^{G}$ induces an equivalence
\begin{equation}\label{sdfbdfbvsfdvfdvsfvwr3r4r4fewesfdvsfdvsf}\Fun_{\otimes/\mathrm{lax}}^{\colim}(\EE^{G},\cC)\xrightarrow{\simeq} \Fun_{\otimes/\mathrm{lax}}^{h,K_{G},\mathrm{ex},\oplus,\mathrm{sfin}}(G\nCalg,\cC)\ .\end{equation}
We refer to \cite[Sec. 3.2]{budu} for the details in the case of discrete groups. An additional nice feature of the equivariant $E$-theory functor $\ee^{G}$  is that it preserves all  filtered colimits. Indeed, its restriction to seperable algebras preserves countable filtered colimits by  \cite[Prop. 3.19]{budu} (see also  the text  before \cite[Lem. 2.7]{MR2193334}), and the general case follows from s-finitarines.

 For the following, we equip $\KK^{G}$ with the maximal tensor product structure.
 The universal property   \eqref{wfrrwefwrfrefvfdsv}  provides a
  canonical left-adjoint symmetric monoidal functor 
\begin{equation}\label{wrfqwdewdewfqef}c^{G}:\KK^{G}\to \EE^{G}\ .
\end{equation}   It equips $\EE^{G}$ with the structure of a commutative algebra over $\KK^{G}$ in $\Pr^{L}_{\mathrm{st}}$.
The functor \eqref{wrfqwdewdewfqef} is a Dwyer-Kan localization forcing exactness and countable sum preserving. Therefore all  the adjunctions and equivalences between functors  in $KK$-theory 
descend to corresponding adjunctions between their descends to $E$-theory
provided the functors descend. This applies to induction and restriction and the maximal crossed products,
 see also \cite[Sec. 3.6]{budu}.

The   $K$-theory functor $K:\nCalg\to \Mod_{KU}(\Sp)$ it known to be an  object of $$\Fun_{\mathrm{lax}}^{h,K,\exa,\mathrm{sfin},\oplus}(\nCalg,\Mod_{KU}(\Sp))\ .$$ By the universal property
\eqref{bvwervopkpsdfvdfvsfdv}, it has an essentially unique  factorization \begin{equation}\label{sadfvasvadscddascadsca}K:\nCalg\xrightarrow{\ee } \EE \xrightarrow{\K^{E}} \Mod_{KU}(\Sp)
\end{equation}
through a  lax symmetric monoidal colimit preserving functor $\K^{E}$.
One can then show that
$\K^{E}\simeq \map_{\EE }(\beins_{\EE },-)$
and that $\K\simeq \K^{E}\circ c$ with $\K$ as in \eqref{hertoijiojgretbth} for the trivial group $G$.
%
%
%

 \begin{prop}[\cref{jeirgowregwerfrefw}]\label{jeirgowregwerfrefw1}
 Assume that $G$ is discrete. If $A$ is in $\langle \Ind^{G}(\KK) \rangle$, then we have an equivalence
\begin{equation}\label{vferopvjopfdewdedq111}\colim_{BG} \K(\widehat \Res^{G}(A)) \simeq \K(A\rtimes_{r}G)\ .
\end{equation}  
 \end{prop}
\begin{proof} 

The idea of the proof is to define a natural transformation 
\eqref{frewfwrewoifjowrfwrefwerf} between colimit preserving functors on $\KK^{G}$ 
and to show that it induces an equivalence on the image of the induction  functor $\Ind^{G}$.
We have already seen by an explicit calculation  in \cref{wtkohpwerfrefwref} that on induced algebras  both sides of this transformation are equivalent.   The main technical problem that we will solve in this proof
is to construct the transformation  in such a way that one can see that it induces the above equivalence.

We first describe the construction of the assembly map \eqref{frewfwrewoifjowrfwrefwerf1} between $\KK$-valued functors. For any $G$-$C^{*}$-algebra $A$ we have a  natural isomorphism in $G\nCalg$
$$A\otimes C_{0}(G)\stackrel{ \cong}{\to} \Res_{G}(\Res^{G}(A)) \otimes C_{0}(G) \ .$$
Identifying elements of the tensor product with functions $f:G\to A$, this isomorphism is given  by
$f\mapsto (g\mapsto g^{-1}f(g))$. Under this isomorphism the $G$-action by right-translations on the $C_{0}(G)$-factor  on the codomain goes to the action $(h,f)\mapsto (g\mapsto h f(gh))$ on the domain. 
On the level of functors from $\KK^{G}$ to $\Fun(BG,\KK^{G})$,  this isomorphism induces an equivalence  \begin{equation}\label{bfsbspokvopsrwgfergvfrrrs} \triv^{G}(-)\otimes C_{0}(G)\simeq  \Res_{G}(  \widehat\Res^{G}(-))\otimes C_{0}(G)\ ,
\end{equation} where $\triv^{G}:\KK^{G}\to \Fun(BG,\KK^{G})$ equips an object with the trivial $G$-action.
 We now have the following chain of equivalences of functors from $\KK^{G}$ to $\Fun(BG,\KK)$: \begin{eqnarray} 
 (\triv^{G}(-)\otimes \kk^{G} (C_{0}(G)))\rtimes_{r}G&\stackrel{\eqref{bfsbspokvopsrwgfergvfrrrs}}{\simeq} &
(\Res_{G}(\widehat \Res^{G}(-))\otimes  \kk^{G}( C_{0}(G)))\rtimes_{r}G\nonumber\\&\stackrel{\eqref{fwerfeiuhieuworf}}{\simeq}&
\widehat \Res^{G}(-)\otimes  (\kk^{G}( C_{0}(G))\rtimes_{r}G)\nonumber\\
&\stackrel{\eqref{qewfiojioefqwdwedewqdedqd}}{\simeq}&
\widehat \Res^{G}(-)\otimes \triv^{G}(\beins_{\KK})\nonumber\\&\simeq&\widehat \Res^{G}(-) . \label{kofpqfdewdfeqwfeqwfed}
\end{eqnarray}

We have a map $*_{G}\C\to C_{0}(G)$ in $\Fun(BG,G\nCalg)$,
where we consider  the free $G$-indexed product $*_{G}\C$ of copies of $\C$ as a $G$-$C^{*}$-algebra for the left action on the index set with the additional 
action by right translations. The map is fixed such that
$\C\xrightarrow{\iota_{g}} *_{G}\C\to C_{0}(G)$ sends $1$ to $\delta_{g}$ for every $g$ in $G$.
The map $\kk^{G}(*_{G}\C)\to  \kk^{G}(C_{0}(G))$ is an equivalence in $\Fun(BG,\KK^{G})$. This follows, e.g., from \cite[Prop. 6.9]{KKG}
applied to the $G$-$C^{*}$-category $\C[G]$ obtained by linearizing the discrete $G$-category $G$.
We  therefore have an equivalence
\begin{equation}\label{vsdfvoijdfvovjiosdfvfdv}\widehat \Res^{G}(-)\stackrel{\eqref{kofpqfdewdfeqwfeqwfed}}{\simeq}  (\triv^{G}(-)\otimes  \kk^{G}( *_{G}\C))\rtimes_{r}G
\end{equation} 
of functors from $\KK^{G}$ to $\Fun(BG,\KK)$.  
We have a morphism \begin{equation}\label{ewfqefdewdqd}\sigma:*_{G}\C\to \triv^{G}(\Res_{G}(\C))
\end{equation}  in $\Fun(BG,G\nCalg)$ that is fixed such that
$\C\xrightarrow{\iota_{g}} *_{G}\C\to \C$ is the identity for all $g$ in $G$. It induces a morphism
$ \kk^{G}( *_{G}\C)\to  \triv^{G}(\beins_{\KK^{G}})$ in $\Fun(BG,\KK^{G})$ and hence, using some obvious identifications,  a natural transformation
$$\widehat \Res^{G}(-)\stackrel{\eqref{vsdfvoijdfvovjiosdfvfdv}}{\simeq}  (\triv^{G}(-)\otimes  \kk^{G}( *_{G}\C))\rtimes_{r}G
\to  \triv^{G}(- \rtimes_{r}G)$$ of functors from $\KK^{G}$ to $\Fun(BG,\KK)$.
Its adjoint is a natural transformation
\begin{equation}\label{frewfwrewoifjowrfwrefwerf1}\colim_{BG}\widehat \Res^{G}(-)\to (-\rtimes_{r}G):\KK^{G}\to \KK \ .
\end{equation} 
We must show that it induces an equivalence after  precomposition with $\Ind^{G}$.
Equivalently, we must show that the adjoint of
$$(\triv^{G}(\Ind^{G}(A))\otimes  \kk^{G}( *_{G}\C))\rtimes_{r}G\to \triv^{G}(\Ind^{G}(A)\rtimes_{r}G)$$ is an equivalence in $\KK$.
Since $-\rtimes_{r}$ preserves colimits it suffices to show that the adjoint of  \begin{equation}\label{fqwefiojwodqwedqewdqewd}
\triv^{G}(\Ind^{G}(A))\otimes  \kk^{G}( *_{G}\C) \to \triv^{G} (\Ind^{G}(A))
 \end{equation}
is an equivalence in $\KK^{G}$. 

For $A$ in $\nCalg$, the induction of $A$ is the $G$-$C^{*}$-algebra  $\Ind^{G}(A)\cong C_{0}(G,A)$
with the $G$-action $(g,f)\mapsto (h\mapsto f(g^{-1}h))$.
We have a natural isomorphism 
 \begin{eqnarray}
\triv^{G}  (C_{0}(G,A))\otimes  *_{G}\C&\cong&
C_{0}(G, A\otimes *_{G}\C)\nonumber\\&\stackrel{f\mapsto \tilde f}{\cong}&
C_{0}(G, A\otimes  *_{G}\C)\nonumber\\&\cong&
\triv^{G}  (C_{0}(G,A))\otimes \Res_{G}(\Res^{G}(*_{G}\C))\label{ejgkowergferfwefd}
\end{eqnarray}
in $\Fun(BG,G\nCalg)$.
Here $\tilde f(g):=g^{-1}f(g)$ for all $g$ in $G$.
The first copy of $C_{0}(G,  A\otimes *_{G}\C)$  is a $G$-$C^{*}$-algebra
with the action $(g,f)\mapsto (h\mapsto gf(g^{-1}h))$ and has the additional $G$-action $(l,f)\mapsto (h\mapsto f(hl))$. The second copy 
is a $G$-$C^{*}$-algebra
with the action $(g,f)\mapsto (h\mapsto f(g^{-1}h))$ and has the additional $G$-action $(l,f)\mapsto (h\mapsto f(hl))$. One checks that
$$\xymatrix{\triv^{G}  (C_{0}(G,A))\otimes  *_{G}\C\ar[rr]_{\eqref{ejgkowergferfwefd}}^{\cong}\ar[dr] &&\ar[dl] \triv^{G}  (C_{0}(G,A))\otimes \Res_{G}(\Res^{G}(*_{G}\C))\\&\triv^{G}  (C_{0}(G,A))&}$$
  commutes, where the diagonal arrows are both induced by $\sigma$ in \eqref{ewfqefdewdqd}.
 Applying $\kk^{G}$ we get a commutative diagram
$$\xymatrix{\triv^{G}  ( \Ind^{G}(A))\otimes  \kk^{G}(*_{G}\C)\ar[rr]^{\cong}\ar[dr]_{\eqref{fqwefiojwodqwedqewdqewd}} &&\ar[dl] \triv^{G}  ( \Ind^{G}(A))\otimes \Res_{G}(\Res^{G}(\kk^{G}(*_{G}\C)))\\&\triv^{G}( \Ind^{G}(A))&}$$
in $\Fun(BG,\KK^{G})$ for any $A$ in $\KK$.
We must show that the adjoint of the right diagonal arrow is an equivalence.
Since $\triv^{G}(\Ind^{G}(A))\otimes \Res_{G}(-)$ preserves colimits it suffices to show that the adjoint of 
\begin{equation}\label{porjfgroepfqfqewfqwdewd} \Res^{G}(\kk^{G}(*_{G}\C))\to \triv^{G}(\beins_{\KK})
\end{equation} is an equivalence in $\KK$.
We have an equivalence $$\Res^{G}(\kk^{G}(*_{G}\C))\simeq \kk(\Res^{G}(*_{G}\C)) \simeq \kk(\Res^{G}(C_{0}(G))\simeq \bigoplus_{G}\beins_{\KK}$$ in $\Fun(BG,\KK)$, where $G$ acts on the index set of the sum by right translations.
Under this identification,  the map  
 in \eqref{porjfgroepfqfqewfqwdewd} goes to the canonical map which induces the identity on each summand.
The adjoint of this map is well-known to be an equivalence.

Applying the colimit preserving functor  $\K:\KK\to \Mod_{KU}(\Sp)$ to the assembly map \eqref{frewfwrewoifjowrfwrefwerf1},
we get the assembly map 
\begin{equation}\label{frewfwrewoifjowrfwrefwerf}\colim_{BG}\K(\widehat \Res^{G}(-))\to \K(-\rtimes_{r}G):\KK^{G}\to \Mod(KU)
\end{equation} 
appearing in the statement of  \cref{jeirgowregwerfrefw}.
Above we have shown that
\eqref{frewfwrewoifjowrfwrefwerf} is an equivalence on the image of $\Ind^{G}:\KK\to \KK^{G}$.
Since it is a natural transformation between colimit preserving functors 
it is then also an equivalence on $\langle \Ind^{G}(\KK)\rangle$.  \end{proof}

\begin{rem}\label{wejkogpefrefewrfewrfewrfrewfw}
For torsion-free groups, the equivalence  asserted in \cref{jeirgowregwerfrefw} can be considered as a special case of a result of \cite{OO} or \cite{MR1836047} after identifying the left-hand side with the usual left-hand side of the Baum-Connes assembly map as in \eqref{fqweqdewdedeq}.

The proof  for  general  discrete groups could be also deduced from
the results of \cite{kranz} or \cite{bel-paschke}. One only needs  to observe  
that the functor $$ K^{G}_{-}:G\nCalg\to \Fun(G\Orb,\Mod_{KU}(\Sp))$$ mentioned in \cref{wergjiwoeferfw}
factorizes over a colimit preserving functor
$$\cK_{-}^{G}:\KK^{G}\to  \Fun(G\Orb,\Mod_{KU}(\Sp))$$  which has the following properties:
\begin{enumerate}
\item $\cK^{G}_{A}(G/H)\simeq A\rtimes_{r}H$ for all subgroups $H$ of $G$.
\item $\cK^{G}_{A}(G)\simeq \K(\widehat \Res^{G}(A))$, where $G$ acts on the left-hand side by right-translations on the $G$-orbit $G$.
\item\label{wrelkgpwerfwefwrefw}  For $B$ in $\KK$, the functor $\cK_{\Ind^{G}(B)}^{G}$ is the left Kan-extension along the inclusion  $BG\to G\Orb$ of
the constant functor on $BG$ with value $\K(B)$.
\end{enumerate}
The existence of such a factorization has been shown in \cite{kranz}, see also  \cite[Sec. 15]{bel-paschke}.
 Using these properties, it is a simple exercise with Kan extensions to deduce
 \cref{jeirgowregwerfrefw}.
 
  The main point of the proof given above
 is to  bypass the construction of the functor $\cK^{G}_{-}$ and, in particular, the verification of the Property \ref{wrelkgpwerfwefwrefw}, that in the references was quite involved. 
 See \cref{ijrogergferwg9} for an analogous result for finite groups $G$.
\hB
\end{rem}

Let $A$ be an object of $\KK^{G}$ and $\rho$ in $\pi_{0}R(G)$ denote the class of the left-regular representation of the finite group $G$.
\begin{prop}[\cref{qirjfofdewdewdqewd}]\label{qirjfofdewdewdqewd1} We assume that $G$ is finite. The following assertions are equivalent.
\begin{enumerate}
\item  \label{qrwijfoqwfqewdqwedew}$\rho$ acts as an equivalence on $A$.
\item\label{qrwijfoqwfqewdqwedew1} We have  $A\in  \langle \Ind^{G}(\KK) \rangle$ and $|G|$ acts an equivalence on $A$.
\end{enumerate}
\end{prop}
\begin{proof} 
We first show that \ref{qrwijfoqwfqewdqwedew}  implies \ref{qrwijfoqwfqewdqwedew1}.
 Since $G\times G$ decomposes as a $G$-set into a union of $|G|$ copies of $G$ we have  the relation $\rho^{2}\simeq |G|\rho$ in $\pi_{0}R(G)$. So if $\rho$ acts as an equivalence on $A$, then so does $|G|$. 
 We now use \cref{wekogpwefgrefwerfwefew} for $H=\{e\}$.
The composition \eqref{feqwdewdadc} for $f=\id_{A}$
provides a factorization of the multiplication map by $\rho$ through
$$A\to \Ind^{G}(\Res^{G}(A))\to A\ .$$
Since $\rho$ acts as an equivalence on $A$, we conclude that $A$ is a retract of $\Ind^{G}(\Res^{G}(A))$.

We now argue that  \ref{qrwijfoqwfqewdqwedew1}  implies \ref{qrwijfoqwfqewdqwedew} .
  Since the fact that $\rho$ acts as an equivalence is preserved under colimits it suffices to show that $\rho$ acts as an equivalence on $A:=\Ind^{G}(B)$ for any $B$ in $\KK$ such that $|G|$ acts as an equivalence on $A$.
By the projection formula \eqref{dfbvdfvdfvdfvsdfv2}  and using $\Res^{G}(\rho)=|G|$, we obtain a chain of equivalences of morphisms 
\begin{eqnarray*}
(\rho:A\to A)&\simeq& (\rho :\Ind^{G}(B)\to \Ind^{G}(B))\\&\simeq&
\Ind^{G}(\Res^{G}(\rho):B\to B)\\&\simeq& ( \Ind^{G}(|G|:B \to  B)) \\&\simeq
& (|G|:A\to A)
\end{eqnarray*}
We can conclude that $\rho:A\to A$ is an equivalence.
\end{proof}

  \subsection{Power functors}\label{wrtogpwferfwerf9}

 In the last subsection \cref{wtkopgwfrefrefw}, we described exact restriction and induction functors which can be used to describe objects in  $\KK^{G}$.  Since these functors fit into adjunctions,  we can understand the resulting objects well.  In the present subsection, we introduce  non-exact multiplicative induction  functors. They produce
 objects of $\KK^{G}$ that do not have an obvious presentation in terms of  exact constructions.

\begin{rem}The  multiplicative induction functors are part of a global symmetric monoidal structure on the global category $\underline{\KK}$ (see also \cref{jrgergoergweferfweff} below) that will be discussed in future work. The present subsection is designed to provide a reference for  all technical background results for this project that require to go  through the $C^{*}$-algebraic details of the construction of the stable $\infty$-categories $\KK^{G}$.  For this reason, it is slightly more general than needed for the purpose of the present paper.
\hB\end{rem}

In order encode the  naturality properties of the multiplicative induction  functors properly,   we will work with finite groupoids.  A finite groupoid $\cG$ is connected if every two objects are connected by a morphisms. Assume that $\cG$ is a connected finite groupoid. If we fix
an object $x_{0}$ in $\cG$ and let $G$ denote its stabilizer group, then the inclusion $BG\to \cG$ is an equivalence of groupoids, where we consider $BG$ as a groupoid with a single object.

\begin{ex}
Let $S$ be a finite  $G$-set. Then we can form the action groupoid $S\curvearrowleft G$. It is connected if and only if the $G$-action on $S$ is transitive. All finite connected groupoids are of this form for some group $G$ and transitive $G$ set $S$. \hB
\end{ex}

\begin{rem} An object $A$ in $\cG\nCalg:=\Fun(\cG,\nCalg)$ is
a collection $(A_{x})_{x\in \Ob(\cG)}$ of $C^{*}$-algebras together with  homomorphisms 
$g:A_{x}\to A_{x'}$ for all $g:x\to x'$ satisfying the associativity rule. 
A morphism $A\to B$ is a family of morphisms
$(f_{x}:A_{x}\to B_{x})_{x\in \Ob(\cG)}$ that is compatible with the morphisms of $\cG$.

If $\cG$ is connected, then the
  restriction along the  equivalence of groupoids $BG\to \cG$ induces an equivalence of categories   $\cG\nCalg\stackrel{\simeq}{\to} G\nCalg$. In particular,
  every object of $\Fun(\cG,\nCalg)$ is determined up to unique isomorphism by
the object $A_{x_{0}}$ in $G\nCalg$, and 
a morphism  is uniquely determined by  
a $G$-equivariant homomorphism $f_{x_{0}}:A_{x_{0}}\to B_{x_{0}}$.
 \hB \end{rem}
 The construction of $\kk^{G}:\nCalg\to \KK^{G}$ from $G\nCalg$ can be adjusted to the construction of a functor 
$$\kk^{\cG}:\cG\nCalg\to \KK^{\cG}$$  such that restriction along $\kk^{\cG}$ 
has the universal property \eqref{sdfbdfbvsfdvfdvsfvsfdvsfdvsf} with $G$ replaced by $\cG$. 

 Generalizing the restriction functors $\Res^{G}_{H}$ from \eqref{vsoijo3rfeffsvfvs}, a morphism of finite connected groupoids $\cH\to \cG$ induces a restriction functor
$\Res^{\cG}_{\cH}:\KK^{\cG}\to \KK^{\cH}$.
In particular,  the restriction
$\Res^{\cG}_{BG}:\KK^{\cG}\to \KK^{BG}\simeq \KK^{G}$ is an equivalence, where the right term
just denote the $G$-equivariant $\KK$-category in the usual notation. \begin{rem}\label{jrgergoergweferfweff}
The association $\cG\to \KK^{\cG}$ is part of the structure of a global category $\underline{\KK}$, see e.g.
\cite[Ex, 5.10]{Cnossen:2023aa}. \hB
\end{rem}

\begin{ex} Let $G$ be a finite group and  $f:S\to T$ be an equivariant map of transitive $G$-sets. In induces a map of finite connected
  groupoids $ S\curvearrowleft G \to T\curvearrowleft G$ that has the special property of being a faithful isofibration, i.e., for every morphism $g:t\to t'$ in $ T\curvearrowleft G$ and $s$ in $f^{-1}(t)$ there exists a unique morphism $s\to s'$ with $f(s\to s')=(t\to t')$. It is implemented by the same group element $g$. \hB
  \end{ex}

From now on, we consider either the minimal or maximal tensor products of $C^{*}$-algebras.

We call $\cG$ is  proper if the automorphism groups of its objects are finite.
Let $\phi:\cH\to \cG$ be a faithful isofibration between proper connected groupoids.
\begin{ddd}\label{jioregoperwgrefwrefwref}We define the multiplicative induction functor
$$\Ind_{\cH}^{\otimes,\cG}:\cH\nCalg\to \cG\nCalg$$ as follows:
\begin{enumerate}
\item objects: It sends $(A_{x})_{x\in \Ob(\cH)}$ to $(\bigotimes_{x\in \phi^{-1}(\{y\}) } A_{x})_{y\in \Ob(\cG)}$. To a morphism 
 $g:y\to y'$  we associate the morphism $$(\otimes_{x\in f^{-1}(\{y\})} g_{x})_{y\in \Ob(\cG)}:(\bigotimes_{x\in  \phi^{-1}(\{y\}) } A_{x})_{y\in \Ob(\cG)}\to
(\bigotimes_{x'\in  \phi^{-1}(\{y'\}) } A_{x})_{y\in \Ob(\cG)}\ ,$$
where $g_{x}:x\to x'$ is the unique lift of $g:y\to y'$ with domain $x$.
\item morphisms: The functor sends a morphism $( f_{x}:(A_{x})_{x\in \Ob(\cH)}\to (B_{x})_{x\in \Ob(\cH)})$ in
$\cH\nCalg$ to the morphism $$(\otimes_{x\in  \phi^{-1}(\{y\})} f_{x})_{y\in \Ob(\cG)}:(\bigotimes_{x\in  \phi^{-1}(\{y\}) } A_{x})_{y\in \Ob(\cG)}\to
(\bigotimes_{x'\in  \phi^{-1}(\{y'\}) } B_{x})_{y\in \Ob(\cG)}\ .$$
\end{enumerate}
\end{ddd}

\begin{ex}The most important example is the map $G\curvearrowleft G\to *\curvearrowleft G\simeq BG$ for a finite group $G$.  We obtain the power  functor
$$(-)^{\otimes G}:\nCalg\simeq \Fun(G\curvearrowleft G,\nCalg)\xrightarrow{\Ind_{ G\curvearrowleft G}^{\otimes,*\curvearrowleft G}} 
 \Fun(*\curvearrowleft G,\nCalg)\simeq G\nCalg\ .$$
 It sends a $C^{*}$-algebra $A$ to the $G$-$C^{*}$-algebra $A^{\otimes G}:=\bigotimes_{g\in G}A$
 with the $G$-action on the index set.  \hB
\end{ex}

The main result of the present subsection is the verification that the multiplicative induction functors descend to $KK$-theory.  

\begin{theorem}\label{weijogwegerfwefwe}
If $\cH\to \cG$ is a faithful isofibration between proper  connected groupoids, then
 there exists an essentially unique factorization
 $$\xymatrix{\cH\nCalg \ar[dr]_{\kk^{\cH} }\ar[r]^{\Ind_{ \cH}^{\otimes,\cG}}&\cG\nCalg \ar[r]^{\kk^{\cG} }&\KK^{\cG}\\&\KK^{\cH} \ar@{..>}[ur]_{\Ind_{ \cH}^{\otimes,\cG}}&}\ .$$
  \end{theorem}
\begin{proof}
In a first step, we reduce the problem to the more common case of groups.
We fix a point $y_{0}$ in $\Ob(\cG)$ and let $G$ denote the stabilizer. Then there exists $x_{0}\in \Ob(\cH)$ such that $f(x_{0})=y_{0}$, and we let $H$ denote the stabilizer of $x_{0}$. Note that $G$ and $H$ are finite  groups by our properness assumption.  In view of the cube 
$$\xymatrix{\KK^{\cH}\ar@{..>}[ddd]^{\Ind_{\cH}^{\otimes,\cG}}\ar[rrr]_{\simeq}^{\Res^{\cH}_{BH}}&&&\KK^{H}\ar@{..>}[ddd]^{(-)^{\otimes G/H}_{\kk}}\\&\ar[ul]^{\kk^{\cH}}\cH\nCalg\ar[r]^{\simeq}\ar[d]^{\Ind_{\cH}^{\otimes,\cG}} & \ar[ur]^{\kk^{H}}H\nCalg\ar[d]^{(-)^{\otimes G/H}} &\\& \ar[dl]^{\kk^{\cG}}\cG\nCalg\ar[r]^{\simeq} & G \nCalg\ar[dr]^{\kk^{G}} &\\\KK^{\cG}\ar[rrr]_{\simeq}^{\Res^{\cG}_{BG}}&&&\KK^{G}} $$
it suffices to construct $(-)^{\otimes G/H}_{\kk}$ on the right.

The explicit description of the symmetric monoidal  functor $$(-)^{\otimes G/H}:H\nCalg\to G\nCalg$$ depends on the choice of  a functor $\cH\to BH$
and a natural isomorphism between $\cH\to BH\to \cH$ and $\id_{\cH}$. For both we choose a collection
of morphisms $h_{x}:x\to x_{0}$ for all $x$ in $\Ob(\cH)$. The  functor 
$\cH\to BH$ is given on morphisms by $(h:x\to x')\mapsto h_{x'}h h_{x}^{-1}$, and the components of the transformation are $h_{x}^{-1}:x_{0}\to x$. The functor $(-)^{\otimes G/H}$ sends a $H$-$C^{*}$-algebra $A$ to 
$$A^{\otimes G/H}:=\bigotimes_{x\in f^{-1}(y_{0})} A\ .$$
The action of $g$ in $G$ on this tensor product ist then given by
\begin{equation}\label{dafasdfqffadfs}\otimes_{x\in f^{-1}(y_{0})} h_{gx}g_{x}h_{x}^{-1}: \bigotimes_{x\in f^{-1}(y_{0})} A\to \bigotimes_{x\in f^{-1}(y_{0})} A\ .
\end{equation} 
Note that $h_{gx}g_{x}h_{x}^{-1}\in H$ for all $g$ in $G$ and $x$ in $f^{-1}(y_{0})$.

 \begin{rem}
 Our main reason to work with groupoids instead of groups is that the multiplicative
 induction functor in \cref{jioregoperwgrefwrefwref} is given by a completely canonical formula. In contrast, the formula  \eqref{dafasdfqffadfs} for the power functor $(-)^{\otimes G/H}$ depends on the choices made in the construction of an inverse
 of the inclusion $BH\to \cH$. \hB
 \end{rem}

We first restrict to separable $C^{*}$-algebras and indicate this by a subscript $\sepa$. 
Recall from \cite[Def. 4.1]{KKG} that a functor $C$ between categories of $C^{*}$-algebras with actions of groups   is $\Ind$-$s$-finitary if for every $A$ the inductive system $(C(A')^{C(A)})_{A'\subseteq_{\sepa} A}$
(see \cref{wtkopgwfrefrefw}.\ref{irwegjoergrffwre} for notation) of images of $C(A')\to C(A)$ is cofinal in $B\subseteq_{\sepa} C(A)$,
and the canonical map $(C(A'))_{A'\subseteq_{\sepa} A}\to (C(A')^{C(A)})_{A'\subseteq_{\sepa} A}$ is an isomorphism of inductive systems.
 
 \begin{prop}\label{weogjrpwegfrfwrfwre} Assume that $G$ is a finite group and $H$ is a  subgroup. Then there  exists an essentially unique symmetric monoidal  factorization
 $$\xymatrix{H\nCalg_{\sepa}\ar[dr]_{\kk^{H}_{\sepa}}\ar[r]^{ (-)^{\otimes G/H}}&G\nCalg_{\sepa}\ar[r]^{\kk^{G}_{\sepa}}&\KK^{G}_{\sepa}\\&\KK^{H}_{\sepa}\ar@{..>}[ur]_{(-)^{\otimes G/H}_{\kk_{\sepa}}}&}\ .$$
 In addition $ (-)^{\otimes G/H}$ is $\Ind$-$s$-finitary
\end{prop}
\begin{proof}
The proof of the second statement is analogous to 
 the proof of \cite[Prop. 3.8]{KKG}. We now  concentrate on the first assertion.
Since $\kk^{H}_{\sepa}$ is a Dwyer-Kan localization at the $\kk_{\sepa}^{H}$-equivalences by \cite[Def. 2.1]{KKG}, it suffices to show that
the composition $\kk_{\sepa}^{G}\circ (-)^{\otimes G/H}$ sends $\kk^{H}_{\sepa}$-equivalences to equivalences.


We will use the construction of the functor $\kk^{H}_{\sepa}:H\nCalg_{\sepa}\to\kk^{H}_{\sepa}$ as a sequence of localizations (see  \cite{Bunke:2023aa} for the non-equivariant case, the equivariant case is completely analogous).  
$$\hspace{-1cm}\kk^{H}_{\sepa}:H\nCalg_{\sepa}\xrightarrow{L_{h}}  H \nCalg_{\sepa,h}\xrightarrow{L_{K_{H}}} L_{K_{H}} H \nCalg_{\sepa,h}
\xrightarrow{L_{\se}}  L_{K_{H}}H \nCalg_{\sepa,h,\se}\xrightarrow{L_{\group}}\KK^{H}_{\sepa} \ .$$
\begin{enumerate} \item
$L_{h} $ presents the symmetric monoidal  left-exact $\infty$-category   $H\nCalg_{\sepa,h} $   as the symmetric monoidal Dwyer-Kan localization of $H\nCalg_{\sepa}$ at the homotopy equivalences. 
The mapping spaces in $H  \nCalg_{\sepa,h}$ are represented by the topological mapping spaces in $ H\nCalg_{\sepa} $, which are
given by the $\Hom$-sets with the point norm topology. 
The localization $L_{h}$ preserves products, coproducts and sends exact sequences  $0\to A\to B\xrightarrow{\pi} C \to 0$
in $ H\nCalg_{\sepa}$ with $\pi$ a Schochet fibration \cite{MR650021} (i.e., $\Hom_{H\nCalg_{\sepa}}(D,\pi)$ is a Serre fibration for every $D$ in $H\nCalg$)  to fibre sequences $L_{h}(A)\to L_{h}(B)\to L_{h}(C)$.
 
\item\label{gwerogpefrefrfw}  $L_{K_{H}}$ is the symmetric monoidal Bousfield localization of $H \nCalg_{\sepa,h}$ at the full subcategory   the $K_{H}$-stable algebras.
The  left-exact  localization functor  $L_{K_{H}}$ sends $A$ to $A \otimes  K_{H}$, where $K_{H}:=K(L^{2}(H)\otimes \ell^{2}) $.  It   presents the pre-additive  $\infty$-category $L_{K_{H}} H \nCalg_{\sepa,h}$ as the Dwyer-Kan localization of $H\nCalg_{\sepa,h}$ at the left-upper corner inclusion
$A \to  A  \otimes K_{H} $  for all $A$ in $\nCalg_{\sepa}$ using the projection $p$ in $K_{H}$ projecting on the subspace of $L^{2}(H)\otimes \ell^{2}$ generated by $\const_{1}\otimes e_{0}$. It exists since we assume that $H$ is finite.
\item   For any    semisplit exact sequence $0\to A\to B\xrightarrow{\pi} C \to 0$ in $H\nCalg$ (i.e., $\pi$ admits  an equivariant completely positive contractive right-inverse),
 we can consider the map $\iota_{\pi}:A \to C(\pi)$ (see Remark \ref{wtiogrtgs} below for notation) representing the   map from $A$ into the homotopy fibre of $\pi$. The left-exact functor $L_{\se}$ presents the symmetric monoidal left-exact $\infty$-category  
$L_{K_{H}}H\nCalg_{\sepa,h,\se}$ as  
 the symmetric monoidal  Dwyer-Kan localization at the closure under $2$-out-of-$3$ and pull-backs of the set of these morphisms $\iota_{\pi}$  for all semisplit exact sequences. The functor $H\nCalg_{\sepa}\to  L_{K_{H}}H\nCalg_{\sepa,h,\se}$ sends
exact sequences  $0\to A\to B\xrightarrow{\pi} C \to 0$ with $\pi$ a Schochet fibration or admitting  an equivariant completely positive contractive split 
to fibre sequences. 
\item   $L_{\group}$ is the symmetric monoidal right Bousfield localization of $ L_{K_{H}}H\nCalg_{\sepa,h,\se }$ at the full subcategory of group objects.  The  left-exact localization functor $L_{\group}$ sends $A$  to its two-fold suspension $S^{2}(A)$. 
 We have a natural   transformation $\beta:S^{2}\to \id$ of endofunctors of 
$ L_{K_{H}}H\nCalg_{\sepa,h,\se}$,  and $L_{\group}$ presents $\KK^{H}_{\sepa}$ as the Dwyer-Kan localization of $ L_{K_{H}}H\nCalg_{\sepa,h,\se}$  at the 
 components of this transformation. 
\end{enumerate}

\begin{rem}\label{wtiogrtgs}
If $ B\xrightarrow{\pi} C $ is  a map of $C^{*}$-algebras, then we   form
the mapping zylinder $$Z(\pi):=\{(f,b)\in C([0,1],C)\oplus B\mid f(0)=\pi(b) \}\ .$$ It comes with a Schochet fibration $
\hat \pi:Z(\pi)\to C$ given by $(f,b)\mapsto f(1)$ and a homotopy equivalence $h:B\to Z(\pi)$ given by $b\mapsto (\const_{\pi(b)},b)$. The kernel $$C(\pi):= \{(f,b)\in Z(\pi)\mid f(1)=0\}$$ of $\hat \pi$ represents the homotopy fibre of $\pi$, and we have  a map $\iota_{A}:A\to C(\pi)$ given by $a\mapsto (0,a)$, where $A:=\ker(\pi)$.

In order to see these statements, note that the Schochet fibrations and homotopy equivalences are part of a fibration category structure on $ H\nCalg_{\sepa}$ which models the localization $L_{h}$, see  \cite{Uuye:2010aa} for the non-equivariant case and \cite[Prop. 2.10]{KKG} for the equivariant case (for discrete groups).
 \hB
\end{rem}

We now consider the following diagram
$$\xymatrix{H\nCalg_{\sepa}\ar[ddrrr]^{\kk^{G}_{\sepa}\circ (-)^{\otimes  G/H}}\ar[d]^{L_{h}}&&& \\H\nCalg_{\sepa,h}  \ar[d]^{L_{K_{H}}}\ar@{.>}[drrr]^{F_{1}}&&&\\ L_{K_{H}}H\nCalg_{\sepa, h} \ar[d]^{L_{\se}}\ar@{.>}[rrr]^{F_{2}}&&& \KK^{G}_{\sepa} 
\\ L_{K_{H}}H \nCalg_{\sepa,h,\se}\ar[d]^{L_{\group}}\ar@{.>}[urrr]^{F_{3}}&&&\\ \KK^{H}_{\sepa}\ar@{..>}[uurrr]_{ (-)_{\kk_{\sepa}}^{\otimes G/(H}}}\ .$$
We will obtain the four  symmetric monoidal   factorizations from the universal properties of the localizations at the left.

 \begin{enumerate}
 \item The functor  $ (-)^{\otimes G/H}:H\nCalg_{\sepa}\to G\nCalg_{\sepa}$  defined above is continuous and hence preserves homotopy equivalences. 
 Since $\kkG_{\sepa}$  is  homotopy invariant, also the composition $\kk^{G}_{\sepa}\circ (-)^{\otimes G/H}:G\nCalg_{\sepa}\to \KK^{G}_{\sepa}$ 
  is homotopy invariant.  We get the factorization $F_{1}$ from the universal property of $L_{h}$.
 \item We consider $u:\C\to K_{H}$ induced by the projection  $p$ described in \cref{gwerogpefrefrfw}.      For $A$ in $H\nCalg_{\sepa,h}$ we get the   left upper corner embedding $\id_{A}\otimes u :A  \to A \otimes K_{H} $. 
  We must show that the induced map $F_{1}(\id_{A}\otimes u):F_{1}(A)\to F_{1}(A\otimes K)$  is an equivalence.
   Since $F_{1}$ is symmetric monoidal we have an  equivalence $F_{1}(\id_{A}\otimes u) \simeq \id_{F_{1}(A)}\otimes F_{1}(u)$ and it suffices to show that $F_{1}(u)$ is an equivalence.  
   
   We have an  isomorphism  $ K_{H}^{\otimes G/H} \cong \bigotimes_{x\in f^{-1}(y_{0})} K_{H}$  in $G\nCalg$ under which the
  homomorphism $  u^{\otimes G/H}$  is induced by the $G$-invariant projection $\otimes_{x\in f^{-1}(\{y_{0}\})} p $.   Since the latter is a one-dimensional, we can conclude from the $K_{G}$-stability of 
  $\kk^{\cG}_{\sepa}$ that  $ \kk^{G}_{\sepa}( u^{\otimes G/H})$ is an equivalence.
  Hence  we get the factorization $F_{2}$ from the universal property of $L_{K_{H}}$.
 \item \label{ejigweporfrefrwfr}Let $0\to A\to B\xrightarrow{\pi}C\to 0$ be a semisplit exact sequence in $H\nCalg_{\sepa}$. Using \cref{wtiogrtgs}, we form the diagram
$$\xymatrix{0\ar[r]&A\ar[r]\ar[d]^{\iota_{\pi}}&B\ar[d]^{h}\ar[r]^{\pi} &C\ar[r]\ar@{=}[d]&0\\
0\ar[r]&C(\pi)\ar[r]&Z(\pi)\ar[r]^{\hat \pi}&C\ar[r] &0}$$ in $H\nCalg_{\sepa}$,
where $h$ is a homotopy equivalence.   Since the target of $F_{2}$ is stable and $F_{2}$ is left-exact 
it suffices (i.e. we do not have to discuss the closure of the set of these maps under pull-back and $2$-out-of-$3$) to show that $F_{2}(\iota_{\pi})$ is an equivalence for all  exact sequences as above.
Equivalently, we must show that the induced map $$ \iota_{\pi}^{\otimes G/H}:  A^{\otimes G/H}
\to   C(\pi)^{\otimes G/H}$$ is a $\kk^{G}_{\sepa}$-equivalence.

For $m$ in $\{0,\dots,|G/H|\}$ we let $I_{m}$ be the $G$-invariant ideal in $B^{\otimes G/H}$ generated by elementary  tensors which have at least $m$ factors in $A$.  The family $(I_{m})_{m=0,\dots,|G/H|}$ is  a   decreasing filtration
of $B^{\otimes G/H}$ by  $G$-invariant ideals such that $I_{|G/H|}\cong A^{\otimes G/H}$, 
$I_{1}=\ker(B^{\otimes G/H}\to C^{\otimes G/H})$, and $I_{0}= B^{\otimes G/H}$.
We introduce the abbreviation   $Z:=f^{-1}(\{y_{0}\})$ for the finite $G$-set given by the fibre of $f$ over $y_{0}$. Then $I_{m}\cong \sum_{F\in \cP_{m}( Z) } A^{\otimes F}\otimes  B^{\otimes  Z\setminus F}$  and   $I_{m}/I_{m+1}  \cong \bigoplus_{F\in \cP_{m}( Z)  } A^{\otimes F}\otimes  C^{\otimes Z\setminus F}$, where $\cP_{m}(Z):=\{F\in \cP(Z)\mid |F|=m\}$ is the set of subsets of $Z$ with $m$ elements. 
Let $s:C\to B$ denote the $H$-equivariant cpc split.
The sequence
\begin{equation}\label{qfjoipqwfevwfoqvfvvfqwdqwdqd}0\to I_{m+1}\to I_{m}\to I_{m}/I_{m+1}\to 0
\end{equation} has an induced $G$-equivariant  completely positive split  given by 
$$ \hat s:\sum_{F\in \cP_{m}( Z) } \id_{A}^{\otimes F}\otimes   s^{\otimes Z\setminus F}   : \bigoplus_{F\in \cP_{m}( Z)  } A^{\otimes F}\otimes  C^{\otimes Z\setminus F} \to  \sum_{F\in \cP_{m}( Z) } A^{\otimes F}\otimes  B^{\otimes  Z\setminus F}\ .$$ 
Using \cite[Rem 2.5]{blackadar} one can modify $\hat s$ to a contractive completely positive split $\tilde s$. Since the construction involves non-canonical choices the latter might not be $G$-equivariant anymore. However, since $G$ is finite, we can finally find an equivariant cpc split by averaging $\tilde s$.\footnote{We thank S. Nishikawa for pointing out a mistake at this point in an earlier version of this paper.}

 We write $(J_{m})_{m=0,\dots,|Z|}$ for the corresponding sequence of ideals for the sequence $0\to C(\pi)\to Z(\pi)\to C\to 0$. We then have a map
$$\xymatrix{0\ar[r]&I_{m+1}\ar[r]\ar[d]^{\iota_{m+1}}&I_{m}\ar[d]^{h_{m}}\ar[r]^{\pi} &I_{m}/I_{m+1}\ar[r]\ar[d]^{q_{m}}&0\\
0\ar[r]&J_{m+1}\ar[r]&J_{m}\ar[r]&J_{m}/J_{m+1}\ar[r]&0}$$
between semisplit exact sequences in $G\nCalg_{\sepa}$ which induces a map of fibre sequences  
in $\KKG_{\sepa}$ upon application of $\kkG_{\sepa}$. We want to conclude by induction over the order of $G$, the number of elements of $Z$,   and  $m$ that the maps $\iota_{m}$  are   $\kkG_{\sepa}$-equivalences.
We start with $m=0$. Then $h_{0}=h^{\otimes Z}$ is a homotopy equivalence and $q_{0}=\id_{C^{\otimes Z}}$ is an isomorphism. Therefore $\iota_{1}$ is a $\kkG_{\sepa}$-equivalence. 
For the induction step from $m$ to $m+1$, we assume that $m<|Z|$. We observe that
\begin{equation}\label{vsdfvoijiowerjovervdfsvfv}q_{m}\cong\bigoplus_{[F]\in \cP_{m}(Z)/G} \Ind_{\sepa,G_{F}}^{G}(\iota_{\pi}^{\otimes F}\otimes \id_{C}^{\otimes Z\setminus F})\ ,
\end{equation}
where $G_{F}$ denotes the stabilizer of $F$ in $G$.
 Note $\iota_{\pi}^{\otimes F}$ is a $\kk^{G_{F}}_{\sepa}$-equivalence by our  induction hypothesis.
It follows that $ \Ind_{\sepa,G_{F}}^{G}(\iota_{\pi}^{\otimes F}\otimes \id_{C}^{\otimes Z\setminus F})$ is a $\kk^{G}_{\sepa}$-equivalence and we can conclude that $\iota_{m+1}$ is a $\kkG_{\sepa}$-equivalence.
Since $\iota_{|Z|}=\iota_{\pi}^{\otimes Z}=\iota^{\otimes G/H}$ this finishes the argument.

  \item  We must show that $F_{3}(\beta_{A})$ is an equivalence in $\KKG_{\sepa}$ for every $A$ in $\nCalg_{\sepa}$. 
       Note that $ C_{0}(\C)^{\otimes G/H}\cong C_{0}(\C^{Z})$ in $G\nCalg_{\sepa}$, where the $G$-action on the right-hand side is implemented by an action on the complex vector space $\C^{Z}$ by $\C$-linear automorphisms. 
  It follows from Kasparov's Bott periodicity that $\kkG_{\sepa}( C_{0}(\C^{Z}))$ is a tensor-invertible object of $\KKG_{\sepa}$.
  Since $F_{3}$ is symmetric monoidal, we have an equivalence
  $F_{3}(\beta_{A}\otimes \id_{C_{0}(\C)})\simeq F_{3}(\beta_{A})\otimes  \kkG_{\sepa}( \id_{C_{0}(\C^{Z})})  $.
  Hence it suffices to show that $F_{3}(\beta_{A}\otimes \id_{C_{0}(\C)})$ is an equivalence.
  But since   
  $\beta_{A}\otimes  \id_{C_{0}(\C)}\simeq \beta_{A\otimes C_{0}(\C)}\simeq \beta_{S^{2}(A)}$ and 
  $S^{2}(A)$ is a group object in $L_{K_{H}} \nCalg_{\sepa,h ,\se}$,  we can conclude that already $\beta_{S^{2}(A)}$ is an equivalence. 
  We get the factorization $(-)^{\otimes G/H}_{\kk_{\sepa}}:\KK^{H}_{\sepa}\to \KKG_{\sepa}$ from the universal property of $L_{\group}$.
     \end{enumerate}
     \end{proof}

 \begin{prop}\label{iorgoergfgsfg} Assume that $G$ is a finite group and $H$ is a subgroup. Then there exists an essentially unique factorization $$\xymatrix{H\nCalg \ar[dr]_{\kk^{H} }\ar[rr]^{\kk^{\cG} \circ  (-)^{\otimes G/H}}&&\KK^{\cG} \\&\KK^{H} \ar@{..>}[ur]_{  (-)^{\otimes G/H}_{\kk}}&}$$
 such that $(-)_{\kk}^{\otimes G/H}$ preserves  all filtered colimits and compact objects.
 A addition
 $(-)_{\kk}^{\otimes G/H}$ has a symmetric monoidal refinement. 
   \end{prop}

\begin{proof} 
By  definition $y^{H}:\KK^{H}_{\sepa}\to \KK^{H}$ and $y^{G}:\KK^{G}_{\sepa}\to \KKG$  present $\KK^{H}$ and $\KKG$ as the $\Ind$-completions of $\KK^{H}_{\sepa}$  and $\KKG_{\sepa}$, and their images consist of the compact objects in $\KK^{H}$ or $\KKG$, respectively.  Since the desired functor 
   $(-)_{\kk}^{\otimes G/H}:\KK^{H}\to \KKG$ is required to preserve filtered colimits and compact objects it is necessarily the   essentially unique such functor    
 such that $$\xymatrix{\KK^{H}_{\sepa}\ar[r]^{(-)^{\otimes G/H}_{\kk_{\sepa}}}\ar[d]^{y} &\KKG_{\sepa} \ar[d]^{y^{G}} \\\KK^{H} \ar[r]^{(-)_{\kk}^{\otimes G/H}} &\KKG } 
 $$ commutes.
 The symmetric monoidal structure of $(-)^{\otimes Z}_{\kk_{\sepa}}$  induces a the symmetric monoidal refinement of
 $(-)_{\kk}^{\otimes Z}$.
 
 We now consider the following diagram:
 \begin{equation}\label{}
  \xymatrix{\ar@/^-2cm/[ddd]_{\incl}H\nCalg_{\sepa}\ar[r]^{(-)_{\sepa}^{\otimes G/H}}\ar[d]^{\kk^{H}_{\sepa}} &G\nCalg_{\sepa} \ar[d]^{\kkG_{\sepa}}\ar@/^2cm/[ddd]^{\incl} \\\KK^{H}_{\sepa}\ar[r]^{(-)_{\kk_{\sepa}}^{\otimes G/H} }\ar[d]^{y^{H}} &\KKG_{\sepa} \ar[d]^{y^{G}} \\\KK^{H} \ar[r]^{(-)_{\kk^{H}}^{\otimes G/H}} &\KKG\\ H\nCalg \ar[u]^{\kk^{H}}\ar[r]^{(-)^{\otimes G/H}} 
  &G\Calg\ar[u]^{\kkG} }\ .  \end{equation}

Our task is to show that the lower square commutes in an essentially unique way. We already know that all other squares commute. They provide a filler of the pentagon
$$\xymatrix{&H\nCalg_{\sepa}\ar[dl]_{\incl}\ar[dr]^{(-)^{\otimes G/H}\circ \incl}&\\ H\nCalg\ar@{..>}[rr]^{(-)^{\otimes G/H} } \ar[d]^{\kk}&&G\nCalg \ar[d]^{\kkG}\\\KK^{H}\ar[rr]^{(-)_{\kk}^{\otimes G/H}}&& \KKG}\ .$$
The   left composition    from $H\nCalg_{\sepa}$ to $\KKG$ is clearly $s$-finitary, and
the right composition is   $s$-finitary by the second assertion of   \cref{weogjrpwegfrfwrfwre} and \cite[Lem. 4.2]{KKG}.
Hence we get an essentially unique filler of the lower square whose composition with the obvious filler of the upper triangle is the filler of the pentagon. \end{proof}

\end{proof}

From now on we will often write $(-)^{\otimes G/H}$ also for $(-)^{\otimes G/H}_{\kk}$.

      \begin{rem}  
In order to show the first assertion of \cref{weijogwegerfwefwe} in the case of a trivial group $H$,
  one could use the characterization of   $\kk_{\sepa}:\nCalg_{\sepa}\to \KK_{\sepa}$ as the universal  functor to a stable $\infty$-category that sends $\kk_{\sepa}$-equivalences to equivalences  \cite{LN},  \cite[Def. 4.2]{KKG}. 
Indeed, one can show that the composition $$H\nCalg_{\sepa}\xrightarrow{(-)^{\otimes G/H}} G\nCalg_{\sepa}\xrightarrow{\kkG_{\sepa}} \KKG_{\sepa}$$
sends $\kk_{\sepa}$-equivalence to equivalences  \cite[Lem. 2.4]{Chakraborty:2022aa}. 
The idea for the filtrations used in the Step \ref{ejigweporfrefrwfr} above is taken  from \cite{Chakraborty:2022aa}. 
\hB
\end{rem}

Let $B$ be in $\KK^{H}$ and $p$ be a prime.
\begin{lem}\label{egrjiweorgerferwfrewf}
If $p$ acts as an equivalence on $B$, then it acts as an equivalence on $B^{\otimes G/H}$.
\end{lem}
\begin{proof}
Since $p$ acts as an equivalence on $B$,  we have an equivalence
$B\simeq B\otimes \beins_{\KK^{H}}[p^{-1}]$ in $\KK^{H}$. Since
 the multiplicative induction functor $(-)^{\otimes G/H}
$ is symmetric monoidal, we obtain an equivalence 
$ B^{\otimes G/H}\simeq B^{\otimes G/H}\otimes \beins_{\KK^{H}}[p^{-1}]^{\otimes G/H}$.
It suffices to show that
$p$ acts as an equivalence on  $\beins_{\KK^{H}}[p^{-1}]^{\otimes G/H}$. We have
$$\beins_{\KK^{H}}[p^{-1}]\simeq \colim \left( \beins_{\KK^{H}}\stackrel{p}{\to}\beins_{\KK^{H}}\stackrel{p}{\to}\beins_{\KK^{H}}\stackrel{p}{\to}\dots \right) \ .$$
 Since $(-)^{\otimes G/H}
$  preserves filtered colimit we get 
$$\beins_{\KK^{H}}[p^{-1}]^{\otimes G/H}\simeq  \colim \left( \beins^{\otimes G/H}_{\KK^{H}}\stackrel{p^{\otimes G/H}}{\to}\beins^{\otimes G/H}_{\KK^{H}}\stackrel{p^{\otimes G/H}}{\to}\beins^{\otimes G/H}_{\KK^{H}}\stackrel{p^{\otimes G/H}}{\to}\dots \right) \ .$$
We now use that $\beins^{\otimes G/H}_{\KK^{H}}\simeq \beins_{\KK^{G}}$ and that the map 
$p^{\otimes G/H}$ is the action of $p^{|G/H|}$. Hence  
$$\beins_{\KK^{H}}[p^{-1}]^{\otimes G/H}\simeq  \colim \left( \beins_{\KK^{G}}\stackrel{p^{|G/H|}}{\to}\beins_{\KK^{G}}\stackrel{p^{ |G/H|}}{\to}\beins_{\KK^{G}}\stackrel{p^{ |G/H|}}{\to}\dots \right)\simeq \beins_{\KK^{G}}[p^{-1}] \ .$$
  \end{proof}

\subsection{A spectral sequence}\label{krwfopqwewfewdffqwf}

Let $G$ be a finite group and  $Z$ be a finite $G$-set. It has a decomposition $Z\cong \bigsqcup_{S\in Z/G} S$ into $G$-orbits and each of the $G$-orbits $S$ in this decomposition is isomorphic to $G/H$ for some subgroup $H$ of $G$.
If $A$ is in $\KK$, then by \cref{iorgoergfgsfg} we can form the power $$A^{\otimes Z}:= \bigotimes_{S\in Z/G} A^{\otimes S}$$ in $\KK^{G}$.
We now assume that $A$ is part of a fibre sequence
\begin{equation}\label{fewdqewdqwdqeded}I\to A\to B
\end{equation} in $\KK$. In this section, we construct a tower in $\KK^{G}$ that calculates $A^{\otimes Z}$ in terms of powers 
of $I$ and $B$. The associated spectral sequence will be used for the calculation in \cref{gwegire0fwrefwrefw}.

For simplicity, we assume that  the sequence belongs to $\KK_{\sepa}$. Then we can realise this fibre sequence 
by a semisplit  exact sequence of $C^{*}$-algebras. So changing notation
we start with 
a
semisplit exact sequence
\begin{equation}\label{vsffsdvfrqfefrf}0\to I\to A\to B\to 0
\end{equation} in $\Calg_{\sepa}$
which gives rise to a fibre sequence
\begin{equation}\label{fqewifhqiofsffeadf}\kk(I)\to \kk(A)\to \kk(B)
\end{equation}
replacing \eqref{fewdqewdqwdqeded}.

 Similarly as in the proof of \cref{weogjrpwegfrfwrfwre},
  we have a decreasing filtration \begin{equation}\label{wregfewr}I^{\otimes Z}=I^{Z}_{|Z|}\subseteq I^{Z}_{|Z|-1}\subseteq\dots I^{Z}_{1}\subseteq I_{0}=A^{\otimes Z}
\end{equation}
of the $G$-$C^{*}$-algebra $A^{\otimes Z}$  by $G$-invariant ideals. Thereby the
 ideal $I^{Z}_{n}$ is generated by  the elementary tensors $\otimes_{z\in Z} a_{z}$ with at least $n$ factors in $I$.
The intermediate quotients  are given by 
\begin{equation}\label{vdfsvrwgweg}\frac{I^{Z}_{n}}{I^{Z}_{n+1}}\cong \bigoplus_{F\in \cP_{|Z|-n}(Z)}    I^{\otimes Z\setminus F} \otimes B^{\otimes F}\ ,
\end{equation} 
where $\cP_{n}(Z)$ is the $G$-set of $n$-element subsets of $Z$.
The action of $G$ on the right-hand side is  induced by the action on $\cP_{|Z|-n}(Z)$
and permutation of tensor functors. 
 For every $G$-orbit $\cF$ in $\cP_{|Z|-n}(Z)/G$ we have a $G$-invariant subalgebra
 \begin{equation}\label{qfewfqwdqfqewfwqef}\bigoplus_{F\in \cF}    I^{\otimes Z\setminus F} \otimes B^{\otimes F}\cong \Ind_{G_{F_{0}}}^{G}( I^{\otimes Z\setminus F_{0}} \otimes B^{\otimes F_{0}})    \ , 
\end{equation} 
 where $G_{F_{0}}$ denotes the stabilizer of a point $F_{0}$ in $\cF$ and acts on 
  $I^{\otimes Z\setminus F_{0}} \otimes B^{\otimes F_{0}}$ via its action on $Z\setminus F_{0}$ and $F_{0}$ by permutation of tensor functors.
The  isomorphism in \eqref{vdfsvrwgweg} is induced by the obvious projections
$\pi^{Z}_{F}:I_{n}^{Z}\to A^{\otimes Z\setminus F} \otimes B^{\otimes F}$ which actually take values in the subalgebras
$I^{\otimes Z\setminus F} \otimes B^{\otimes F}$.
We further use the notation  $$\pi^{Z}_{n}:= \oplus_{F\in \cP_{|Z|-n}(Z)} \pi^{Z}_{F}:I^{Z}_{n}\to \bigoplus_{F\in \cP_{|Z|-n}(Z)} I^{\otimes Z\setminus F} \otimes B^{\otimes F}$$
for the projection map.
  By  
$$\partial^{Z}_{n}:  \kk^{G}(   \bigoplus_{F\in \cP_{|Z|-n}(Z)/G } I^{\otimes Z\setminus F} \otimes B^{\otimes F})\to \Sigma \kk^{G}(I^{Z}_{n+1})$$ we denote the boundary operator of the semisplit  exact sequence
$$0\to I^{Z}_{n+1}\to I^{Z}_{n}\to \bigoplus_{F\in \cP_{|Z|-n}(Z)} I^{\otimes Z\setminus F} \otimes B^{\otimes F}\to 0\ .$$

The filtration \eqref{wregfewr} gives rise to a tower \begin{equation}\label{ewqfqwoijoiqwedewdq}\xymatrix{0\ar[d]&\\
\kk^{G}(I^{Z}_{|Z|})\ar[r]^{\pi^{Z}_{|Z|}}\ar[d]&\kk^{G}(I^{Z}_{|Z|})\ar@{..>}[ul]_{\partial^{Z}_{|Z|}}\\ \kk^{G}(I^{Z}_{|Z|-1})\ar[r]\ar[r]^-{\pi^{Z}_{|Z|-1}}\ar[d]& \kk^{G}(\bigoplus_{F\in P_{1}(Z)}I^{\otimes Z\setminus F} \otimes B^{\otimes F})\ar@{..>}[ul]_{\partial^{Z}_{|Z|-1}}\\\vdots\ar[d]&\vdots\\ \kk^{G}(I^{Z}_{2})\ar[d]\ar[r]^-{\pi_{2}^{Z}}&   \kk^{G}(\bigoplus_{F\in P_{|Z|-2}(Z)}I^{\otimes Z\setminus F} \otimes B^{\otimes F})\ar@{..>}[ul]_{\partial^{Z}_{2}}\\
\kk^{G}(I^{Z}_{1})\ar[d]\ar[r]^-{\pi_{1}^{Z}}&  \kk^{G}(\bigoplus_{F\in P_{|Z|-1}(Z)}I^{\otimes Z\setminus F} \otimes B^{\otimes F})\ar@{..>}[ul]_{\partial_{1}^{Z}}\\ \kk^{G}(A^{\otimes Z})\ar[r]^{\pi_{0}^{Z}}&
\kk^{G}(B^{\otimes Z})
\ar@{..>}[ul]_{\partial^{Z}_{0}}}
\end{equation}
in $\KK^{G}$.

\begin{rem}
If $\Phi:\KK^{G}\to \Sp$ is an exact functor, e.g. $ \map_{\KK^{G}}(\beins_{\KK^{G}},-)$, then we 
obtain a tower in $\Sp$, giving rise to a spectral sequence converging to
$\pi_{*}\Phi(A^{\otimes Z})$. \hB
\end{rem}

 In the present section we will apply the functor $\K(\widehat \Res^{G}(-)):\KK^{G}\to \Fun(BG,\Mod_{KU}(\Sp))$. 
In this case, the  tower provides a spectral sequence  $(E^{Z,r}_{i,j},d^{Z,r}_{i,j})_{r=1,\dots}$
  for the  $K$-theory groups
$K_{*}(A^{\otimes Z}) $ 
 with the induced $G$-action.
 We use the Adams grading so that $$d^{Z,r}_{i,j}:E^{Z,r}_{i,j}\to E^{Z,r}_{i-1,j+r}\ .$$ 
 
 More precisely, $K_{i}(A^{\otimes Z})$ for fixed $i$ in $\Z$, has a decreasing $|Z|$-step   filtration by $C_{p}$-invariant submodules 
whose associated graded is given by the graded group
$E^{Z,|Z|+1}_{i,*}$.

Our goal in this section is to calculate the first page of the spectral sequence including the differentials.
    We have
\begin{equation}\label{wrgeoigowerpfwerfwerf}E^{Z,1}_{i,j}\cong \bigoplus_{F\in \cP_{|Z|-j}(Z)}K_{i}(I^{\otimes Z\setminus F}\otimes B^{\otimes F})\ .
\end{equation}The first differential $d^{Z,1}_{*,j}$ is induced by the composition
\begin{equation}\label{fvsvsdfvsfdv}\hspace{-0.3cm}d^{Z,1}_{j}:=\pi^{Z}_{j+1}\circ \partial^{Z}_{j}:\bigoplus_{F'\in P_{|Z|-j}(Z)}K(I^{\otimes Z\setminus F'} \otimes B^{\otimes F'}) 
\to \Sigma K(I^{Z}_{j+1})\to  \bigoplus_{F\in P_{|Z|-j-1}(Z)}\Sigma K(I^{\otimes Z\setminus F} \otimes B^{\otimes F})\ .
\end{equation}  We must calculate the  matrix coefficients $d^{Z,1}_{j,F,F'}:K(I^{\otimes Z\setminus F'}\otimes B^{\otimes F'})\to \Sigma K(I^{\otimes Z\setminus F} \otimes B^{\otimes F}) $ of this differential.


\begin{ddd}\label{wekogpwefgrewfwre} The bootstrap class   $\UCT$  is defined to be  the localizing subcategory of $\KK$ generated by the tensor unit $\beins_{\KK}$. \end{ddd}
We will assume that $I$ and $B$ belong to bootstrap class $\UCT$. The restriction of $\K:\KK^{G}\to \Mod_{KU}(\Sp)$ to $\UCT$ is symmetric monoidal so that we can pull-out the tensor powers.
Then 
$d^{Z,1}_{j,F,F'}$ becomes a map
$$d^{Z,1}_{j,F,F'}:K(I)^{\otimes Z\setminus F'}\otimes K(B)^{\otimes F'}\to \Sigma K(I)^{\otimes Z\setminus F} \otimes K(B)^{\otimes F}\ . $$
Let $\partial: K(B)\to \Sigma K(I)$ denote the boundary map associated to the exact sequence \eqref{vsffsdvfrqfefrf}.
\begin{prop}\label{rkogpwergwerfwefwerfwref}
If  $F\subseteq F'$, then  $d^{Z,1}_{j,F,F'}$ is given by the composition 
\begin{eqnarray*}
K(I)^{\otimes Z\setminus F'}  \otimes  K(B)^{\otimes F'}&\simeq&
K(I)^{\otimes Z\setminus F'}\otimes K(B)^{\otimes F'\setminus F}  \otimes  K(B)^{\otimes F}\\&\simeq&
K(I)^{\otimes Z\setminus F'}\otimes K(B)  \otimes  K(B)^{\otimes F}\\&\xrightarrow{\id_{K(I)^{\otimes Z\setminus F'}}\otimes \partial \otimes \id_{K(B)^{\otimes F}}} &K(I)^{\otimes Z\setminus F'}\otimes \Sigma K(I)  \otimes  K(B)^{\otimes F}\\&\simeq&K(I)^{\otimes Z\setminus F'}\otimes \Sigma K(I)^{\otimes F'\setminus F}  \otimes  K(B)^{\otimes F}\\&\simeq&\Sigma K(I)^{\otimes Z\setminus F} \otimes  K(B)^{\otimes F}\ .
\end{eqnarray*}
Otherwise 
we have $d^{Z,1}_{j,F,F'}=0$.
\end{prop}
\begin{proof}
 In the case that $F\subseteq F'$ we used that $F'\setminus F$ is a one-point set and therefore
 $K(B)^{\otimes F'\setminus F}\simeq K(B)$ and $ \Sigma K(I)^{\otimes F'\setminus F}\simeq  \Sigma K(I)$ canonically.




We have 
a map of exact sequences given by the bold part of the following diagram.
$$\hspace{-2cm}
\xymatrix{
0\ar[r]&I^{\otimes Z\setminus F'}\otimes I^{F'}_{1}\ar@/_0.5cm/@{..>}[ddrr]^(0.3){\id_{I^{\otimes Z\setminus F'}}\otimes \pi^{F'}_{1} }\ar[d]\ar[r]&I^{\otimes Z\setminus F'}\otimes A^{\otimes F'}\ar[r]\ar[d]&I^{\otimes Z\setminus F'}\otimes B^{\otimes F'}\ar[r]\ar[d]&0\\
0\ar[r]&\ar@{..>}[d]^{\pi^{Z}_{j+1}}I^{Z}_{j+1}\ar[r]&I^{Z}_{j}\ar[r]&\bigoplus_{F\in \cP_{|Z|-j}(Z)} I^{\otimes Z\setminus F} \otimes B^{\otimes F}\ar[r]&0\\&\ar@{..>}[d]\bigoplus_{F'''\in \cP_{|Z|-j-1}(Z)} I^{\otimes Z\setminus F'''} \otimes B^{\otimes F'''}&& I^{\otimes Z\setminus F'}\otimes \bigoplus_{F''\in \cP_{|F'|-1}(F')} I^{\otimes F'\setminus F''}\otimes  B^{\otimes F''}\ar@{.>}[ll]^{}&\\&I^{\otimes Z\setminus F}\otimes B^{\otimes F}&&&}\ ,
$$
where the restriction of lower dotted horizontal arrow
to the summand $F''$ is given by the composition
$$ I^{\otimes Z\setminus F'}\otimes I^{\otimes F'\setminus F''}\otimes  B^{\otimes F''}\cong 
I^{\otimes Z\setminus F''} \otimes  B^{\otimes F''}\to \bigoplus_{F\in \cP_{|Z|-j-1}(Z)} I^{\otimes Z\setminus F} \otimes B^{\otimes F}$$
of the canonical identification and inclusion of a summand.
 Note that $|F''|=|F'| -1 =|Z|-j-1$.
 The lower left vertical map is the projection into the summand indexed by $F$.

  We read off  the following commutative square in $\Mod_{KU}(\Sp)$
\begin{equation}\label{fqewdqewdewdqd}\hspace{-1.5cm}\xymatrix{\ar@/^1cm/@{..>}[dddr]\ar@/^-3cm/[dd]_{\id_{K(I)^{\otimes Z\setminus F'}}\otimes d^{F',1}_{0}}K(I)^{\otimes Z\setminus F'}\otimes K(B)^{\otimes F'}\ar[r]^{ }\ar[d]^{\id_{K(I)^{\otimes Z\setminus F'}}\otimes \partial^{F'}_{0}} &\bigoplus_{F\in P_{|Z|-j}(Z)}K(I)^{\otimes Z\setminus F}\otimes K(B)^{\otimes F}\ar[dd]^{d_{j}^{Z,1}} \\ K(I)^{\otimes Z\setminus F'}\otimes \Sigma K(I^{ F'}_{1} ) \ar[d]^{\id_{K(I)^{\otimes Z\setminus F'}}\otimes \pi^{F'}_{1}}& \\ K( I)^{\otimes Z\setminus F'}\otimes \bigoplus_{F''\in P_{|F'|-1}(F')} \Sigma K(  I)^{\otimes F'\setminus F''}\otimes  K(B)^{\otimes F''}\ar[r]^-{ } &\ar[d]\bigoplus_{F'''\in P_{|Z|-j-1}(Z)}\Sigma K(I)^{\otimes Z\setminus F'''}\otimes K(B)^{\otimes F'''}\\ &\Sigma K(I)^{\otimes Z\setminus F}\otimes K(B)^{\otimes F}&}\ ,
\end{equation}  
where $d^{Z,1}_{F,F'}$ is induced by the dotted arrow.
This diagram expresses $d^{Z,1}_{F,F'}$ in terms of $d^{F',1}_{0}$. 
We therefore must understand $d^{F',1}_{0}$. To the end we consider the   map of exact sequences $$
\xymatrix{
0\ar[r]&I_{1}^{F'}\ar@/^-3cm/@{..>}[dd]_{\pi^{F'}_{F}}\ar[d]^{\pi^{F'}_{1}}\ar[r]&A^{\otimes F'}\ar[r]\ar[dd]&B^{\otimes F'}\ar[r]\ar@{=}[dd]&0\\&\bigoplus_{F''\in \cP_{|F'|-1}(F')}I^{F'\setminus F''}\otimes B^{F''}\ar[d]&&&\\
0\ar[r]&I^{F'\setminus F}\otimes B^{F}\ar[r]&A^{F'\setminus F}\otimes B^{F'}\ar[r]&B^{F'\setminus F}\otimes B^{F}\ar[r]&0}
$$
where $F$ is an element of $\cP_{|F'|-1}(F')$.  We read off the following commutative square in $\Mod_{KU}(\Sp)$
$$\hspace{-0.5cm}\xymatrix{K(B)^{\otimes F'} \ar[d]^{d^{F',1}_{0}} &\ar[l]_{\simeq} K(B)^{\otimes F'\setminus F}\otimes K(B)^{\otimes F}  \ar[r]^{\simeq}&\ar[d]^{\partial \otimes \id_{K(B)^{\otimes F}}} K(B) \otimes K(B)^{\otimes F}\\\bigoplus_{F''\in \cP_{|F'|-1}(F')}\Sigma K(I)^{\otimes F'\setminus F''}\otimes K(B)^{\otimes F''} \ar[r] &  \Sigma K(I)^{\otimes F'\setminus F}\otimes K(B)^{\otimes F} \ar[r]^{\simeq}& 
  \Sigma K(I )\otimes K(B)^{\otimes F}}$$
  which describes the components of $d^{F',1}_{0}$ in terms of $\partial$.
We  insert this calculation into \eqref{fqewdqewdewdqd} and read off the asserted formula for $d^{Z,1}_{j,F',F}$.
\end{proof}

\subsection{Phantoms}\label{rtijowhpwgrfrefwerfw}

One application of homotopy theory in the present paper  is to
extend the validity of a statement about $G$-$C^{*}$-algebras from a small
  initial set of examples  to a much bigger set.  For example,
 if the statement is that a natural transformation between functors is an equivalence,    and  the functors have  colimit preserving factorizations over $\KK^{G}$, then we can obviously extend to the colimit closure of  the set of examples in $\KK^{G}$.

This idea was applied in the proof of \cref{jeirgowregwerfrefw},
where the initial set of examples was the image of the induction functor $\Ind^{G}$.
As we shall see in \cref{qewfoiihqowedewqdeqwd}, the Rokhlin property for $A$ in $G\nCalg$  does not imply that $\kk^{G}(A)$ belongs to  the colimit closure
of the image of $\Ind^{G}$ itself, but it does imply a slightly weaker statement, namely being in the closure of the the image   of $\Ind^{G}$ under summable phantom retracts. In order to formalizes this
and eventually show \cref{gkerpogkpwergwerfwerf}, in the present section we recall some general $\infty$-categorical notions. The language developed here can be used in general  to capture asymptotic properties
of $C^{*}$-algebras and maps   homotopically. We will   indicate the main idea in \cref{jioregowegrewe9}.

Let $\bC$ be a cocomplete stable $\infty$-category.
\begin{ddd}
A  morphism $C\to D$ in $\bC$ is called compact\footnote{
We took this definition from a lecture series of Dustin Clausen at IHES, 2023 on Efimov's K-theory, where he attributed it  to Akhil Mathew.  Compact morphisms play an important role in the intrinsic characterization of dualizable objects in $\Pr^{L}_{\mathrm{st}}$.}
 if for any filtered  diagram $Z:I\to  \bC$
with $\colim_{I}Z\simeq 0$ and morphism  $D\to Z(i)$ for some $i$ in $I$ there exists a morphism $ i\to j$ in $I$ such that  the composition $C\to D\to Z(i) \to Z(j)$ vanishes.
\end{ddd}

\begin{rem}\label{rtkjglprtegergrgedf}
It might be useful to record the following equivalent characterization of the compactness of a morphism $C\to D$:
For every filtered diagram $Z:I\to \bC$ and map $D\to \colim_{I}Z$, there exists an $i$ in $I$ and a commutative square
\begin{equation}\label{}\xymatrix{C\ar[r]\ar[d] &Z(i) \ar[d] \\ D\ar[r] & \colim_{I}Z } \ .
\end{equation}
\hB
\end{rem}

Note that an object $C$ is compact    if $\id_{C}:C\to C$ is a compact morphism.
If at least one of $C$ or $D$ is a compact object, then every morphism $C\to D$ is compact.

\begin{ddd}\mbox{}
\begin{enumerate}\item 
A morphism $D\to E$ in $\bC$ is called a strong phantom map if the composition $C\to  D\to E$ vanishes for
every compact map $C\to D$. \item An object $D$ is called a strong phantom object if $\id_{D}$ is a strong phantom map.
\end{enumerate}
\end{ddd}

\begin{rem}
Classically, a morphism $D\to E$ is called a phantom map, if for every morphism $C\to D$ from a compact object $C$
the composition $C\to D\to E$ vanishes.  An object $D$ is called phantom if $\id_{D}$ is phantom. 
A strong phantom morphism is a phantom morphism, and a strong phantom object is a phantom object.
If $\bC$ is compactly generated, then the notions of a strong phantom morphism and phantom morphism  are actually equivalent.   Moreover, zero is the only phantom object. 

In the other extreme, if $0$ is the only compact object of $\bC$, then every morphism or object is phantom, but not necessarily strong phantom.  So the strong versions of phantom objects and morphisms are  interesting  replacements of the classical notions in the case of   cocomplete  stable $\infty$-categories that are not compactly generated. \hB
\end{rem}

We consider an adjunction $$L:\bC\leftrightarrows \bD:R$$  between cocomplete stable $\infty$-categories. The following is straightforward to prove.
\begin{lem}\label{kopgwefgrewfwerf} 
If $R$ preserves colimits, then it preserves strong phantom morphisms and $L$ preserves compact morphisms.
 \end{lem}
%
%
%
%


In our applications, we encounter colimit preserving functors that are not right adjoints. 
In order to control their interaction with the strong phantom  morphisms, we introduce the following stronger condition.

\begin{ddd}
A morphism $C\to D$ in $\bC$  is a universal (strong) phantom morphism, if
for every colimit preserving functor $F:\bC\to \bD$ the morphism $F(C)\to F(D)$ is also
a (strong) phantom morphism.
\end{ddd}

We will actually encounter a potentially even stronger condition.  
\begin{ddd}\mbox{}\begin{enumerate}
\item\label{erijogoewrpgfrefw} A sequence $(f_{i})_{i\in \nat}$ of morphisms $f_{i}:C\to D$ is called summable if 
there exists a morphism $\hat f:C\to \bigoplus_{\nat}D$ with the components $(f_{i})_{i\in \nat}$.
\item A morphism $f:C\to D$ is called a summable   map if the constant sequence $(f)_{i\in \nat}$ is summable. \end{enumerate}
\end{ddd}
In \ref{erijogoewrpgfrefw} we say that $\hat f$  witnesses  the summability of the sequence.
 
 \begin{rem} If $(f_{i})_{i\in \nat}$ is a summable sequence, then the choice of $\hat f$  fixes a choice of a sum $\sum_{i\in \nat} f_{i}:C\to D$, employing the universal property of $\bigoplus_{\nat}$ as the coproduct in $\bC$. This motivates the word {\em summable}. 
The sum may depend on the choice of  $\hat f$.   \hB \end{rem}

\begin{rem}
If the category $\bC$ is a dualizable presentable $\infty$-category, then one can show
that every strong phantom morphism between $\aleph_{1}$-compact objects is summable  \cite[Prop. 2.46]{budu}. As a consequence of   \cite[Thm. 1.1]{budu}
 we know that $\EE^{G}$ is dualizable. Note that
 $\EE^{G}_{\sepa}$ is precisely the subcategory of $\aleph_{1}$-compact
 objects in $\EE^{G} $. Therefore every strong phantom between objects in 
 $\EE^{G}_{\sepa}$ (i.e. $\EE^{G}$-classes of  separable $G$-$C^{*}$-algebras) is summable.
  \hB
\end{rem}

\begin{lem}\label{kwgjerovfdsbf}\mbox{}
\begin{enumerate}
\item A summable morphism is a  universal strong phantom morphism.
\item Countable colimit preserving functors preserve summable sequences.
\end{enumerate}
\end{lem}
\begin{proof}
Let $f:C\to D$ be a summable   morphism. We first show that $f$ is a strong phantom morphism.
Let $F\to C$ be a compact morphism. Then we can find (see \cref{rtkjglprtegergrgedf}) an integer $k$ and  the dotted extension of the diagram
to 
$$\xymatrix{F\ar@{-->}@/^1cm/[rr]^{0}\ar@{..>}[r]\ar@{..>}[dr] &C\ar[r]^{f}\ar[dr]^{\hat f}&D \\&\bigoplus_{i=0}^{k}D\ar@{..>}[r]&\bigoplus_{\nat}D\ar[u]_{\pr_{n}}}\ . $$
We then choose $n$ bigger than $k$ and get the dashed part of the diagram showing that the composition $F\to C\to D$ vanishes. This shows that $f$ is a strong phantom.

If $R$ is a countable colimit preserving functor and $(f_{i})_{i\in I}$ is a summable sequence of morphisms $f_{i}:C\to D$ with summability  witnessed by $\hat f$, then $R(\hat f)$ witnesses  $(R(f_{i}))_{i\in \nat}$ as  summable.
It then follows from the first part that a summable  map is  actually a universal strong phantom map.
 \end{proof}

Let $D$ be an object of $\bC$.
\begin{ddd}
An object $C$ is called a  (summable, universal, strong) phantom retract of $D$ if there exists
a map $C\to D$ such that $\Fib(C\to D)\to C$ is a (summable, universal, strong) phantom map.
\end{ddd}

A (strong) phantom object is a (strong) phantom retract of $0$.

\begin{ddd} If $\bU$ is a full subcategory of $\bC$, then by $\overline{\bU}$  we denote the closure of $\bU$ under taking universal  strong phantom retracts. \end{ddd}

For any colimit preserving functor $R:\bC\to \bD$ between cocomplete stable $\infty$-categories we have the inclusion
$R(\overline{\bU})\subseteq \overline{R(\bU)}$, where $R(\bU)$ denotes the full subcategory of $\bD$ on the image of $R$.
\begin{kor} \label{ioerjgoiwergreffwreferwfwr}Assume that $\bD$ is compactly generated. If $R$ is colimit preserving  functor  which   annihilates $\bU$,  then $R$ annihilates 
 $ \overline{\bU}$.
 \end{kor}

 Here is a general method to detect  summable   retracts.
 We assume that $\bC$ admits countable products.
 \begin{lem}\label{wekorgpwegfwerferwfwrefw}
If $C$ fits into a diagram of the form
 \begin{equation}\label{asdvjasvopascsdacadc}\xymatrix{\prod_{\nat}C& C\ar[l]_{\diag}\ar[d]^{d}&\\\bigoplus_{\nat}C\ar[r]\ar[u]^{\can}& P\ar[ul]\ar[r]^{q}&Q}
\end{equation}  with a horizontal fibre sequence,
 then $C$ is a  summable phantom retract of $Q$.
 \end{lem}
\begin{proof}
For any $n$ in $\nat$ we can
  extend the diagram as follows:
 \begin{equation}\label{idhuewhdiuewqdeqwdqwed}\xymatrix{ F\ar@{..>}[dr]^(0.3){\iota}\ar@/^-1cm/@{-->}[dd]&&\\\prod_{\nat}C \ar@/^1cm/[rr]^{\pr_{n}}& C\ar@{..>}[dr] \ar@{=}[r]\ar[l]_{\diag}\ar[d]^{d}  &C \\\bigoplus_{\nat}C\ar[r]\ar[u]^{\can}& P\ar[ul]\ar[ur]_(0.7){p_{n}}\ar[r]^{q}&Q}\ ,
\end{equation}  
where the lower diagonal dotted arrow is the composition $q\circ d$, $\iota$ represents the fibre of this map,
and the dashed arrow is obtained from the universal property of the  lower horizontal fibre sequence and the 
vanishing of $F\xrightarrow{d\circ \iota} P\xrightarrow{q} Q$. The dashed arrow witnesses $\iota$ as summable  map.
%
%
%
\end{proof}

One can use the diagram \eqref{asdvjasvopascsdacadc} also to detect  summable maps.

\begin{kor}\label{wkotgpwetgferfwfref}
If $f:D\to C$ is a map such that $D\to C\xrightarrow{q\circ d} Q$ vanishes, then
$f$ is a summable phantom map.
\end{kor}
\begin{proof}
The assumption implies that $f$ has a factorization $D\to F\xrightarrow{\iota}C$.
\end{proof}

\begin{ex}\label{jioregowegrewe9}
In this example,  we give a simple application of the language developed above.
Let $A,B$ be  $C^{*}$-algebras and $f,h:A\to B$ be two morphisms.
One says that $f$ and $h$ are unitarily isomorphic if there exists a unitary $u$  in $UM(B)$ such that $h(-)=u^{*}f(-)u$. If $f$ and $h$ are unitarily isomorphic, then $\kk (f)\simeq \kk (h)$ and therefore also $\ee (f)\simeq \ee (h)$.

More generally, we say that $f$ and $h$ are approximately  unitarily isomorphic 
if  there exists a sequence of unitaries $(u_{n})_{n\in \nat}$  in $UM(B)$ such that $h(a)=\lim_{n\in \nat}u_{n}^{*}f(a)u_{n}$ for every $a$ in $A$.

\begin{lem}
If  $f$ and $h$ are approximately  unitarily isomorphic, then $\ee (f)-\ee (h)$ is a summable phantom map in $\EE $. 
\end{lem}
\begin{proof} 
Let $(u_{n})_{n\in \nat}$ be a sequence in $UM(B)$ 
exhibiting the approximate unitary isomorphism between $f$ and $h$.
 We consider the exact sequence
\begin{equation}\label{vfoijewoivfvsdfvsdfvsdfv}  0\to \bigoplus_{\nat} B\to \prod_{\nat }B \to B^{\infty}\to 0\ .\end{equation}
The sequence $(u_{n})_{n\in \nat}$ provides unitaries $u$ in $ UM( \prod_{\nat }B )$ and 
$\bar u$ in $UM(B^{\infty})$.
We let $$\bar f:A\xrightarrow{f}B\xrightarrow{\diag} \prod_{\nat}B \to B^{\infty}$$ and define $\bar h$ similarly.
The assumptions imply that $\bar f(-)=\bar u^{*}h(-)\bar u$.

We form the diagram
 \begin{equation}\label{asdvjasvopascsdrrrrrracadc}\xymatrix{\prod_{\nat}\ee (B)& \ee (B)\ar[l]_{\diag}\ar[d]^{d}&\\\bigoplus_{\nat}\ee (B)\ar[r]\ar[u]^{\can}& \ee (\prod_{\nat}B)\ar[ul]\ar[r]^{q}&\ee (B^{\infty})}
\end{equation}
Since $\ee $ preserves countable sums the lower sequence is a fibre sequence.
 Let $\delta:=\ee (f)-\ee (h)$.
Then 
$q\circ d\circ \delta\simeq \ee (\bar f)-\ee (\bar h):\ee (A)\to  \ee (B^{\infty})$ vanishes since
$\bar f$ and $\bar h$ are unitarily isomorphic.
By \cref{wkotgpwetgferfwfref}, we conclude that $\delta$ is a  summable phantom map.
\end{proof}

%


\end{ex}

\subsection{The Rokhlin property}\label{qewfoiihqowedewqdeqwd}



 The Rokhlin property of $G$-$C^{*}$-algebra $A$ a is usually defined 
as an approximative condition. But using the techniques from \cref{jioregowegrewe9}, one can reformulate the condition as follows.

\begin{ddd}[{\cite[Def. 3.1]{Izumi_2004}}]\label{rktohoprthrgtgegetrg} $A$ has the   Rokhlin property if   there exists  a family 
$(p_{g})_{g\in G}$ of orthogonal projections in $M(A^{\infty})$ such that:
\begin{enumerate}
\item We have $\sum_{g\in G} p_{g}=1_{M(A^{\infty})}$.
\item  For  for all $g,h$ in $G$ we have $hp_{g}=p_{hg}$.
\item For every $g$ in $G$ and $a$ in $A$ we have
$p_{g}a'=a' p_{g}$, where $a'$ denotes the image of $a$ under 
 $A\stackrel{\diag}{\to} \prod_{\nat} A\to A^{\infty}$.   
 \end{enumerate}\end{ddd}
 Here, $A^{\infty}$ is defined as in \eqref{vfoijewoivfvsdfvsdfvsdfv}.
The first two conditions say that $(p_{g})_{g\in G}$ is an equivariant partition of unity. 
Note that the definition as stated above is a natural extension of Izumi's definition to non-unital algebras.

We let $A_{r}^{\infty}$ denote the  subalgebra of $A^{\infty}$ generated by the products $p_{g}a'$ for all $g$ in $G$ and $a$ in $A$.  We note that the projections $p_{g}$  belong to    $M(A_{r}^{\infty})$ and commute with all element of $A_{r}^{\infty}$.
We let $ A_{r}^{\infty,G}$ denote the subalgebra of $G$-fixed points in $ A_{r}^{\infty}$.
\begin{lem}  \label{oekgopwergerferwfwref} We have an isomorphism $A^{\infty}_{r}\cong \Ind^{G}(A_{r}^{\infty,G})$ in $G\nCalg$.
\end{lem}
\begin{proof}
We have an isomorphism of algebras $A_{r}^{\infty,G}\to p_{e}A_{r}^{ \infty}$, $a\mapsto p_{e}a$ with inverse
$a\mapsto \sum_{g\in G} p_{g}ga$.
We identify $\Ind^{G}(A_{r}^{\infty,G})\cong C(G,p_{e}A_{r}^{ \infty}  )$ and define a $G$-invariant isomorphism
$$ C(G, p_{e}A_{r}^{\infty})\to A_{r}^{\infty}\ , \quad f\mapsto \sum_{g\in G} p_{g} gf(g)\ .$$
Its inverse sends $a$ to the function $g\mapsto p_{e}g^{-1} a$. 
 \end{proof}

It follows from \cite[Prop. 3.60]{budu} that the adjunction \eqref{dfbvdfvdfvdfvsdfv}   descends through the comparison  \eqref{wrfqwdewdewfqef} between $KK$ and  $E$-theory and provides  an adjunction
$$\Ind_{H}^{G}: \EE^{H} \leftrightarrows  \EE^{G} : \Res_{H}^{G}\ .$$
We will use the case where $H$ is the trivial group.
Recall that $\overline{\Ind^{G}(\EE )}$ denotes the closure of the image of $\Ind^{G}$ under universal strong phantom retracts.

\begin{prop}\label{woketpgwergwergw9}
If $A$ has the Rokhlin property, then $\ee^{G} (A)$ belongs to $\overline{\Ind^{G}(\EE )}$.
\end{prop}
\begin{proof}
We retain the notation from \cref{oekgopwergerferwfwref}.
We define the $G$-$C^{*}$-algebra $P$ as the pull-back
$$\xymatrix{P\ar@{^{(}->}[r]^{i}\ar[d]^{q} & \prod_{\nat}A\ar[d] \\  A_{r}^{\infty}\ar@{^{(}->}[r] & A^{\infty}}\ . $$ It fits 
into a commutative diagram
$$\xymatrix{&\prod_{\nat}A& A\ar[l]_{\diag}\ar[d]^{d}&&\\0\ar[r]&\bigoplus_{\nat}A\ar[r]\ar[u]^{\can}& P\ar[ul]^{i}\ar[r]^{q}&A^{\infty}_{r}\ar[r]&0}\ ,$$
where the lower horizontal sequence is exact.
We now apply $\ee^{G} $ and use that it preserves countable sums and sends exact sequences to fibre sequences in order to get the commutative diagram
$$\xymatrix{\ee^{G} (\prod_{\nat}A)& \ee^{G} (A)\ar[l]_{\diag}\ar[d]^{d}&\\ \bigoplus_{\nat}\ee^{G} (A)\ar[r]\ar[u]^{\can}&\ee^{G}  (P)\ar[ul]^{i}\ar[r]^{q}&\ee^{G} (A^{\infty}_{r}) }\ ,$$
 where  lower horizontal sequence is  part of a fibre sequence.
  \cref{wekorgpwegfwerferwfwrefw} implies that $\ee^{G} (A)$  is a universal phantom retract of
  $\ee^{G} (A^{\infty}_{r}) $.  By \cref{oekgopwergerferwfwref}, the latter object belongs to 
  $\Ind^{G}(\EE )$.
%
\end{proof}

%
%

  \begin{prop}[\cref{gkerpogkpwergwerfwerf}]\label{gkerpogkpwergwerfwerf1}
If $G$ is finite and $A$ is a $G$-$C^{*}$-algebra with  the Rokhlin property, then the assembly map is an equivalence
$$\colim_{BG} K(A)\stackrel{\simeq}{\to}  K(A\rtimes_{r}G)\ .$$
\end{prop}
\begin{proof} 
%
%

We extend the natural transformation \eqref{frewfwrewoifjowrfwrefwerf} 
to a cofibre sequence  of functors with cofibre denoted by $\Cof$.
 We have equivalences  \begin{equation}\label{fqopwjeopfqwdqewdqe}\colim_{BG}\K(\widehat \Res^{G}(\kk^{G}(-)))\simeq \colim_{BG}K( -) \ , \quad \K(\kk^{G}(-)\rtimes_{r}G)\simeq K(-\rtimes_{r}G)
\end{equation}of functors from $G\nCalg$ to 
 $\Fun(BG,\Mod_{KU}(\Sp))$. These functors are homotopy invariant, $G$-stable, exact, $s$-finitary
 and preserve countable sums. By \eqref{bvwervopkpsdfvdfvsfdv} they factorize over
 $\ee^{G} $. We therefore also have a colimit preserving factorization $\Cof$ through colimit preserving functors 
 $$\KK^{G}\xrightarrow{c^{G}} \EE^{G}\xrightarrow{\Cof^{E}} \Mod_{KU}(\Sp)\ .$$
Since $\langle c^{G}(\langle \Ind^{G}(\KK) \rangle)\rangle = \langle \Ind^{G}(\EE )\rangle $,  
 it  follows from \cref{jeirgowregwerfrefw} that   $\Cof^{E}$ annihilates   the localizing subcategory generated by  $\Ind^{G}(\EE)$.  Since $ \Mod_{KU}(\Sp)$ is compactly generated, we can conclude by \cref{ioerjgoiwergreffwreferwfwr} that
$\Cof^{E}$ annihilates $\overline{\langle \Ind^{G}(\EE)\rangle}$. Thus $$\Cof(\kk^{G}(A))\simeq \Cof^{E}(\ee^{G}(A))\simeq 0$$ by \cref{woketpgwergwergw9} and   therefore 
$$\colim_{BG}   \K(\widehat \Res^{G}(\kk^{G}(A)))\stackrel{\simeq}{\to} \K(\kk^{G}(A)\rtimes_{r}G)\ . $$
In view of the equivalences \eqref{fqopwjeopfqwdqewdqe} this is equivalent to the assertion of  \cref{gkerpogkpwergwerfwerf}.
   \end{proof}

%

  \subsection{Projectivity }\label{oierjgoiegwergwerg} 

We assume that $G$ is discrete and consider $M$ in $ \Fun(BG,\Mod_{KU}(\Sp))$.
For example, $M$ could be the $K$-theory spectrum $K(A)$ with the $G$-action induced by functoriality of a $G$-$C^{*}$-algebra $A$. As we have seen above,  the calculation of the $K$-theory groups of crossed products is often reduced to the  calculation of 
  $\pi_{*}(\colim_{BG}M)$. We want to understand these homotopy groups in terms of
the homotopy groups $\pi_{*}(M)$ with the $G$-action induced by functoriality. 
In general, we have a spectral sequence with second term
$$E^{2}_{p,q}\cong H_{p}(G,\pi_{q}(M))$$ converging to $\pi_{*}(\colim_{BG}M)$. 
However, even for simple groups like $G=C_{p}$, this spectral sequence turns out to be too complicated to be useful this for explicit calculations. 

In this section, we discuss the condition of projectivity of $M$ which simplifies this calculation considerably. 
If $M$ is projective, then this spectral sequence degenerates since $E^{2}_{p,q}=0$ for $p\not=0$ and
the groups $\pi_{*}(\colim_{BG}M)$ have a simple algebraic description.

We then discuss conditions on $C^{*}$-algebras $A$ (or objects $A$ in $\KK^{G}$) that imply that $K(A)$ (or $\K(\widehat \Res^{G}(A))$) is projective.


\begin{ddd}\mbox{}
\begin{enumerate}
\item   $M$ is called free  if it is equivalent to  $\bigoplus_{G}N$  with the $G$-action by left translations on the index set for some $N$ in $\Mod_{KU}(\Sp)$.
\item $M$ is called projective if it is a retract of a free object.
\end{enumerate}
  \end{ddd}
  
  \begin{lem}\label{kojerwgpwergerfewrfref} \mbox{}
  \begin{enumerate}
\item  If $M$ is projective, then we have an isomorphism \begin{equation}\label{vsadoivjadsoivcadscsadc1}
\colim_{BG}\pi_{*}(M) \cong \pi_{*}(\colim_{BG} M)\ . \end{equation} 
\item If $G$ is finite and $M$ is projective, then we have  an  equivalence (or isomorphism, respectively)
 \begin{equation}\label{vsadoivjadsoivcadscsadc}   \colim_{BG}M\simeq \lim_{BG}M\ , \quad 
    \pi_{*}(\lim_{BG}M) \cong   \lim_{BG}\pi_{*}(M) \ .\end{equation}  
   \end{enumerate}
%
%
%
\end{lem}
\begin{proof} 
The assertions of this lemma are certainly well-known, but let us add the argument for completeness.
For any cocomplete $\infty$-category $\cC$, we have adjunctions
\begin{equation}\label{fdsvpoksfopvosfdvsfdvsdvsv}\bigoplus_{G}:\cC\leftrightarrows \Fun(BG,\cC) :\Res^{G}\ , \quad \colim_{BG}:\Fun(BG,\cC)\leftrightarrows \cC:\Res_{G}\ ,
\end{equation} 
where $\Res^{G}$ forgets the $G$-action and $\Res_{G}$ introduces the trivial $G$-action.
The second adjunction is the definition of the colimit, and the first is the precise description of the functor $\bigoplus_{G}$. We have an equivalence \begin{equation}\label{veervwercecsfv}\Res^{G}\circ \Res_{G}\simeq \id
\end{equation}  which induces the equivalence of adjoints
 \begin{equation}\label{verpokpwevefsfvfdv}
 \colim_{BG}\bigoplus_{G}\stackrel{\eqref{veervwercecsfv}}{\simeq}
  \colim_{BG}\bigoplus_{G}\Res^{G}\Res_{G}\xrightarrow{counits,\simeq} \id_{\cC}\ ,
\end{equation}
  where the arrow marked with $counits$ is given by the composition of the counits of the two adjunctions.  We apply the equivalence  \eqref{verpokpwevefsfvfdv} to the natural transformation \begin{equation}\label{fwerfwerfsfff}\colim_{BG}\pi_{*}\to \pi_{*}\colim_{BG}
\end{equation}   of functors from $\Fun(BG,\Mod_{KU}(\Sp))$ to $\Ab^{\Z}$, precomposed with $\bigoplus_{G}$, and get the commutative square
\begin{equation}\label{}\xymatrix{
\colim_{BG}\pi_{*}(\bigoplus_{G}(-))\ar[r]^{!}\ar[d]^{\cong}& \pi_{*}(\colim_{BG} \bigoplus_{G}(-))\ar[d]^{\cong}   \\\pi_{*}(-)\ar@{=}[r]&\pi_{*}(-)}\ .
\end{equation}
 This shows that the marked transformation is an isomorphism.
It provides the isomorphism \eqref{vsadoivjadsoivcadscsadc1}  for free objects $M$. 
Since 
  a retract of an isomorphism is again an isomorphism,
we conclude that the natural transformation \eqref{fwerfwerfsfff} more generally induces the isomorphism \eqref{vsadoivjadsoivcadscsadc1} for projective $M$.

If $\cC$ is also complete, then we have
adjunctions
\begin{equation}\label{fdsvpoksfopvosfdvsfdvsdvsv1} \Res^{G}:\cC\leftrightarrows \Fun(BG,\cC):\prod_{G}\ , \quad \Res_{G}:\cC\leftrightarrows \Fun(BG,\cC): \lim_{BG}\ .
\end{equation} 
The  analogue of \eqref{verpokpwevefsfvfdv}  is an equivalence \begin{equation}\label{foqwejofpewdewdqwedqewdqe} \id_{\cC} \xrightarrow{\simeq}   \lim_{BG}\prod_{G} \ .
\end{equation}
We now assume   that $\cC$ is additive. 
We have a  natural transformation $\id_{\cC}\to \Res^{G}\prod_{G}$ that sends $M$ to the map
$(f_{g})_{g\in G}:M\to \prod_{G}M$ given by  $$f_{g}:=\left\{\begin{array}{cc} \id_{M}&g=e\\ 0& g\not=e \end{array} \right.\ .$$  If $G$ 
 is finite, then the adjoint of this  transformation is an equivalence
 \begin{equation}\label{ewrogpwefrefewfr}\Nm^{\Res^{G}}:\bigoplus_{G}\xrightarrow{\simeq} \prod_{G}\ .
\end{equation}
To see this, we observe that after applying $\Res^{G}$ we obtain the canonical equivalence between a the coproduct and the product in an additive category, and we use that $\Res^{G}$ is conservative.

\begin{rem}
$(\Res^{G},\Nm^{\Res^{G}})$ is a tensor iso-normed functor in the sense of \cite[Def. 2.3.1]{zbMATH07526076}.\hB
\end{rem}

 The analog of the argument above  now shows that the natural transformation
$$\pi_{*}\lim_{BG} \to \lim_{BG}\pi_{*} $$ induces the second  isomorphism in \eqref{vsadoivjadsoivcadscsadc}
on   projective objects in $\Fun(BG,\Mod_{KU}(\Sp))$. 

In order to see the first equivalence in  \eqref{vsadoivjadsoivcadscsadc}, we
use the norm map
\begin{equation}\label{fdsvsdfvfdcfcqrfv}\Nm^{\Res_{G}}:\colim_{BG}\to \lim_{BG}
\end{equation}  
between functors from $\Fun(BG,\cC)$ to $\cC$
which again exists for finite $G$ and additive $\cC$.

We claim that it induces and equivalence when applied to free objects. 
In order see this we could cite \cite[Lemma 1.3.8]{MR3904731}. Alternatively, we could use \cite[4.2.4]{HL13} in order to see that the composition 
$$   \colim_{BG}\circ \bigoplus_{G} \xrightarrow{\Nm^{\Res_{G}}  \bigoplus_{G}}  \lim_{BG}\circ \bigoplus_{G}\xrightarrow{\lim_{BG} \Nm^{\Res^{G}}} \lim_{BG}\prod_{G}\ , $$  
  under   the equivalences \eqref{foqwejofpewdewdqwedqewdqe} and \eqref{verpokpwevefsfvfdv},
 is equivalent to the identity functor of $\cC$.  
 Since  $\Nm^{\Res^{G}}$ from \eqref{ewrogpwefrefewfr} is known to be an equivalence,
 we see that $\Nm^{\Res_{G}}  \bigoplus_{G}$
 is an equivalence, too.
    The claim implies that the norm is also an equivalence on projective objects.
   %
%
%
%
\end{proof}

\begin{rem}
If $\cC$ is stable,  then the norm map \eqref{fdsvsdfvfdcfcqrfv} is a natural transformation between finite limit preserving functors. It follows that the equivalence 
 $ \colim_{BG} M \stackrel{\simeq}{\to} \lim_{BG}M$
 persists  to compact objects $M$ in $\Fun(BG,\cC)$, i.e. objects which can presented as finite colimits of projectives.  \hB
\end{rem}

%

We assume that $G$ is finite  and consider $M$ in $\Fun(BG,\Mod_{KU}(\Sp))$.
\begin{lem}\label{ekorgpefreqwdqedqed} If $|G|$ acts as an equivalence on $M$, then $M$ is projective. \end{lem}
\begin{proof}
We have the counit and the unit maps  $$\bigoplus_{G}\Res^{G}M\to M\ , \quad M\to \prod_{G}\Res^{G}M$$ of the first adjunctions in \eqref{fdsvpoksfopvosfdvsfdvsdvsv} and \eqref{fdsvpoksfopvosfdvsfdvsdvsv1}, respectively.
One can check that the composition
\begin{equation}\label{bsdfsvfdvdfsvsdfv}M\to  \prod_{G}\Res^{G}M\stackrel{\eqref{ewrogpwefrefewfr}}{\simeq} \bigoplus_{G}\Res^{G}M\to M
\end{equation} 
is given by multiplication by $|G|$. Since $|G|$ acts as an equivalence on $M$, the composition  \eqref{bsdfsvfdvdfsvsdfv} exhibits $M$ as a retract of the free object $\bigoplus_{G}\Res^{G}M$.
\end{proof}

\begin{rem}
We note that \cref{qirjfofdewdewdqewd} and \cref{ekorgpefreqwdqedqed} can be considered as instances of
the higher Maschke Theorem as stated in \cite[Thm. 2.4.4]{zbMATH07526076}.
\end{rem}

\begin{ex} Assume that $G$ is finite and consider the UHF algebra  \begin{equation}\label{fwmfpqefdqewdqewde}M_{|G|^{\infty}}:=\colim_{n\in \nat} \bigotimes_{i=0}^{n}\Mat_{|G|}(\C)
\end{equation} in $\nCalg$. 
  For any $G$-$C^{*}$-algebra $A$ we have   an equivalence 
$$\kk^{G}(A\otimes M_{|G|^{\infty}})\simeq \kk^{G}(A)[|G|^{-1}]\ .$$
If $A$ is a $G$-$C^{*}$-algebra  which tensorially absorbs $M_{|G|^{\infty}}$, i.e.,   
$A\cong A\otimes M_{|G|^{\infty}}$, then  also $\kk^{G}(A) \simeq \kk^{G}(A)[|G|^{-1}]$.

\begin{kor}
If $A$  tensorially absorbs $M_{|G|^{\infty}}$, then    $K(A) $ satisfies the assumption of  \cref{ekorgpefreqwdqedqed}   and is therefore projective in $\Fun(BG,\Mod_{KU}(\Sp))$. \end{kor}

Assume now that $A$ has the trivial $G$-action and tensorially absorbs $M_{|G|^{\infty}}$.
Then, on the one hand, we have
$$\pi_{*}(\colim_{BG} K(A))\cong  \colim_{BG} K_{*}( A) \cong K_{*}(A)\ .$$
On the other hand, since $K(C_{r}^{*}(G))\simeq  \bigoplus_{|\hat G|}KU$ and $A\rtimes_{r} G\cong A\otimes C_{r}^{*}(G)$ (see \eqref{fwerfeiuhieuworf} for the $KK$-theory level version),
we have
$$ K_{*}(A\rtimes_{r} G)\cong \bigoplus_{|\hat G |} K_{*}(A)\ .$$
If $G$ is non-trivial, then the map
$$\pi_{*}(\colim_{BG} K(A))\to K_{*}(A\rtimes_{r} G)$$ is not an isomorphism.
This shows that   condition
\cref{qirjfofdewdewdqewd}.\ref{qirjfofdewdewdqewdf} can not be replaced by the weaker  condition
that $|G|$ acts as an equivalence.
%
%
%
%
%
%
%
%
%
%
 \hB
  \end{ex}

\begin{rem}\label{weotkgperfgrefrefwrefwf}
Assume that $A$ is  a unital Kirchberg algebra that absorbs  $M_{|G|^{\infty}}$ and assume that
$G\to \KK_{0}(A,A)^{-1}$ (the superscript $-1$ indicates invertible elements) is a homomorphism whose images fix the class of the identity. Then by \cite[Rem. 2.12]{MR3572256},  this homomorphism 
 comes from a refinement of $A$ to a $G$-$C^{*}$-algebra with the Rokhlin property.   
 By \cref{ekorgpefreqwdqedqed}, we know that $ K( A)$ is  a projective object in $\Fun(BG,\Mod(KU))$. 
  We can therefore
 deduce  from
   \cref{gkerpogkpwergwerfwerf} and \eqref{vsadoivjadsoivcadscsadc1} that \begin{equation}\label{qeroijofeweqdewdq1}\colim_{BG} K(A)\simeq K(A\rtimes_{r}G)\ , \qquad \colim_{BG} K_{*}( A)\simeq K_{*}(A\rtimes_{r}G)\ .
\end{equation}
   \hB
\end{rem}

 The following is consequence of   \cref{woketpgwergwergw9}.
\begin{kor}\label{rtlkhpwtrgrtgtge} Assume that $G$ is finite.
If $A$ has the Rokhlin property, then $K( A)$ is a universal strong phantom retract of a free object.
\end{kor}
\begin{proof}
We use that the functor $\K^{E}\circ \widehat \Res^{G}$ preserves colimits and therefore universal strong phantom retracts. It furthermore 
sends the objects of $\Ind^{G}(\EE )$ to free objects.
By   \cref{woketpgwergwergw9}, we conclude that $ \K^{E} ( \widehat \Res^{G}(\ee^{G} (A)))\simeq  K(A)$ is a universal strong phantom retract of a free object.
 \end{proof}
 
 Note that \cref{rtlkhpwtrgrtgtge} does not imply that the object $K(A)$ is  a projective object
in  $\Fun(BG,\Mod_{KU}(\Sp))$ since   in this category  not every 
 universal strong phantom retract  is a retract. 
   But we can obtain projectivity under additional assumptions, different from the property that $|G|$ acts as an equivalence.

Let $A$ be in $G\nCalg$.
 \begin{prop}\label{wetkogpwefrwefrwgrg}
  Assume that $G$ is finite.    If $A$ is simple, has the Rokhlin property,   and 
   $K(  A )$  is compact  $\Fun(BG,\Mod_{KU}(\Sp))$, 
 then $K( A )$ is projective.  
 \end{prop}
 \begin{proof}
By Lemma \eqref{oekgopwergerferwfwref} and the observations made in \cref{wtkohpwerfrefwref}, we have an equivalence 
$K( A^{\infty}_{r})\simeq \bigoplus_{G} K( A^{\infty,G}_{r})$  in $\Fun(BG,\Mod_{KU}(\Sp))$, where the $G$-action on the sum is given by left-translations 
on the index set. 
The idea of the proof of \cref{wetkogpwefrwefrwgrg} is to present $K( A)$ as a retract of the free object $K( A^{\infty}_{r})$.

By the proof of \cref{woketpgwergwergw9}, we have a fibre sequence $$F\xrightarrow{\iota} \ee^{G}(A)\to  \ee^{G}(A_{r}^{\infty})$$ in $\EE^{G} $, where $\iota$ is a universal strong phantom morphism.
Since $\K^{E}\circ \widehat \Res^{G}$  preserves colimits and therefore preserves universal strong phantom morphism\footnote{It is also a right-adjoint and therefore preserves strong phantoms.}, we conclude that
$\K^{E}(\widehat\Res^{G}(\iota)):\K^{E}(\widehat \Res^{G}(F))\to K( A )$ is a universal strong phantom morphism.

By assumption,  $K( A)$ is a compact object on $\Fun(BG,\Mod_{KU}(\Sp))$.
We use simpleness of $A$ to the effect that  the non-equivariant homomorphism $A\to A^{\infty}_{r}$, $a\mapsto p_{e}a'
$ is injective and thus identifies $A$ with the corner $p_{e}A^{\infty}_{r}$. In view of the proof of \cref{oekgopwergerferwfwref}, we conclude that  $K( A^{\infty}_{r})\simeq \bigoplus_{G}\Res^{G}(K(A))$ is also compact in  
$\Fun(BG,\Mod_{KU}(\Sp))$. This implies that the left entry of the fibre sequence
\begin{equation}\label{fqwefewfewdqed}K^{E} (\widehat\Res^{G}(F))\stackrel{\K^{E}(\widehat\Res^{G}(\iota))}{ \to}   K( A) \to K( A^{\infty}_{r})
\end{equation}
 is also compact. But then   $\K^{E}(\widehat\Res^{G}(\iota))\simeq 0$ 
since it is a strong phantom morphism from a compact object.
We conclude that the fibre sequence \eqref{fqwefewfewdqed} splits and
$$K( A^{\infty}_{r})\simeq K( A )\oplus \Sigma \K^{E} (\widehat\Res^{G}(F))$$ in $\Fun(BG,\Mod_{KU}(\Sp))$.
\end{proof}

\begin{rem}
In the statement of  \cref{wetkogpwefrwefrwgrg}, we could replace the simpleness of $A$ with the assumption that $A\to A^{\infty}_{r}p_{e}$ is injective, or simply with the assumption that $K(A^{\infty}_{r})$ is a compact object in 
$\Fun(BG,\Mod_{KU}(\Sp))$. But note that these assumptions do not only depend on $A$ but also on the choice of the additional
data given by the Rokhlin property.

Note that the compactness of 
 $K( A )$  in $\Fun(BG,\Mod_{KU}(\Sp))$ is a strong assumption, much stronger than compactness of $\Res^{G}(K  (A)) $ in $\Mod_{KU}(\Sp)$.
\hB
\end{rem}


Combining \cref{wetkogpwefrwefrwgrg}  with  \cref{gkerpogkpwergwerfwerf},   
\eqref{vsadoivjadsoivcadscsadc1} and  \eqref{vsadoivjadsoivcadscsadc} we conclude:
\begin{kor}\label{wtjkopegfwerferfwrefrefw}
 Under the assumptions of    \cref{wetkogpwefrwefrwgrg},   we have   isomorphisms
\begin{equation}\label{frepofkpqwedewdqewdq}\colim_{BG} K_{*}(A)\cong K_{*}(A\rtimes_{r}G)\cong \lim_{BG}K_{*}(A)\ .
\end{equation}  
\end{kor}

\begin{rem}In this remark we explain how    \cref{wtjkopegfwerferfwrefrefw} is related with the results of \cite{Izumi_2004}.
Let $A$ be a $G$-$C^{*}$-algebra for a finite group $G$ and let $A^{G}$ be the subalgebra of $G$-fixed points.
By \cite[Prop. 4.3]{Rieffel_1980} we have   an isomorphism  of $C^{*}$-algebras
\begin{equation}\label{aslvjaodcadscadscadsc}(A\otimes K(L^{2}(G))^{G}\cong A\rtimes_{r}G\ .
\end{equation} 
 By \cite[Lemma. 3.11]{Izumi_2004}, if $A$ is unital and has the  Rokhlin property,
then   $A$ is stable, i.e., we have an isomorphism $A\cong A\otimes K(L^{2}(G))$.
We conclude  under these assumptions that 
\begin{equation}\label{inter0}A^{G}\cong A\rtimes_{r}G\ .\end{equation}
If $A$ is in addition simple, then \cite[Thm. 3.13]{Izumi_2004} states that the inclusion $A^{G}\to A$ induces an isomorphism
\begin{equation}\label{inter1}K_{*}(A^{G})\cong \lim_{BG} K_{*}(A)\ . \end{equation}
 If we combine  \eqref{inter0}, \eqref{inter1}, then we can conclude that if $A$ is simple, unital, and satisfies   the Rokhlin property, then
\begin{equation}\label{inter3}K_{*}(A\rtimes_{r}G)\cong \lim_{BG} K_{*}(A)\end{equation}
%
%
%

 We thus get the second isomorphism in the statement of \cref{wtjkopegfwerferfwrefrefw}.
 
 Turning things around,  up to the  additional compactness assumption  on $K(A)$,
\cref{wtjkopegfwerferfwrefrefw} provides a new proof of   \cite[Thm. 3.13]{Izumi_2004}, which seems to be  very different from Izumi's argument.
%
\hB
 \end{rem}

%
%


 \section{The K-theory functor on the orbit category}

 \subsection{A KK-valued functor on the orbit category}
 
 Let $G$ be a second countable locally compact group.
 A $G$-orbit is a  locally compact topological $G$-space  
 that is isomorphic to a $G$-space of the form $G/H$ for some closed subgroup $H$ of $G$. 
 We let  $G\Orb$   denote the $\infty$-category presented by the full topologically enriched subcategory of $G\Top$ on $G$-orbits. If $G$ is discrete, then $G\Orb$ is an ordinary category.
 
We let   $G\Orb^{\comp}$ or $G\Orb^{\fin}$  denote the full subcategories of $G\Orb$ on $G$-orbits which are compact or finite, respectively. For a family of subgroups $\cF$, i.e.
a set of closed subgroups of $G$ which is closed under conjugation and taking subgroups, we let $G_{\cF}\Orb$ denote the full subcategory of $G\Orb$ on $G$-orbits isomorphic to $G/H$ with $H$ in $\cF$.
 
 In this section, we define a functor $$V^{G}:G\Orb^{\comp}\to \KK^{G}$$
 that will be used as an ingredient of the definition of the functor $\K^{G}_{A}$ in \cref{wtgijowergrwefwrefwef} below.

  We let  $G\LCH_{+}$   denote the category of locally compact Hausdorff spaces with continuous $G$-action   and
 partially defined proper maps.  We have a functor $$C_{0}:G\LCH_{+}^{\op}\to G\nCalg$$ that sends a locally compact Hausdorff $G$-space $X$  to the $G$-$C^{*}$-algebra $C_{0}(X)$ of complex functions on $X$ vanishing at $\infty$.
By Gelfand duality,  it identifies $G\LCH_{+}^{\op}$ with the full subcategory of $ G\nCalg$ of commutative $G$-$C^{*}$-algebras.

 
 The category $ G\LCH_{+}$ has a topological enrichment such that for any compact Hausdorff space
 $W$ we have
 $$\Hom_{\Top}(W,\underline{\Hom}_{G\LCH_{+}}(X,Y))\cong\Hom_{G\LCH_{+}}(W\times X,Y)\ ,$$
 where $\underline{\Hom}_{G\LCH_{+}}  (X,Y)$ denotes the topological refinement of the set $ \Hom_{G\LCH_{+}}(X,Y)$.
 The enrichment provides a notion of homotopy equivalence.  We  
 let $$L_{h} : G\LCH_{+}  \to G\LCH_{+,h}$$ denote the Dwyer-Kan localization of $G\LCH_{+}  $ at the homotopy equivalences. 
By the homotopy invariance of the functor $\kk^{G}$ the composition 
  $$ G\LCH_{+}^{\op}\xrightarrow{C_{0}}G\nCalg\xrightarrow{\kk^{G}}\KK^{G}$$ 
  sends homotopy equivalences to equivalences. The universal property of the Dwyer-Kan localization $L_{h}$  then provides a factorization of this composition through a functor 
   indicated by the dotted arrow in 
     \begin{equation}\label{vrevjeoivfevsdfvdfvsdv}\xymatrix{G\LCH_{+}^{\op}\ar[r]^-{C_{0}}\ar[dr]_{L_{h}}&G\nCalg\ar[r]^-{\kk^{G}}& \KK^{G} \\&G\LCH_{+,h}^{\op}\ar@{..>}[ur]_{\kk^{G}C_{0}(-)}& }\ .\end{equation} 
 

     We   let $$ \underline{\KK^{G}}:\KK^{G,\op}\times \KK^{G}\to \KK^{G}$$ denote the internal $\Hom$-functor of the closed symmetric monoidal $\infty$-category $\KK^{G}$. 
 We can  restrict the functor $\kk^{G}C_{0}(-)$  along the  fully faithful  inclusion $  G\Orb^{\comp}\to  G\LCH_{+,h}$. 
  \begin{ddd}\label{wijotghwrthggbdh}
We define the functor \begin{equation}\label{vfsw3fsfsfg}V^{G}(-):= \underline{\KK^{G}}(\kk^{G}C_{0}(-),\beins_{\KK^{G}}):G\Orb^{\comp}  \to \KK^{G}\ .
\end{equation} \end{ddd}

\begin{rem}\label{9qriuzghiqergrfqrff}
If $G$ is not compact, then not every  map between $G$-orbits is   proper.
But the restriction to   compact orbits in \cref{wijotghwrthggbdh}  ensures this properness. 
The other option would be to restrict along the inclusion $G_{\mathcal{C}\mathrm{omp}} \Orb\to G\LCH_{+,h}$, where $\mathcal{C}\mathrm{omp}$ denotes the family of compact subgroups.
In the case of discrete groups the resulting functor 
$$ \underline{\KK^{G}}(\kk^{G}C_{0}(-),\beins_{\KK^{G}}):G_{\mathcal{C}\mathrm{omp}} \Orb \to \KKG$$ turns out to be equivalent, as a consequence of Paschke duality, to the functor
denoted by $\kk^{G}_{\nCcat}\circ \C[-]$ in \cite[Sec. 16]{bel-paschke}. 
The latter is actually defined on the whole orbit category $G\Orb$.
For discrete $G$ and $A$ in $\KK^{G}$, the functor  
$$K^{G}_{A}:=\K^{G}(A\otimes \kk^{G}_{\nCcat}\circ \C[-]):G\Orb\to \Mod(KU)$$ is the functor
featuring the Baum-Connes conjecture in \cref{wergjiwoeferfw}.

If $G$ is not compact, then $G_{\mathcal{C}\mathrm{omp}} \Orb\cap G\Orb^{\comp}=\emptyset$.
Hence, we can not use the results from \cite[Sec. 16]{bel-paschke} in order to extend the functor $V^{G}$.
\hB
\end{rem}

If $S$ in $G\Orb^{\comp}$ is such that $\kk^{G}C_{0}(S)$ is dualizable in $\KK^{G}$, then by definition
$V^{G}(S)$ is the dual of $\kk^{G}C_{0}(S)$. The following two propositions refine this observation.

Note the restriction to finite orbits in the following statement.
 \begin{prop}\label{wetijgowpgferfwerferwfw}
For every $S$ in $G\Orb^{\fin}$, the object $\kk^{G}C_{0}(S)$ of $\KK^{G}$ is selfdual. \end{prop}
 \begin{proof}
 We consider 
  $S=G/H
$ in $G\Orb^{\fin}$. Note that $C_{0}(S)\cong \Ind_{H}^{G}(\Res_{H}(\C))$, $\beins_{\KK^{H}}\simeq \kk^{H}(\Res_{H}(\C))$, and therefore 
$\kk^{G}C_{0}(S)\simeq \Ind_{H}^{G}(\beins_{\KK^{H}})$. Using this observation, the selfduality of
$\kk^{G}C_{0}(S)$ is a formal consequence of the some of the  formulas listed in \cref{wtkopgwfrefrefw}.
We   have the following   chain of equivalences which is natural in $A,B$ in $\KK^{G}$:
 \begin{eqnarray*}
 \KK^{G}(B, \kk^{G}C_{0}(S)\otimes A) 
&\simeq&\KK^{G}(B, \Ind_{H}^{G}(\beins_{\KK^{H}})\otimes A)\\
&\stackrel{\eqref{dfbvdfvdfvdfvsdfv2}}{\simeq}& \KK^{G}(B, \Ind_{H}^{G}( \Res^{G}_{H}(A)))\\
&\stackrel{\eqref{dfbvdfvdfvdfvsdfv1}}{\simeq}& \KK^{H}(\Res^{G}_{H}(B),  \Res^{G}_{H}(A))\\
&\stackrel{\eqref{dfbvdfvdfvdfvsdfv}}{\simeq}& \KK^{G}(\Ind_{H}^{G}(\Res^{G}_{H}(B)),  A )\\
&\stackrel{\eqref{dfbvdfvdfvdfvsdfv2}}{\simeq}& \KK^{G}(\Ind_{H}^{G}(\beins_{\KK^{H}})\otimes  B ,  A )\\
& \simeq&  \KK^{G}( \kk^{G}C_{0}(S)\otimes  B ,  A ) 
\end{eqnarray*}
Note that we must assume that $S$ is finite since this is required for   the adjunction  \eqref{dfbvdfvdfvdfvsdfv1} used in the chain above. This condition also implies that $H$ is open which is necessary for the application of 
\eqref{dfbvdfvdfvdfvsdfv}.\end{proof}

\begin{rem}
By \cref{wekogpwefgrefwerfwefew}, we have a tensor normed functor 
$(\Res^{G}_{H},\Nm^{G}_{H})$.  We could conclude the self-duality of  $\Ind_{H}^{G}(\Res_{H}^{G}(\beins_{\KK}))\simeq \kk^{G}C_{0}(G/H)$ from  \cite[2.3.4, (3)]{zbMATH07526076} (which has essentially the same proof as given above).  We write out   the details since they will be used    in the proof   \cref{qerogfpqrfewqdewdq}
asserting an   equivariant refinement of this self-duality.
\hB 
\end{rem}

 By functoriality, we can consider    $V^{G}(S)$ and $\kk^{G}C_{0}(S)$ both  as  objects in $\Fun(B\Aut_{G\Orb}(S),\KK^{G})$.

\begin{prop} \label{qerogfpqrfewqdewdq} For every $S$ in $G\Orb^{\fin}$, we have a canonical equivalence of functors 
\begin{equation}\label{dsfvdsfvdfvsdfv}V^{G}(S)\simeq \kk^{G}C_{0}(S): B\Aut_{G\Orb}(S)\to \KK^{G}\ .
\end{equation}  
\end{prop}
\begin{proof}
Note that $V^{G}(S) $ is by  \eqref{vfsw3fsfsfg} the dual to $\kk^{G}C_{0}(S)$. If one ignores the  $B\Aut_{G\Orb}(S)$-actions, then 
such an equivalence \eqref{dsfvdsfvdfvsdfv} is given by \cref{wetijgowpgferfwerferwfw}.
Our task is to refine the selfduality of $\kk^{G}C_{0}(S)$ to a $B\Aut_{G\Orb}(S)$-equivariant selfduality.
The given proof of \cref{wetijgowpgferfwerferwfw} does not make this equivariance apparent.

  The selfduality of $\kk^{G}C_{0}(S)$ is equivalent to an adjunction
\begin{equation}\label{sdfvdfsvr3fevfvsdfv}-\otimes \kk^{G}C_{0}(S):\KK^{G}\leftrightarrows  \KK^{G}: -\otimes \kk^{G}C_{0}(S)\ .
\end{equation} 
Such an adjunction  is fixed by  the unit  which is a natural transformation 
$$\epsilon:\id_{\KK^{G}}\to  (-\otimes  \kk^{G}C(S))\otimes \kk^{G}C(S)$$
of endofunctors of $\KK^{G}$. In the present case, this unit is given
  by the tensor product with   a morphism
$$\beins_{\KK^{G}}\to  \kk^{G}C_{0}(S)\otimes  \kk^{G}C_{0}(S)$$ in $\KK^{G}$.
Equivariance of the adjunction is then equivalent to a refinement of this  morphism     to a morphism in $\Fun(B\Aut_{G\Orb}(S),\KK^{G})$.

Below we will show that the morphism given by the proof of  \cref{wetijgowpgferfwerferwfw} 
is induced by
 the map of $C^{*}$-algebras
\begin{equation}\label{vfsdvsijowrjivewvfvsdfvsfdvs}\C\to C_{0}(S\times S)\cong C_{0}(S)\otimes C_{0}(S)\ , \quad  \lambda\mapsto \lambda \chi_{\diag(S)}\ ,
\end{equation} 
where $ \chi_{\diag(S)}$ is the characteristic function of the diagonal in $S$.
This map  is obviously  $B\Aut_{G\Orb}(S)$-equivariant.
 

We assume that $S=G/H$. 
We already know that $\kk^{G}C_{0}(G/H)$ is self dual by \cref{wetijgowpgferfwerferwfw}. The unit of this selfduality is the image of $\id_{\kk^{G}C_{0}(G/H)}$ under the duality 
$$\KK^{G}(\kk^{G}C_{0}(G/H),\kk^{G}C_{0}(G/H))\simeq \KK^{G}(\beins_{\KK^{G}},\kk^{G}C_{0}(G/H)\otimes \kk^{G}C_{0}(G/H))$$
given by the chain \begin{eqnarray*}
\KK^{G}( \kk^{G}C_{0}(G/H),\kk^{G}C_{0}(G/H))&\simeq&
\KK^{G}( \Ind_{H}^{G}(\beins_{\KK^{H}}),\kk^{G}C_{0}(G/H))\\&\simeq&
\KK^{H}(  \beins_{\KK^{H}},\Res^{G}_{H}(\kk^{G}C_{0}(G/H)))\\&\simeq&
\KK^{H}( \Res^{G}_{H} (\beins_{\KK^{G}}),\Res^{G}_{H}(\kk^{G}C_{0}(G/H)))\\&\simeq&
\KK^{G}( \beins_{\KK^{G}}, \Ind_{H}^{G}(\Res^{G}_{H}(\kk^{G}C_{0}(G/H))))\\&\simeq&
\KK^{G}( \beins_{\KK^{G}}, \kk^{G}C_{0}(G/H) \otimes \kk^{G}C_{0}(G/H))
\end{eqnarray*} of equivalences in the proof of \cref{wetijgowpgferfwerferwfw}.
We must show that this image of $\id_{\kk^{G}C_{0}(G/H)}$ is induced by the map
which sends $\lambda$ in $\C$ to $\lambda \chi_{\diag(G/H)}$.
 As explained in \cref{wtkopgwfrefrefw}, all these  equivalences are induced by explicit maps on the level of morphims sets of $C^{*}$-algebras:
\begin{eqnarray*}
\Hom_{G\nCalg}(  C_{0}(G/H), C_{0}(G/H))&\cong&
\Hom_{G\nCalg}( \Ind_{H}^{G}(\Res_{H}(\C)), C_{0}(G/H)) \\&\stackrel{1}{\to}&
\Hom_{H\nCalg}(  \Res_{H}(\C) ,\Res^{G}_{H}( C_{0}(G/H)))\\&\cong&
\Hom_{H\nCalg}( \Res^{G}_{H}(  \Res_{G}(\C)) ,\Res^{G}_{H}( C_{0}(G/H)))\\&\stackrel{2}{\to}&
\Hom_{G\nCalg}(    \Res_{G}(\C) , \Ind_{H}^{G}(\Res^{G}_{H}( C_{0}(G/H))))\\&\cong&
\Hom_{G\nCalg}(    \Res_{G}(\C),  C_{0}(G/H \times G/H))
\end{eqnarray*}
The map marked by $1$ is given by applying $\Res^{G}_{H}$ and precomposing
with the unit $\Res_{H}(\C)\to \Res^{G}_{H}(\Ind^{G}_{H}(\Res_{H}(\C)))$ of the  $(\Ind_{H}^{G},\Res_{H}^{G})$-adjunction which happens to exist on the level of algebras. This
unit sends $\lambda$ to the function $\lambda \chi_{\{H\}}$.
So the image of $\id_{C_{0}(G/H)}$ in the target of $1$ is the map 
$\lambda\mapsto \lambda \chi_{\{H\}}$. The map marked by $2$ is
given by applying $\Ind_{H}^{G}$ and then precomposing with the unit $\Res_{G}(\C)\to \Ind_{H}^{G}(\Res_{H}^{G}(\Res_{G}(\C)))$ of the
$(\Res_{H}^{G},\Ind_{H}^{G})$-adjunction which again happens to exist on the level of
algebras. This map sends  $\lambda$ to the constant function $\lambda \chi_{G/H}$.
The image of $\id_{C_{0}(G/H)}$ in the target of this map is therefore
the map $\lambda\mapsto(gH\mapsto (g'H\mapsto  \lambda \chi_{\{H\}}(g'H)))$.
The isomorphism $\Ind_{H}^{G}(\Res^{G}_{H}( C_{0}(G/H)))\cong C_{0}(G/H)\times C_{0}(G/H)$
sends the function $gH\mapsto (g'H\mapsto f(gH)(g'H))$ to the function
$\hat f:(gH,g'H)\mapsto f(gH)(g^{-1}g'H)$. In particular, the function $(gH,g'H)\mapsto \chi_{\{H\}}(g^{-1}g'H)$
is the characteristic function $\chi_{\diag(G/H)}$.
\end{proof}

 If $G$ is compact, then since  the left-adjoint   of the adjunction \eqref{bsdfvkspdfvqreve} has a symmetric monoidal refinement, we obtain a   lax symmetric monoidal refinement of the right-adjoint crossed product functor $-\rtimes_{r}G$. 
 The image of the group $C^{*}$-algebra $C_{r}^{*}(G)\cong \C\rtimes_{r}G$ in $\KK$ therefore refines to  a commutative algebra object 
  \begin{equation}\label{erwfeqewdweqwd}\tR(G):=\beins_{\KK^{G}}\rtimes_{r} G
\end{equation} in $ \KK$, and 
the reduced crossed product functor itself   refines to a functor
 \begin{equation}\label{erwfeqewdweqwd1}-\rtimes_{r}G:\KK^{G}\to \Mod_{\tR(G)}(\KK)\ .\end{equation}

If   $G$ is finite, then  $G$ with the action by left translations belongs to in $G\Orb^{\fin}$ and we have an isomorphism $\Aut_{G\Orb}(G)\cong G^{\op}$ with $G^{\op}$   acting on $G$ by multiplication from the right.  
By functoriality,
  $V^{G}(G)$ becomes an object of $\Fun(BG^{\op},\KK^{G})$.
We define the functor \begin{equation}\label{sdvwrfgeerfr}\widetilde \Res^{G}:\KK^{G}\to \Fun(BG^{\op},\Mod_{\tR(G)}(\KK)) \ ,\quad  \widetilde \Res^{G}(-):=(-\otimes V^{G}(G))\rtimes_{r}G\ .
\end{equation} 
 
The following result shows that $\widetilde \Res^{G}$ refines the restriction functor $\widehat \Res^{G}:\KK^{G}\to \Fun(BG,\KK)$ from \eqref{qfqwefqwedqwdq}.
Let \begin{equation}\label{zhtherthrthe} e:\Fun(BG^{\op},\Mod_{\tR(G)}(\KK))\to \Fun(BG,\KK)
\end{equation} be the functor which forgets the  $\tR(G)$-module structure
and precomposes with the inversion map $BG\to BG^{\op}$.
\begin{lem}\label{gwegerfwerfwref} In the case of a finite group $G$,
 we have an equivalence of functors $$e\circ \widetilde \Res^{G}\simeq \widehat \Res^{G}:\KK^{G}\to \Fun(BG,\KK)\ .$$
\end{lem}
\begin{proof}
For a $G$-$C^{*}$-algebra $A$, we have a  natural isomorphism in $G\nCalg$
$$A\otimes C_{0}(G)\stackrel{ \cong}{\to} \Res^{G}(A) \otimes C_{0}(G) \ .$$
Identifying elements of the tensor product with functions $f:G\to A$, this isomorphism is given  by
$f\mapsto (g\mapsto g^{-1}f(g))$. Under this isomorphism the $G^{\op}$-action   on the left goes to the action $(h,f)\mapsto (g\mapsto h^{-1} f(gh^{-1}))$ on the right. 
On the level of $\Fun(BG^{\op},\KK^{G})$,  this isomorphism induces an equivalence of functors
\begin{equation}\label{bfsbspokvopsrwgfergvfs}-\otimes C_{0}(G)\simeq e^{',*} \widehat \Res^{G}(-)\otimes C_{0}(G)\ ,
\end{equation} 
where $e^{\prime}:BG^{\op}\to BG$ is induced by the inversion map. 

We have the following chain of equivalences in $\Fun(BG,\KK)$: \begin{eqnarray*}
e\circ \widetilde \Res^{G}(-)&\simeq &e(-\otimes V^{G}(G))\rtimes_{r}G\\&\stackrel{\eqref{dsfvdsfvdfvsdfv}}{\simeq}&
e (-\otimes \kk^{G} C_{0}(G))\rtimes_{r}G\\&\stackrel{\eqref{bfsbspokvopsrwgfergvfs}}{\simeq} &
\widehat \Res^{G}(-)\otimes e(\kk^{G} C_{0}(G)\rtimes_{r}G)\\&\stackrel{\cref{wtkorhopweferwfwregw}}{\simeq}&
\widehat \Res^{G}(-)\ .
\end{eqnarray*}
\end{proof}

 \subsection{The K-theory functor}\label{rthkpwgfrefref}
 
 We let   \begin{equation}\label{rfweferfwrfwfwrf}\K:\KK\to \Mod_{KU}(\Sp)\ , \quad B\mapsto \map_{\KK}(\beins_{\KK}, B) \end{equation}
denote the lax symmetric monoidal functor corepresented by the tensor unit $\beins_{\KK}$ of $\KK$.
The usual $KU$-module valued $K$-theory functor   for $C^{*}$-algebras is then given by the composition
\begin{equation}\label{rfweferfwrfwfwrf1}K:\nCalg\xrightarrow{\kk} \KK\xrightarrow{\K}  \Mod_{KU}(\Sp)\end{equation} is   (see  \cite{blp}  or \cite{Bunke:2023aa} for details).
Similarly, for     a second countable locally compact group $G$,
we define
 the  lax symmetric monoidal equivariant $K$-theory functor  \begin{equation}\label{werfwerfreffwefrefwrefwref}
\K^{G}:\KK^{G}\to \Mod_{R(G)}(\Sp)\ , \quad B\mapsto \map_{\KK^{G}}( \beins_{\KK^{G}} ,B)\ .
 \end{equation}
  For a compact group $G$, 
 using $\beins_{\KK^{G}}\simeq \Res_{G}(\beins_{\KK})$,    we obtain the equivalence $$\K(\tR(G))\stackrel{ \eqref{erwfeqewdweqwd}}{\simeq}  \map_{\KK}(\beins_{\KK}, \beins_{\KK^{G}}\rtimes_{r}G)  \stackrel{\eqref{bsdfvkspdfvqreve}  }{\simeq} \map_{\KK^{G}}(\beins_{\KK^{G}} ,\beins_{\KK^{G}} )\stackrel{\eqref{lekjkqrlfqwedqewdqde}}{\simeq}R(G)$$ in $\CAlg(\Mod_{KU}(\Sp))$.
 This leads to a refinement of  the functor \eqref{rfweferfwrfwfwrf} to a    functor $$\K^{G}: \KK^{G}\to \Mod_{R(G)}(\Sp)$$ between module categories.
 In view of the    Green-Julg theorem \eqref{bsdfvkspdfvqreve}  and \eqref{erwfeqewdweqwd1} we furthermore get   the equivalence  of functors  \begin{equation}\label{kejtgprtgwegwregw}\K^{G}(-)\simeq \K\circ (-\rtimes_{r}G):  \KK^{G}\to \Mod_{R(G)}(\Sp)\ .
\end{equation}    
In this note the $R(G)$-module structures on the values of $K$-theory are an important detail 
since we want to consider localizations and completions at elements of $\pi_{0}R(G)$.

Let $G$ be a second countable locally compact group and recall \cref{wijotghwrthggbdh} of the functor $V^{G}$.
\begin{ddd} \label{wtgijowergrwefwrefwef}We define the functor
\begin{equation}\label{getgoijif0ewrfrwef}\K_{-}^{G}:\KK^{G}\to \Fun(G\Orb^{\comp},\Mod_{R(G)}(\Sp))\ , \quad A\mapsto (S\mapsto   \K^{G}(A\otimes V^{G}(S)))\ .
\end{equation} 
\end{ddd}

%
%

 Since   $V^{G}(*)\simeq \beins_{\KK^{G}}$ we get the following consequence of \eqref{kejtgprtgwegwregw}.   \begin{kor} \label{fkjqnrhiofqwefewfqd1}
If $G$ is compact, then we have a canonical equivalence of functors
$$\K^{G}_{-}(*)\simeq \K(-\rtimes_{r} G):\KK^{G}\to \Mod_{R(G)}(\Sp)\ .$$
 \end{kor}

The following corollary says that for a finite group $G$ and $A$ in $\KK^{G}$  the value of the functor $\K^{G}_{A}$ on the $G$-orbit $G$ with the induced $\Aut_{G\Orb^{\fin}}(G)$-action is equivalent, under the natural identifications, with $\K(A)$ with the $G$-action via functoriality.  Similarly as in \eqref{zhtherthrthe}, let \begin{equation}\label{} e:\Fun(BG^{\op},\Mod_{R(G)}(\Sp))\to \Fun(BG, \Mod_{KU}(\Sp))
\end{equation} be the functor which forgets the  $R(G)$-module structure
and precomposes with the inversion map $BG\to BG^{\op}$.

\begin{kor}\label{rtkohpertgregertgetg}In the case of a finite group $G$,
we have an equivalence of functors $$e(\K^{G}_{-}(G))\simeq \K(\widehat\Res^{G}(-)):\KK^{G}\to
\Fun(BG,\Mod_{KU}(\Sp))\ .$$
\end{kor}
\begin{proof} 
We apply $\K$ to the equivalence asserted in \cref{gwegerfwerfwref}. Unfolding definitions this gives precisely the assertion.
\end{proof}

Recall \cref{wekogpwefgrewfwre} of the bootstrap class $\UCT$.
Since
$\K(\beins_{\KK})\simeq KU$ and $\K$ preserves colimits, the restriction of the monoidal structure of $\K$ induces  an equivalence of functors  \begin{equation}\label{erthoijoeprthgetgtrgetrgetg}\K(-)\otimes_{KU}\K(-)\stackrel{\simeq}{\to} \K(-\otimes -):\UCT\times \KK\to \Mod_{KU}(\Sp)\ .
\end{equation}  
\begin{lem} \label{qrejigoqrgregwergwerg}Assume that  $G$ is a compact group,   $B$ is in $\KK$ and $A$ is in $\KK^{G}$, and that one of the following conditions is satisfied:
\begin{enumerate}
\item $B$ is in $\UCT$
\item \label{rkgofpewgefrfw} $\Res^{G}_{L}(A)\rtimes L$ is in $\UCT$ for  all closed subgroups $L$ of $G$.
\end{enumerate}
Then we have an equivalence of functors 
$$\K^{G}_{\Res_{G}(B)\otimes A}\simeq \K(B)\otimes_{KU}\K^{G}_{A} :G\Orb \to \Mod_{R(G)}(\Sp)\ .$$
\end{lem}
\begin{proof} The desired equivalence is given by the following chain: \begin{eqnarray*}
\K^{G}_{\Res_{G}(B)\otimes A}&\stackrel{\eqref{getgoijif0ewrfrwef}}{\simeq}&
\K^{G}(\Res_{G}(B)\otimes A\otimes V^{G})\\&\stackrel{\eqref{kejtgprtgwegwregw}}{\simeq}&
\K((\Res_{G}(B)\otimes A\otimes V^{G})\rtimes_{r} G)\\&\stackrel{\eqref{fwerfeiuhieuworf}}{\simeq}&
\K( B \otimes ((A\otimes V^{G})\rtimes_{r} G))\\& \stackrel{\eqref{erthoijoeprthgetgtrgetrgetg}}{\simeq}&
\K(B)\otimes_{KU} \K((A\otimes V^{G})\rtimes_{r} G)\\&\stackrel{\eqref{kejtgprtgwegwregw}}{\simeq} &
\K(B)\otimes_{KU} \K^{G}(A\otimes V^{G})\\
&\stackrel{\eqref{getgoijif0ewrfrwef}}{\simeq}&
\K(B)\otimes_{KU} \K^{G}_{A}.
\end{eqnarray*}
 The equivalence marked by \eqref{erthoijoeprthgetgtrgetrgetg}  is clear if $B$ is in $\UCT$. 
 In the other case we use that 
  $(A\otimes V^{G}(G/L))\rtimes_{r}G\simeq \Res^{G}_{L}(A) \rtimes_{r} G$ for all closed subgroups $L$.
 \end{proof}

\subsection{Group change}\label{wiegjhowgwegewg9}

In this subsection, we study the compatibility of $V^{G}$ and $\K^{G}$ with induction.
This will later be used in induction arguments on families of subgroups of $G$.
Let $G$ be a second countable locally compact group. 
 Let $L$ be a closed subgroup of $G$ such that $G/L$ is  compact, and  let $i:L\Orb^{\comp}\to G\Orb^{\comp}$ be the induction functor that is given on objects by $S\mapsto G\times_{L}S$. 
 
 \begin{lem}\label{weijogpwergwergf9}
 If $G/L$ is finite,
then we have a canonical equivalence of functors 
$$i^{*}V^{G}\simeq \Ind_{L}^{G}(V^{L}) : L\Orb^{\comp}\to \KK^{G}\ .$$
 \end{lem}
 \begin{proof} We start with providing
  a natural isomorphism of functors  
 \begin{equation}\label{fqoijfiofjp14orleqfwqf}C_{0}(G\times_{L}-)\cong \Ind_{L}^{G}(C_{0}(-))
\end{equation}on $L \Orb^{\comp}$ with values in  $G\nCalg$. 
Here we only need compactness of $G/L$.
 
 Let $S$ be in $L\Orb^{\comp}$. 
Since $G/L$ is compact the support condition mentioned in the text after \eqref{asdcadscaacd} becomes irrelevant
and we have  an isomorphism $\Ind_{L}^{G}(C_{0}(S))\cong C(G,C_{0}(S))^{L}$. Since $S$ is compact we can identify $C(G,C_{0}(S))\cong C(G\times S)$. Under this isomorphism 
  the $L$-invariance condition requires that $f(gl,l^{-1}s)=f(g,s)$ for all $l$ in $L$. But the functions satisfying this condition
  are precisely the functions which descend to the quotient $G\times  S\to G\times_{L}S$. Since
  $G\times_{L}S$ is compact we have the equality
  $C_{0}(G\times_{L}S)=C(G\times_{L}S)$.
  This gives the desired isomorphism \eqref{fqoijfiofjp14orleqfwqf}
  
  We now  get the asserted equivalence from 
   \begin{eqnarray*}
i^{*}V^{G}&\stackrel{\eqref{vfsw3fsfsfg}}{\simeq}& \underline{\KK}^{G}(\kk^{G}C_{0}(G\times_{L}-),\beins_{\KK^{G}})\\&\stackrel{\eqref{fqoijfiofjp14orleqfwqf}}{\simeq}&\underline{\KK}^{G}( \Ind_{L}^{G}(\kk^{L}C_{0}(-)),\beins_{\KK^{G}})\\&\stackrel{!}{\simeq} & \Ind_{L}^{G}(\underline{\KK}^{L}(  \kk^{L}C_{0}(-)),\beins_{\KK^{L}})\\&\stackrel{\eqref{vfsw3fsfsfg}}{\simeq}&
\Ind_{L}^{G}(V^{L})\ ,
\end{eqnarray*}
where the marked equivalence uses the equivalence
$$\underline{\KK^{G}}(\Ind_{L}^{G}(-_{1}),-_{2})\simeq \Ind_{L}^{G} (\underline{\KK^{L}}(-_{1},\Res^{G}_{L}(-_{2}))$$ which is a formal consequence of \eqref{dfbvdfvdfvdfvsdfv} (here we need that $L$ is closed and open in $G$, or equivalently,   that $G/L$ is discrete), \eqref{dfbvdfvdfvdfvsdfv2}, the equivalence $\beins_{\KK^{L}}\simeq \Res^{G}_{L}(\beins_{\KK^{G}})$, and 
the defining properties of the internal $\Hom$-functors.
  \end{proof}

We retain the assumption on $L$  and  calculate the restriction of $\K^{G}_{-}$ along $i$
 in terms of $\K^{L}_{-}$.
 
 Note that $\K^{G}_{A}$ for $A$ in $\KK^{G}$ takes values in $R(G)$-modules while $\K^{L}_{B}$ for $B$ in $\KK^{L}$ takes values in $R(L)$-modules. However, since
  $\beins_{\KK^{L}}\simeq \Res^{G}_{L}(\beins_{\KK^{G}})$, the action of $\Res^{G}_{L}$ on the endmorphism algebras of the tensor units yields
  a morphism $\Res^{G}_{L}:R(G)\to R(L)$ in $\CAlg(\Mod_{KU}(\Sp))$.   We will use the corresponding restriction morphism $$r^{*}:\Mod_{R(L)}(\Sp)\to \Mod_{R(G)}(\Sp)$$ in order to turn the values of $\K^{L}_{B}$ into $R(G)$-modules.   
 
 \begin{rem}
  If $G$ is compact, then the induced map
 $\pi_{0}\Res^{G}_{L}:\pi_{0}R(G)\to \pi_{0}R(L)$ sends the class $[V]$ of a finite dimensional representation $V$ of $G$ to the class $[V_{|L}]$ of its restriction   to $L$, see \cref{wrjtohipwgwregfwerfw}. \hB
 \end{rem}

 \begin{lem}\label{weijogpwergwergf91t}
 If $G/L$ is finite, then  an equivalence of functors $$i^{*}\K^{G}_{-}\simeq r^{*} \K^{L}_{\Res^{G}_{L}(-)}: \KK^{G}\times L\Orb^{\comp} \to \Mod_{R(G)}(\Sp)\ .$$   
 \end{lem}
\begin{proof}
The desired equivalence is given by the following chain:
\begin{eqnarray*}
i^{*}\K^{G}_{-}&\stackrel{\eqref{getgoijif0ewrfrwef}}{\simeq}& \K^{G}( -\otimes i^{*}V^{G} )\\ 
&\stackrel{\cref{weijogpwergwergf9}}{\simeq}
&\K^{G}(-\otimes \Ind_{L}^{G}(V^{L}) )\\
&\stackrel{\eqref{dfbvdfvdfvdfvsdfv2}}{\simeq}&
\K^{G}(\Ind_{L}^{G}(V^{L} \otimes \Res^{G}_{L}(-))\\&\stackrel{\eqref{werfwerfreffwefrefwrefwref}}{\simeq}&
\map_{\KK^{G}}(\beins_{\KK^{G}},  \Ind_{L}^{G}(V^{L} \otimes \Res^{G}_{L}(-)))\\
&\stackrel{\eqref{dfbvdfvdfvdfvsdfv1}}{\simeq}&
r^{*} \map_{\KK^{L}}(\beins_{\KK^{L}}, V^{L} \otimes \Res^{G}_{L}(-)  )\\
&\stackrel{\eqref{werfwerfreffwefrefwrefwref}}{\simeq}& r^{*} \K^{L}_{ \Res^{G}_{L}(-) }\ .
\end{eqnarray*}
Note that we must assume that $G/L$ is finite  since this is required for the adjunction
\eqref{dfbvdfvdfvdfvsdfv1} and \cref{weijogpwergwergf9} used in the chain above.
\end{proof}

If we combine  \cref{weijogpwergwergf91t} and \cref{fkjqnrhiofqwefewfqd1}, then we can calculate some values of the functor $\K^{G}_{A}$.
 \begin{kor}\label{qerjigoiewrgwerfrefw} If $G$ is compact and 
 $G/L$ is finite, then  for every $A$ in $\KK^{G}$ we have an equivalence $$\K^{G}_{A}(G/L)\simeq  r^{*}\K(\Res^{G}_{L}(A)\rtimes_{r} L)$$ of $R(G)$-modules.
 \end{kor}

%

We consider a closed  subgroup $L$ of $G$ and $B$ in $\KK^{L}$.  
\begin{prop}\label{ijrogergferwg9} If $G$ is finite, then we have an equivalence of functors
$$\K^{G}_{\Ind_{L}^{G}(B)}\simeq i_{!}\K^{L}_{B} :G\Orb \to \Mod_{KU}(\Sp)\ .$$
\end{prop}\begin{proof}
We have the following chain of equivalences:
\begin{eqnarray*}
\K^{G}_{\Ind_{L}^{G}(B)}&\stackrel{\eqref{getgoijif0ewrfrwef}}{\simeq}& \K^{G}(\Ind_{L}^{G}(B)\otimes V^{G})\\
&\stackrel{\eqref{dfbvdfvdfvdfvsdfv2}}{\simeq}&
\K^{G}(\Ind_{L}^{G}(B\otimes \Res^{G}_{L}(V^{G})))\\&\stackrel{\eqref{dfbvdfvdfvdfvsdfv1}}{\simeq}&
\K^{L}( B\otimes \Res^{G}_{L}(V^{G})) \ .
\end{eqnarray*}
We can apply  \eqref{dfbvdfvdfvdfvsdfv1} since $G$ is finite.

Let $j^{G}:G\Orb\to G\Set^{\fin}$ be the inclusion of the orbit category into the category of finite $G$-sets and equivariant maps.   We consider the functor
  $$\hat V^{G}:G\Set^{\fin}\to \KK^{G} \ , \quad S\mapsto  \underline{\KK}^{G}(\kk^{G}C_{0}(-),\beins_{\KK^{G}})\ .$$
Then we have a canonical equivalence $j^{G,*}\hat V^{G}\simeq V^{G}$, i.e., the functor $\hat V^{G}$ is an extension of $V^{G}$ from transitive $G$-sets   to all finite $G$-sets.

We have an adjunction
$$\hat i:L\Set^{\fin} \leftrightarrows G\Set^{\fin}: \res\ ,$$ where $\res$ restricts the $G$-action on an $L$-action, and
$\hat i(S):=G\times_{L}S$ is the extension of $i$ to finite $G$-sets. This implies an equivalence $\hat i_{!}\simeq \res^{*}$.

 The isomorphism $\Res^{G}_{L} (C_{0}(-))\cong j^{G,*} \res^{*} C_{0}(-)$ of functors from $G\Orb^{\fin}$ to $L\nCalg$ induces an equivalence  $ \Res^{G}_{L}(V^{G})\simeq  j^{G,*} \res^{*}\hat  V^{L}$, and consequently
\begin{eqnarray*}
\K^{L}( B\otimes \Res^{G}_{L}(V^{G}))&\simeq&  \K^{L}(B\otimes  j^{G,*} \res^{*}\hat  V^{L})\\&\simeq&  j^{G,*} \res^{*} \K^{L}(B\otimes  \hat  V^{L})\\&\simeq& j^{G,*}\hat i_{!} \K^{L}(B\otimes  \hat  V^{L})\\&\simeq&   i_{!}j^{L,*} \K^{L}(B\otimes  \hat  V^{L}) \\&\simeq &
i_{!}K^{L}_{B}
\ .
\end{eqnarray*}
 Here we have used that  the canonical map $i_{!}j^{L,*}F\to  j^{G,*}\hat i_{!}F$ is an equivalence  
 for any coproduct preserving functor $F$ on $G\Set^{\fin}$.
 \end{proof}

   \begin{prop}[\cref{werjoigpwegerfwref}]\label{werjoigpwegerfwref1} If $G$ is finite and  $A$ is in $\langle \Ind_{\prp}^{G} \rangle$, then  the assembly map is an equivalence
   \begin{equation}\label{vsdfvsdcscdsf111}\colim_{G_{\prp}\Orb} \K^{G}_{A}\xrightarrow{\simeq} \K(A\rtimes_{r}G)\ .
\end{equation}
  \end{prop}

 \begin{proof}  Since the assembly map is the evaluation at $A$ of a natural transformation between colimit preserving functors  on $\KK^{G}$, it suffices to verify the equivalence \eqref{vsdfvsdcscdsf} for
 $A=\Ind_{L}^{G}(B)$ for any proper subgroup $L$ of $G$ and $B$ in $\KK^{L}$.
  We use the factorization
 \begin{equation}\label{qefewdedqewqf}i:L\Orb\xrightarrow{i'} G_{\prp}\Orb\stackrel{i^{\prp}}{\to} G\Orb\ .
\end{equation}   The assembly map in the assertion of \cref{werjoigpwegerfwref} factorizes as the following chain of equivalences:
 \begin{eqnarray*}
\colim_{G_{\prp}\Orb} \K^{G}_{\Ind_{L}^{G}(B)}&\simeq&
(i^{\prp}_{!} \K^{G}_{\Ind_{L}^{G}(B)})(*)\\&\stackrel{\cref{ijrogergferwg9}}{\simeq}&
 (i^{\prp}_{!} i'_{!} \K^{L}_{ B})(*)\\&\stackrel{\eqref{qefewdedqewqf}}{\simeq}&(i_{!} \K^{L}_{ B})(*)\\&\stackrel{\cref{ijrogergferwg9}}{\simeq}&
 \K^{G}_{\Ind_{L}^{G}(B)}(*)\\&\stackrel{\cref{fkjqnrhiofqwefewfqd1}}{\simeq}& \K(B\rtimes_{r}G)\ ,
\end{eqnarray*}
where the first equivalence is an instance of the pointwise formula for the left Kan-extension $i^{\prp}_{!}$
and the identification of $G_{\prp}\Orb   $ with the slice of $i^{\prp}$ over $*$ .
  \end{proof}

We consider a finite group $G$. Recall the functor $$\Cof^{G}:\KK^{G}\to \Mod_{R(G)}(\Sp)$$ from \eqref{vwvssfdvsfvs}. If $H$ is a subgroup of $G$, then we obtain the faithful isofibration of finite
connected groupids $G/H\curvearrowleft G\to *\curvearrowleft G$. We then have, by \cref{iorgoergfgsfg}, a multiplicative induction functor  $$(-)^{\otimes  G/H}_{\kk}:\KK^{H}\to \KK^{G}\ .$$
Note that this induction functor is not necessarily exact, but preserves filtered colimits.

\begin{prop}\label{ewgojwoepgfrefwfwerfwrf}
The composition $$\Cof^{G}\circ  (-)^{\otimes  G/H}_{\kk}:\KK^{H}\to \Mod_{R(G)}(\Sp)$$
preserves all colimits.
\end{prop} 
\begin{proof}
The composition preserves filtered colimits since $\Cof^{G}$ and $ (-)^{\otimes  G/H}_{\kk}$ do so.
In particular, it sends zero algebras to zero objects. It remains to show
  exactness. It  suffices to show that the composition
$$ P: \KK^{H}_{\sepa}\xrightarrow{ (-)^{\otimes  G/H}_{\kk_{\sepa}}} \KK^{G}_{\sepa}\xrightarrow{y^{G}} \KK^{G}\xrightarrow{\Cof^{G}}  \Mod_{R(G)}(\Sp)$$ is exact. To this end it suffices to show that it sends all semisplit  
 exact sequences $0\to A\to B \stackrel{\pi}{\to} C\to 0$ in $H\nCalg_{\sepa}$ to  fibre sequences.
To this end we must show that the canonical map
$$P(A)\to \Fib(P(B)\to P(C))$$ is an equivalence. We will employ the filtration  $(I_{m})_{m=0,\dots,|G/H|}$ of $B^{\otimes G/H}$
from Step  \ref{ejigweporfrefrwfr} of the proof of \cref{weogjrpwegfrfwrfwre}.
Since $$0\to I_{1}\to B^{\otimes G/H} \to C^{\otimes G/H} \to 0$$ is exact and semisplit,  we have a fibre sequence
 $$\kk_{\sepa}^{G}(I_{1})\to \kk^{G}_{\sepa}(B^{\otimes G/H})\to \kk^{G}_{\sepa}(C^{\otimes G/H})\ .$$ 
  Applying  $\Cof^{G}\circ y^{G}$, we see that we must show that the map 
  $ P(A)\to \Cof^{G}(\kk^{G}(I_{1}))$ is an equivalence.  
  Since $ A^{\otimes G/H}\cong I_{|G/H|}$, the map
  $P(A)\to   \Cof^{G}(\kk^{G}(I_{|G/H|}))$ is an equivalence.
  We now argue by downward induction on $m$ that
  $ P(A)\to \Cof^{G}(\kk^{G}(I_{m}))$  is an equivalence for all $m$ in $\{1,\dots,|G/H|\}$. The case of $m=1$ is the desired assertion.  For $1\le m< |G|$ the map in question has a factorization 
  $$ P(A)\to \Cof^{G}(\kk^{G}(I_{m+1}))\to   \Cof^{G}(\kk^{G}(I_{m}))\ .$$
  The first map is an equivalence by the induction hypothesis. In order to show that the second map
  is an equivalence we
  use the semisplit exact  exact sequence  \eqref{qfjoipqwfevwfoqvfvvfqwdqwdqd} which
  induces a fibre sequence 
   $$\Cof^{G}(\kk^{G}(I_{m+1}))\to \Cof^{G}(\kk^{G}(I_{m}))\to\Cof^{G}(\kk^{G}(I_{m}/I_{m+1}))$$
  Using the notation from \eqref{vsdfvoijiowerjovervdfsvfv} we observe that \begin{equation}\label{qefwhiuqewuhiqewdewdqd}I_{m}/I_{m+1}\cong  \bigoplus_{[F]\in \cP_{m}(Z)/G} \Ind_{\sepa,G_{F}}^{G}(A^{\otimes F}\otimes C^{\otimes Z\setminus F})\ .
\end{equation}
  Since $1\le m <|G/H|$, all stabilizer subgroups $G_{F}$ are proper subgroups of $G$. Thus
  $\kk^{G}(I_{m}/I_{m+1})$ belongs to  $\langle \Ind_{\prp}^{G} \rangle$. By \cref{werjoigpwegerfwref}, we have
  $ \Cof^{G}(\kk^{G}(I_{m}/I_{m+1}))\simeq 0$ which completes the induction step.
  \end{proof}

%
%
%
%

 \section{KK-valued homology theories}   
  \subsection{Localization and completion}\label{tirwjgowgerwgwerg9}
 
 In this section, we recall some generalities on completion and localizations.   
  
 Let $\bC$ be any presentably symmetric monoidal   stable $\infty$-category. Our examples for $\bC$  are   $\Sp$, $\KK^{G}$, or $\Mod_{R(G)}(\Sp)$.  We write $\beins_{\bC}$ for the tensor unit of $\bC$. Then $R:=\map_{\bC}(\beins_{\bC},\beins_{\bC})
$ is a commutative algebra in $\Sp$.

 Let $\bD$ be a $\Mod_{R}(\Sp)$-module in  $\Pr_{\st}^{L}$. The basic example is $\bC$ itself. 
In the case of spectra $\Sp$, the tensor unit $\beins_{\Sp}$ and $R$ are both the sphere spectrum $S$ and  
$\Sp\simeq \Mod_{S}(\Sp)$.

 If we  apply this to the case $\bC=\KK^{G}$, then  $R\simeq R(G)$. 
 If $H \to G$ is a homomorphism, then we can consider
 $\KK^{H}$ as a $\Mod_{R(G)}(\Sp)$-module in place of $\bD$ using the restriction map $R(G)\to R(H)$.

We fix  a class  $\xi$   in $\pi_{0}R$.
We write $\xi:\id_{\bD}\to \id_{\bD}$ also for the  endomorphism of the identity functor $\id_{\bD}$   whose component at the object $D$ is 
given by $$D\simeq \beins_{\bC}\otimes D\xrightarrow{\xi\otimes \id_{D}}  \beins_{\bC}\otimes D\simeq D\ ,$$
where $\otimes$ denotes   the   structure map  $ \bC\times \bD \to \bD$ of the module structure and the two equivalences are given  the unit constraint.


 \begin{ddd}\label{weijgowiergijrefrewf9} \mbox{}\begin{enumerate} \item \label{regijwerogerfwrefrwefw}We define the $\xi$-torsion  functor by 
$$S_{\xi}:=\colim_{n\in \nat }\Fib(\id_{\bD}\stackrel{\xi^{n}}{\to} \id_{\bD}):\bD\to \bD\ .$$ 
 \item \label{regijwerogerfwrefrwefw11} We define the   inversion of $\xi$ by
$$(-)[\xi^{-1}]:=\Cofib(S_{\xi}\to \id_{\bD})\simeq \colim_{n\in \nat} (\id_{\bD}\xrightarrow{\xi} \id_{\bD}\xrightarrow{\xi} \id_{\bD}\xrightarrow{\xi} \dots):\bD\to \bD\ .$$ \end{enumerate}
 \end{ddd}
 By construction, we have a fibre sequence of endofunctors \begin{equation}\label{dvsdvfdvfdvsfvsdfvdfsvdfsverererer}S_{\xi}\to \id_{\bD}\to (-) [\xi^{-1}]
\end{equation}
 of $\bD$.

 The commutative algebra $R$, considered as an  $R$-module, is an object of $\Mod_{R}(\Sp)$ so that we can consider $S_{\xi}(R)$ in $\Mod_{R}(\Sp)$.
 Let $D$ be an object  of $\bD$. 
\begin{ddd}\mbox{} 
\begin{enumerate}
\item   $D$  is called $\xi$-acyclic if $S_{\xi}(R)\otimes D\simeq 0$.
\item  $D$ is called $\xi$-complete if $\map_{\bD}(D',D)\simeq 0$ for all $\xi$-acyclic objects $D'$ of $\bD$.
\end{enumerate}
\end{ddd}
 
 The following assertion is well-known. The case of $\bC=\Sp$ and $\xi$ a prime in $\pi_{0}S\cong \Z$ is considered in  \cite{zbMATH03649651}, and the general case is completely analogous.
\begin{prop}\label{erjigowergfregfw}\mbox{}
\begin{enumerate}
\item\label{werkjgowerferwfqewdqewdacf} There exists a Bousfield localization
$$L_{\xi}:\bD\leftrightarrows  L_{\xi}\bD:\incl\ ,$$
where $L_{\xi}\bD$ consists of the $\xi$-complete objects, and where the kernel of  $L_{\xi}$ consists precisely of the $\xi$-acyclic objects.
\item \label{wekgopwegrfwerfrf}  An object $D$ in $\bD$ is $\xi$-acyclic if and only if $\xi$ acts on $D$ as an equivalence.
\item\label{erkhooeprgtrgertg}  For any object  $D$ in $\bD$ we have an equivalence
$$L_{\xi}(D)\simeq \lim_{n\in \nat} D/\xi^{n}\ ,$$
where $D/\xi^{n}:=\Cofib(D\xrightarrow{\xi^{n}} D)$. \end{enumerate}
\end{prop}

 \subsection{Two homological functors}
 
 Let $G$ be a second countable locally compact  topological group and
 recall the functor $\kk^{G}C_{0}:G\LCH_{+}^{\op}\to \KK^{G}$ from \eqref{vrevjeoivfevsdfvdfvsdv}. This functor is homotopy invariant and excisive for  invariant split-closed  decompositions   in the sense that, for any decomposition   of $X$ in $ G\LCH_{+}^{\op}$    into two invariant split-closed subspaces  $Y$ and $Z$   (see \cite[Def. 1.11]{KKG}), we have  a pull-back square
 $$ \xymatrix{\kk^{G}C_{0}(X)\ar[r]\ar[d] &\kk^{G}C_{0}(Y) \ar[d] \\ \kk^{G}C_{0}(Z)\ar[r] &\kk^{G}C_{0}(Y\cap Z) }\ . $$ These assertions follow from \cite[Thm. 5.2]{KKG}.    
\begin{ddd}[{\cite[Def. 1.14]{KKG}}]We   define  the  $\KK^{G}$-valued  analytic $K$-homology functor  by   $$L^{G}:=  \underline{\KK^{G}}(\kk^{G}C_{0}(- ) ,\beins_{\KK^{G}}):    G\LCH_{+}  \to \KK^{G} \ .$$ \end{ddd}
   It is  homotopy invariant and excisive for invariant split-closed decompositions (with the obvious modification of the definition for covariant functors).

  Let  $\ell:\Top \to \Spc$ be the usual presentation of the $\infty$-category of spaces as the Dwyer-Kan localization of $\Top$ at the weak equivalences. We consider $G\Top$ as a topological enriched category. If we insert a $G$-space $X$ into the second argument  functor
  $$\Map_{G\Top}(-,-):G\Top^{\op}\times G\Top\to \Top\ ,$$ then we get a continuous functor $\Map_{G\Top}(-,X):G\Top\to \Top$.
 After composition with $\ell$,  its restriction to $G$-orbits induces
  a functor
 \begin{equation}\label{dcscqcwdcqda}\ell^{G}:G\Top\to \PSh(G\Orb)\ ,  \quad X\mapsto (S\mapsto \ell\Map_{G\Top}(S,X))\ .
\end{equation} 
 Recall that an equivariant weak equivalence in $G\Top$ is a map $X\to Y$ that induces a weak equivalence
  $X^{H}\to Y^{H}$ of fixed point spaces  for all closed subgroups $H$ of $G$.  By Elmendorf's theorem, the functor  $\ell^{G}$ presents its target as the Dwyer-Kan localization of $G\Top$ at the weak equivalences.
      We can consider $\ell^{G}$ as the universal (unstable) equivariant homology theory.

   Note that $\KK^{G}$ is cocomplete and \cref{wijotghwrthggbdh} of $V^{G}$.
By the universal property of the Yoneda embedding
$y$, we  
 get a canonical factorization
\begin{equation}\label{sdvmdfvafvfdfvsdfvsdfv}\xymatrix{G\Orb^{\comp}\ar[dr]_{y}\ar[rr]^{V^{G}}&&\KK^{G}\\&\PSh(G\Orb^{\comp})\ar@{..>}[ur]_{\hat V^{G}}&}
\end{equation} 
 where  $\hat V^{G}$ preserves colimits.
   
   We now assume that $G$  is compact so that $G\Orb^{\comp}\simeq G\Orb$.
  \begin{ddd} \label{weokgpweferfwerf}If $G$ is compact, then we define the homotopy theoretic  $\KK^{G}$-valued analytic  $K$-homology  as the composition 
  $$\tilde L^{G}:G\Top \xrightarrow{\ell^{G} } \PSh(G\Orb) \xrightarrow{\hat V^{G} } \KK^{G}  \ .$$\end{ddd}
 This functor is  in particular 
  homotopy invariant and excisive for invariant open decompositions.
  
  \begin{rem}
  In \cref{weokgpweferfwerf}, we must restrict to compact groups since the functor $V^{G}$ is only defined on compact orbits. 
  In the case of discrete groups which are not necessarily finite,  we have an alternative definition of  an equivariant homotopy theoretic $\KK^{G}$-valued homology theory  using the functor  $\kk^{G}_{\nCcat}\circ \C[-]$ mentioned in \cref{9qriuzghiqergrfqrff} in place of $\kk^{G}C_{0}(-)$. This version is considered in \cite{bel-paschke}. \hB
    \end{rem}

%
%


 Let $G\LCH_{+}^{\hfin}$ denote the full subcategory of $G\LCH_{+}$ of spaces  which are homotopy  retracts of finite $G$-$CW$-complexes. The following results says that the analytic version of the $\KK^{G}$-valued equivariant  $K$-homology $L^{G}$  coincides with its homotopy theoretic counterpart $\tilde L^{G}$, at least on homotopy finite spaces.

\begin{prop}\label{wrtkjhoprgwerfwrefwerf}
We assume that $G$ is compact. We have an equivalence of functors
$$L^{G}_{|G \LCH_{+}^{\hfin}}  \simeq \tilde L^{G}_{|G \LCH_{+}^{\hfin} }\ .$$
\end{prop}
\begin{proof}
Both functors $L^{G}$ and $\tilde L^{G}$ are homotopy invariant and excisive for cell attachments.
The restrictions of the two functors to $G\Orb$ are equivalent by definition.
 This equivalence  therefore extends essentially uniquely to an equivalence of functors on $G \LCH_{+}^{\hfin}$, see  \cite[Cor. 10.9]{bel-paschke}.
\end{proof}

\subsection{A model for the classifying space of proper subgroups in homotopy theory}\label{okgpwefrfrweffsd}
\newcommand{\odd}{\mathrm{odd}}

We assume that $G$ is a finite group. We let $E_{\prp}G$ in $\PSh(G\Orb)$ be the classifying space of proper subgroups of $G$. It is characterized by the property that
\begin{equation}\label{bfgbfdbtrgrdgbbd}E_{\prp}G(G/H)\simeq \left\{\begin{array}{cc} *& H\subsetneqq G	 \\ \emptyset&H=G  \end{array} \right.  \ .
\end{equation}  In the following,  we provide an explicit model for $E_{\prp}G$.
This is a well-known construction \cite{zbMATH03341701}.

Let $V$ be the orthogonal complement of the trivial representation in the complex Hilbert space $L^{2}(G)$. Let $V^{n}:=\bigoplus_{i=1}^{n}V$ and consider the canonical embeddings $V^{n}\to V^{n+1}$. For any finite-dimensional unitary representation $W$ of $G$, we consider the  unit sphere $S(W)$ as a compact $G$-topological space in $G\Top$. Recall the functor $\ell^{G}$ from \eqref{dcscqcwdcqda}.
\begin{prop}\label{wrtijgoiwergerfgw9} For finite $G$, 
we have an equivalence  $$E_{\prp}G:=\colim_{n\in \nat} \ell^{G}(S(V^{n}))\ .$$
\end{prop}
\begin{proof}
%

 If  $H$ is a proper subgroup of $G$,
then the subspace $V^{H}$ of $H$-fixed vectors is non-zero.
For example,  the function  $\chi_{H}- \frac{|H|}{|G|}$ in $L^{2}(G)$  is non-zero  and belongs to $V^{H}$, where $\chi_{H}$ is the characteristic function of $H$. This implies that $S(V^{n})^{H}$ is a sphere of dimension bigger than $n$. Since $\colim_{n\in \nat} \ell(S^{n})\simeq *$, we conclude that 
  $\colim_{n\in \nat }\ell(S(V^{n})^{H})\simeq *$. On the other hand, we have  $(V^{n})^{G}\cong \emptyset$ so that  $
\colim_{n\in \nat} \ell(S(V^{n})^{G})\simeq \emptyset$.

Hence $\colim_{n\in \nat} \ell^{G}(S(V^{n}))$ satisfies the condition
characterizing $E_{\prp}G$ in \eqref{bfgbfdbtrgrdgbbd}.   
   \end{proof}

  \subsection{The KK-theory model for the classifying space of proper subgroups.}\label{ekgowpegwerfrewfwrf}
  
 Recall the functor $\hat  V^{G}:\PSh(G\Orb)\to \KK^{G}$ from \eqref{sdvmdfvafvfdfvsdfvsdfv}.
 In this section we provide a formula \eqref{adsvoijqoirfffsg} for $\hat V^{G}(E_{\prp}G)$ in $\KK^{G}$ in terms of  language introduced in \cref{tirwjgowgerwgwerg9}.

  Note that $V$ is a finite-dimensional complex representation of $G$.
  We consider the element \begin{equation}\label{dfvsfdvwerfc}\xi:=\Lambda_{-1}(V):=\sum_{j=0}^{\infty} (-1)^{j}[\Lambda^{j}(V)]
\end{equation}  in $\pi_{0}R(G)$, 
where for any finite-dimensional $G$-representation $W$, we let    
$[W]$ denote its class in the representation ring $\pi_{0}R(G)$, see \cref{wrjtohipwgwregfwerfw}.

\begin{prop} \label{wtjkigowpegfrferwferfwerf}We assume that $G$ is finite.
We have  an equivalence \begin{equation}\label{adsvoijqoirfffsg}\hat V^{G}(E_{\prp}G)\simeq S_{\xi}(\beins_{\KK^{G}})\ .
\end{equation}
\end{prop}
\begin{proof}
 We start with analysing the effect of the  inclusion maps   $ S(V^{n})\to S(V^{n+1})$ in $\KK^{G}$.
\begin{lem}\label{ojhkorepzhjkoeprhtrgetr}
We have a sequence of fibre sequences \begin{equation}\label{twgregegrewwe}\xymatrix{\vdots\ar[d]&\vdots\ar[d]&\vdots\ar[d]\\\beins_{\KKG}\ar[r]^{\xi^{n+1}}\ar[d]^{\xi}&\beins_{\KKG}\ar[r]\ar[d]&\ar[d]\kkG(C(S(V^{n+1})))\ar[d]\\\beins_{\KKG}\ar[d]\ar[r]^{\xi^{n}}&\ar[d]\beins_{\KKG}\ar[r]&\kkG(C(S(V^{n})))\ar[d]\\\vdots&\vdots&\vdots}\end{equation}
\end{lem}
\begin{proof}

We have a map between semi-split exact   sequences
$$\xymatrix{0\ar[r]&  C(B(V^{n+1}),S(V^{n+1})) \ar[d]\ar[r] &C(B(V^{n+1}))\ar[r] \ar[d]&  \ar[d] C(S(V^{n+1}))\ar[r]& 0\\
0\ar[r]& C(B(V^{n}),S(V^{n}))\ar[r] &C(B(V^{n}))\ar[r]&  C(S(V^{n}))\ar[r]& 0}$$
in $G\nCalg$ that after application of $\kkG$ and some identifications explained below provides the desired map of  horizontal fibre sequences.

The restriction along the embedding of the origin of $B(V^{n})$ induces an equivalence
$$\kkG(C(B(V^{n})))\stackrel{\simeq}{\to} \kkG(\C) \simeq \beins_{\KKG}$$ in $\KKG$. Furthermore, we have a Bott equivalence
\begin{equation}\label{vwervcewrcsdfd}\beins_{\KKG}\stackrel{\simeq}{\to} \kkG(C(B(V^{n}),S(V^{n})))\ .
\end{equation} 
Explicitly, consider the  total alternating power $\Lambda^{*}V^{n}=\Lambda^{*}_{+}V^{n}\oplus \Lambda^{*}_{-}V^{n}$ with the even-odd grading.  We consider the equivariant bundle $B(V^{n})\times  \Lambda^{*}V^{n}\to B(V^{n})$ and 
  the  equivariant  odd bundle
endomorphism $\phi$, considered as a map   $\phi:B(V^{n})\to \End(\Lambda^{*}V^{n})^{\odd}$ , given by $\phi(v):=\epsilon_{v}-i_{v}$, where  $v$ is in $B(V^{n})$, $\epsilon_{v}$ is the exterior product with $v$, and $i_{v}$ is the interior product using the unitary scalar product.
Then $\phi(v)^{2}=-\|v\|^{2}$. In particular $\phi$ is invertible on $S(V^{n})$ and $(\Lambda^{*}V^{n},\phi)$ defines a class $\beta_{V^{n}}$ in $\KKG_{0}(\C,C(B(V^{n}),S(V^{n})))$.  It has been shown by Kasparov  \cite{kasparovinvent} that this class induces the Bott equivalence \eqref{vwervcewrcsdfd}.  
The composition \begin{equation}\label{ogksgpbsggregw}\beins_{\KKG}\to \kkG(C(B(V^{n}),S(V^{n}))) \to  \kkG(C(B(V^{n}))\simeq \beins_{\KKG}
\end{equation} 
is then given by the multiplication by the class $\Lambda_{-1}(V^{n})$ of $\Lambda^{*}V^{n}$ in $R(G)$. Using the fact that $\Lambda_{-1}
$ turns sums into products we conclude that
\eqref{ogksgpbsggregw} multiplication by $\xi^{n}=\Lambda_{-1}(V)^{n}$.
Using  the isomorphism $\Lambda^{*}V^{n+1}\cong \Lambda^{*}V\otimes \Lambda^{*}V^{n}$,  the restriction of the  relative class $\beta_{V^{n+1}}$ of $(\Lambda^{*}V^{n+1},\phi)$ along
$(B(V^{n+1}),S(V^{n+1}))\to (B(V^{n}),S(V^{n}))$ can be identified on the level of representatives with the product of the class $\beta_{V^{n}}$ of 
$(\Lambda^{*}V^{n},\phi)$ by the class of $\Lambda^{*}V$ in $R(G)$, hence  by $\xi$.
This explains the identification of the
  left vertical map in \eqref{twgregegrewwe}.
  \end{proof}

We have a chain of equivalences (natural for $n$ in the poset $\nat$)
\begin{equation}\label{erwgwpojpoewferfrewf} \underline{\KK^{G}}( \kk^{G}(C(S(V^{n}))),\beins_{\KK^{G}})  \simeq L^{G}(S(V^{n}) ) \stackrel{\cref{wrtkjhoprgwerfwrefwerf}}{\simeq} \tilde L^{G}(S(V^{n}))  \ .
\end{equation}  
By \cref{wrtijgoiwergerfgw9},  we obtain
$$ \hat V^{G}(E_{\prp}G) \simeq \colim_{n\in \nat } \underline{\KK^{G}}( \kk^{G}(C(S(V^{n}))),\beins_{\KK^{G}}) \ .$$
Applying $\underline{\KK^{G}}(-,\beins_{\KK^{G}})$ to the diagram \eqref{ojhkorepzhjkoeprhtrgetr}
and using \eqref{erwgwpojpoewferfrewf},  we obtain a diagram of fibre sequences
\begin{equation}\label{twgregegrewwwwe}\xymatrix{\vdots\ar[d]&\vdots\ar[d]&\vdots\ar[d]\\ \tilde L^{G}( V^{n}) 
\ar[r]^{}\ar[d]&\beins_{\KK^{G}}\ar[r]^{\xi^{n}}\ar[d]^{\id} &\ar[d]^{\xi}\beins_{\KK^{G}} \ar[d]\\\tilde L^{G}( V^{n+1}) \ar[d]
\ar[r]^{ } &\ar[d]\beins_{\KK^{G}}\ar[r]^{\xi^{n+1}}&\beins_{\KK^{G}}\ar[d]\\\vdots\ar[d]&\vdots\ar[d]&\vdots\ar[d]\\ 
 \hat V^{G}(E_{\prp}G) \ar[r]&\beins_{\KK^{G}}\ar[r]&\colim_{\nat} (\beins_{\KK^{G}}\xrightarrow{\xi} \beins_{\KK^{G}}\xrightarrow{\xi}\dots)}\end{equation}
 In order to identify the lower left corner with $S_{\xi}(\beins_{\KK^{G}})$, we use \eqref{dvsdvfdvfdvsfvsdfvdfsvdfsverererer} and  \cref{weijgowiergijrefrewf9}.\ref{regijwerogerfwrefrwefw11}.
 \end{proof}





%
%

\begin{theorem}[\cref{wekopgwerfrewfw}]\label{wekopgwerfrewfw1}
 We assume that $G$ is finite.
The cofibre sequence \eqref{vwvssfdvsfvs} is equivalent to the cofibre sequence
\begin{equation}\label{sfvsdfvsfdcssc111}S_{\xi}(\K^{G}_{-}(*))\to \K^{G}_{-}(*)\to  \K^{G}_{-}(*)[\xi^{-1}]\end{equation}  of functors from $\KK^{G}$ to $\Mod_{R(G)}(\Sp)$.
\end{theorem}
\begin{proof} 
Evaluating  \eqref{dvsdvfdvfdvsfvsdfvdfsvdfsverererer} on $\beins_{\KK^{G}}$, we obtain the fibre sequence 
   \begin{equation}\label{bwkopwerf}S_{\xi}( \beins_{\KK^{G}})\to  \beins_{\KK^{G}}\to  \beins_{\KK^{G}}[\xi^{-1}]\ .
\end{equation} 
We apply the  functor $\K^{G}(-\otimes A)$ to the sequence
\eqref{bwkopwerf}, where $\K^{G}$ is the equivariant $K$-theory 
 from \eqref{werfwerfreffwefrefwrefwref}.    
Since  $\K^{G}(-\otimes A)$ is a morphism between $\Mod_{R(G)}(\Sp)$-modules in $\Pr^{L}_{\st}$,  it commutes with the endofunctors $(-)[\xi^{-1}]$ and $S_{\xi}$ and we obtain the fibre sequence \eqref{sfvsdfvsfdcssc}.

 It remains to identify the first map in this fibre sequence with the first map in \eqref{vwvssfdvsfvs}.
 As seen in the proof of \cref{wtjkigowpegfrferwferfwerf}, the fibre sequene 
  \eqref{bwkopwerf} is equivalent to the fibre sequence \begin{equation}\label{sdfbsdvfdvfsdvfsdvsdfv}\hat V^{G}( E_{\prp}G) \to  \beins_{\KK^{G}}\to \colim_{\nat} ( \beins_{\KK^{G}}\xrightarrow{\xi}  \beins_{\KK^{G}}\xrightarrow{\xi}\dots)
\end{equation}
  If we apply $\K^{G}(-\otimes A)$ to the fibre sequence \eqref{sdfbsdvfdvfsdvfsdvsdfv}, then we obtain  a fibre sequence which is equivalent to \eqref{sfvsdfvsfdcssc}, and whose first map ist
  $$\K^{G}(\hat V^{G}(E_{\prp}G)\otimes A)\to \K^{G}_{A}(*)\ .$$
  We now use that 
 $E_{\prp}G\simeq \colim_{G_{\prp}\Orb}y $, where $y:G\Orb \to \PSh(G\Orb)$ is the Yoneda embedding.
 Since $\hat V^{G}\circ y\simeq V^{G}$ by definition, 
 this implies that   $$\K^{G}( \hat V^{G}(E_{\prp}G)\otimes A)\simeq \colim_{G_{\prp}\Orb}\K^{G}( V^{G}\otimes A)
\stackrel{\eqref{getgoijif0ewrfrwef}}{\simeq}  \colim_{G_{\prp}\Orb}\K^{G}_{A}\ .$$\end{proof}

%
%


\subsection{Reduction to cyclic subgroups}\label{kowgpwgwrgwr9}

In this subsection, we show \cref{kopthethegtrg1} and \cref{geropfrfwrefrefrwfwrf}.
\begin{lem} \label{kopthethegtrg1}
If $G$ is finite and not cyclic, then $\xi=0$.
\end{lem}
\begin{proof} 
 Let $\Cyc(G)$ be the set of cyclic subgroups of $G$. 
The map
$$\pi_{0}R(G)\to \bigoplus_{C\in \Cyc(G)} \pi_{0}R(C)\ , \quad \rho\mapsto \oplus_{C\in \Cyc(G)}\rho_{|C}$$
is injective. Indeed, if $\rho$ is in the kernel of this map, then $\rho_{|C}=0$ for all cyclic subgroups. 
Since every $g$ in $G$ generates a cyclic subgroup, we have
  $\tr \rho(g)=0$ for all $g$ in $G$. This implies $\rho=0$.

We have $\Lambda_{-1}(\C^{n})=0$ for $n\ge 1$.
If $W$ is a complex representation of $G$ with $W^{G}\not=0$, then
$$\Lambda_{-1}(W)=\Lambda_{-1}(W^{G}\oplus (W^{G})^{\perp})=\Lambda_{-1}(W^{G})\Lambda_{-1}((W^{G})^{\perp})=0\ .$$

If $G$ is not cyclic, then any   $C$ in $\Cyc(G)$ is proper and $V^{C}\not=0$.
This implies $\Lambda_{-1}(V)_{|C}=\Lambda_{-1}(V_{|C})=0$ for every $C$ in $\Cyc(G)$.
Hence $\Lambda_{-1}(V)=0$.
\end{proof}

%
%
%
%
%

\begin{prop}[\cref{geropfrfwrefrefrwfwrf}]\label{geropfrfwrefrefrwfwrf1} If $G$ is finite, then
we have an equivalence
$$\colim_{G_{\Cyc}\Orb}\K^{G}_{A}\xrightarrow{\simeq} \K(A\rtimes_{r}G)$$ for every $A$ in $\KK^{G}$.
\end{prop}
\begin{proof} 
 For a family $\cF$ of subgroups of $G$, we let $p^{\cF}:G_{\cF}\Orb\to *$ denote the canonical functor. Then $\colim_{G_{\cF}\Orb}$  is the left Kan extension functor $p^{\cF}_{!}$  
along $p^{\cF}$.

If $\cF$, $\cF'$ are two families of subgroups such that $\cF'\subseteq \cF$, then we let $p_{\cF}^{\cF'}:G_{\cF'}\Orb \to  G_{\cF}\Orb $ denote the inclusion.   
 For the family  $ \All_{G}$ of all subgroups of $G$, we obtain
 $$p^{\All_{G}}_{!}\K^{G}_{A}\simeq \K^{G}_{A}(*)\simeq \K(A\rtimes_{r}G)\ .$$
 For any family $\cF$ of subgroups of $G$, the map
 $$c_{G}^{\cF}:p^{\cF}_{!} \K^{G}_{A} \to p^{\All_{G}}_{!}\K^{G}_{A}$$  (in the domain we omitted restriction of the functor to $G_{\cF}\Orb$) is induced by the inclusion
 $G_{\cF}\Orb\to G_{\All_{G}}\Orb$. If $\cF'$ is a second family such that $\cF'\subseteq \cF$, then
  we have a factorization  of $c_{G}^{\cF'}$ as 
$$c^{\cF'}_{G}:p^{\cF'}_{!} \K^{G}_{A} \xrightarrow{c^{\cF'}_{\cF}} p^{\cF}_{!} \K^{G}_{A} 
\xrightarrow{c_{G}^{\cF}} p^{\All_{G}}_{!}\K^{G}_{A}\ .$$

 We consider a family $\cF$ of   subgroups of $G$ containing $\Cyc$. The set of conjugacy classes of subgroups in $\cF$ is a poset. If $[H]$ is a maximal element in this poset  that is not cyclic, then we can define a new such family $\cF\setminus {[H]}$ by removing $[H]$.

 We now assume by induction that
 $$c_{G}^{\cF}:p^{\cF}_{! }\K^{G}_{A} \to p^{\All_{G}}_{!}\K^{G}_{A}$$ is an equivalence and let $\cF':=\cF\setminus [H]$ for some maximal  conjugacy class $[H]$ in $\cF$.
  In order to show that $c_{G}^{\cF'}$ is an equivalence, it suffices to show that $c^{\cF'}_{\cF}$ is an equivalence.
 By the pointwise formula for the left Kan extension, we must show for all $G/L$  in $G_{\cF}\Orb$ that   $\colim_{G_{\cF'}\Orb_{/(G/L)}} \K^{G}_{A}\to \K^{G}_{A}(G/L)$ is an equivalence. The functor $i:L\Orb\to G\Orb$, $S\mapsto G\times_{L}S$ induces 
  an equivalence $  L_{\cF'\cap L}\Orb \xrightarrow{\simeq}G_{\cF'}\Orb_{/(G/L)} $, and using 
 \cref{weijogpwergwergf91t} we can identify 
  the map in question  with  $c_{L}^{\cF'\cap L}:p^{\cF'\cap L}_{!}\K^{L}_{A'}  \to 
 \K^{L}_{A'}(*)$, where $A':=\Res^{G}_{L}(A)$. If $L\not\in [H]$, then $ \cF'\cap L=\All_{L}$ and $c^{\cF'\cap L}$ is clearly an equivalence.
 If $L\in [H]$, then $\cF'\cap L$ is the family of proper subgroups of $L$ and $c^{\cF'\cap L}$ is an equivalence by 
 \cref{gjkwegokwerpferfwerf}.
\end{proof}

\section{The case of cyclic p-groups}
 
\subsection{p-completion}

Calculating up to $\xi$-completion, for $\xi$ in $\pi_{0} R(G)$  defined by    \eqref{dfvsfdvwerfc}, using
  \cref{wjogpwergerfrefwerfw} might still be a non-trivial task. If $p$ is a prime and one is interested in $p$-torsion phenomena, then it 
is more natural to calculate up to $p$-completion. Here we view $p$ as a class in
$\pi_{0}R(G)$ represented by the trivial $p$-dimensional representation of $G$. This allows to apply the language developed in \cref{tirwjgowgerwgwerg9}.

Recall that the $p$-completion of an  abelian group $A$  is defined by
$$\widehat{A}_{p}:= \lim_{n\in \nat} (\dots \to A/p^{2}A\to A/pA\to A)\ .$$ 

Let $E$ be in $\Mod_{R(G)}(\Sp)$. The following is well-known.
\begin{lem}
If $\pi_{*}(E)$ is a degree-wise   finitely generated abelian group, then we 
  have an isomorphism
$$\pi_{*}(L_{p}(E))\cong \widehat{\pi_{*}E}_{p}\ .$$
\end{lem}
\begin{proof} By \cref{erjigowergfregfw}.\ref{erkhooeprgtrgertg},
 one knows that the $p$-completion functor is given by the formula
$$L_{p}(E)\simeq \lim_{n\in \nat} (\dots\to E/p^{2}\to E/p\to E)\ .$$
Applying $\pi_{*}$, we get the Milnor exact sequence
\begin{equation}\label{qwefpojwopfdqewdeqwdeqwdq}0\to \lim{}^{1}_{n\in \nat} \pi_{*+1}(E/p^{n})\to \pi_{*}(L_{p}(E))\to \widehat{\pi_{*}E}_{p}\to 0\ .
\end{equation}
If $\pi_{*}E$ is finitely generated, then $\pi_{*}E/p^{n}$ is finite for every $n$ and, consequently,  the $\lim^{1}$-term vanishes.
\end{proof}

Let $\xi$ be any class in $\pi_{0}R(G)$ and consider the endofunctor $(-)[\xi^{-1}]$ of a $\Mod_{R(G)}(\Sp)$-module $\bC$.
\begin{lem}\label{werijgowergwrefwerfw9}
If there exists $n$ in $\nat$ such that $\xi^{n}$ is divisible by $p$, then $L_{p}( (-)[\xi^{-1}])\simeq 0$.
 \end{lem}
\begin{proof}  By assumption  there exists  $\eta$ in $\pi_{0}R(G)$ such that 
  $ \xi^{n}=p\eta$. Let $E$ be an object of $\bC$. Then
  $p$ acts as an equivalence on $E[\xi^{-1}]$.  But then
   $E[\xi^{-1}]$ is $p$-acyclic by \cref{erjigowergfregfw}.\ref{wekgopwegrfwerfrf}. \end{proof}


Let now $\xi$ be defined by \eqref{dfvsfdvwerfc}.
 \begin{kor}\label{thijwoptrgrgwergrgwrf}
If there exists $n$ in $\nat$ such that $\xi^{n}$ is divisible by $p$,
then the assembly map $$\colim_{G_{\prp}\Orb}\K^{G}_{A}\to  \K(A\rtimes_{r}G)$$
becomes an equivalence after $p$-completion. 
\end{kor}
\begin{proof}  
We apply $L_{p}$ to the cofibre sequence \eqref{sfvsdfvsfdcssc}. By \cref{werijgowergwrefwerfw9} and \cref{wekopgwerfrewfw},
 we get an equivalence
$$ L_{p}(\colim_{G_{\prp}\Orb}\K^{G}_{A})\simeq L_{p}(S_{\xi}  (\K^{G}_{A}(*)))\stackrel{\simeq}{\to} L_{p}(\K^{G}_{A}(*))\simeq L_{p}(\K(A\rtimes_{r}G))\ .$$
 \end{proof}

\subsection{Cyclic groups}\label{werogkpwergerwfwf}

In this subsection, we study the divisibity properties of $\xi$ defined in \eqref{dfvsfdvwerfc}
for the cyclic group $C_{n}:=\Z/n\Z$ for an integer $n$.
The result will be used to show \cref{erijgowergrefwrefr}.

 The (algebraic) representation ring of $C_{n}$ is given by
$$\pi_{0}R(C_{n})\cong \Z[x]/(x^{n}-1)\ ,$$ where $x$ corresponds to the representation $G\to \Aut(\C)$ that sends the generator 
$[1]$ to the  primitive $n$'th root of unity $\exp(\frac{2\pi i}{n})$. 

\begin{prop}\label{wtjigowefrfwref}
If $n=p^{k}$ for some prime $p$ and $k\in \nat\setminus\{0\}$, then $\xi^{2}=p\xi$. Otherwise
$\xi^{2}=\xi$ and $\xi$ is not divisible by any prime.
\end{prop}
\begin{proof}
Recall from \cref{okgpwefrfrweffsd} that $V$ denotes the complement of the trivial representation of $G$ in $L^{2}(G)$. Under the identification of the representation ring above
we have $[V]=\sum_{\ell=1}^{n-1}x^{\ell}$,  and hence, $\xi=\Lambda_{-1}(V)=\prod_{\ell=1}^{n-1} (1-x^{\ell})$, see \eqref{dfvsfdvwerfc}.
We have an injective map of algebras
$$j:\Z[x]/(x^{n}-1)\to \C^{C_{n}}$$ that sends the class of $f$ in $\Z[x]$ to the function 
$j(f):[k]\mapsto f(\exp(\frac{2\pi i k}{n}))$. 
Then $\xi$ goes to the function
$$j(\xi):([q]\mapsto \prod_{\ell=1 }^{n-1}(1-\exp(\frac{2\pi i q \ell}{n}))\ .$$
We have
$$j(\xi)(c)=\left\{\begin{array}{cc}0 & \mbox{$c$ does not generate $C_{n}$}\\  \Phi_{n}(1)& \mbox{$c$  generates $C_{n}$} \end{array} \right.\ ,$$
where
$\Phi_{n}\in \Z[x] $ is the $n$th cyclotomic polynomial.
 \begin{lem}\label{tkohgprthrheht}
We have
$$\Phi_{n}(1)=\left\{\begin{array}{cc}  p&\mbox{$n=p^{k}$ for some prime $p$}\\1
&\mbox{$n$ is not a prime power} 
 \end{array} \right.\ .$$
\end{lem}
 \begin{proof}
 The proof uses some well-known identities  for cyclotomic polynomials.
 If $n=p^{k}r$ such that $p$ does not divide $r$, then
 $\Phi_{n}(x)=\Phi_{pr}(x^{p^{k-1}})$.
 In particular we get 
 $\Phi_{n}(1)= \Phi_{pr}(1)$.
 We conclude that
 $\Phi_{n}(1)=\Phi_{q}(1)$, where $q$ is the product of the  prime divisors of $n$.
  
 If $q$ is odd and of the form $pr$ with $r\not=1$, then we have
   $ \Phi_{q}(x)=\frac{\Phi_{r}(x^{p})}{\Phi_{r}(x)}$ which implies that 
 $\Phi_{q}(1)=1$. We now consider the case
  $q=2r$ for $r\not=1$ in which $\Phi_{q}(x)=\Phi_{r}(-x)$. If $r$ is an odd prime, then
  $\Phi_{r}(-1)=1$. Otherwise
  $r=pr'$ and $\Phi_{r}(-1)= \frac{\Phi_{r'}((-1)^{p})}{\Phi_{r'}(-1)}=1$.

  If   $n=p^{k}$ for some $k$ and  prime $p$, then $q=p$ and
 $\Phi_{n}(1)=p$.  
  \end{proof}

 \cref{tkohgprthrheht} implies that $$j(\xi^{2})=\left\{\begin{array}{cc} pj(\xi) & \mbox{$n=p^{k}$ for some prime $p$ and $k$ in $\nat\setminus \{0\}$}\\ j(\xi) & \mbox{else} \end{array} \right. \ .$$
 Since $j$ is injective, we can conclude that
 $$\xi^{2}=\left\{\begin{array}{cc} p\xi & \mbox{$n=p^{k}$ for some prime $p$ and $k$ in $\nat\setminus \{0\}$}\\ \xi & \mbox{else} \end{array} \right. \ .$$
In the second case, if $p$ would divide $\xi$, then $\xi=p\eta$ for some $\eta\in \Z[x]/(x^{n}-1)$.
Then we can conclude that $p^{2^{\ell}-1}\eta^{2^{\ell}}=\eta$ for any $\ell$ in $\nat$ and therefore $\eta$ would be $p$-divisible.  But    $\Z[x]/(x^{n}-1)$ does not contain non-trivial $p$-divisible elements.
This completes the proof of \cref{wtjigowefrfwref}.
\end{proof}

For a prime $p $    and  positive integer $k$ in $\nat$, we consider the group  $G:=C_{p^{k}}$. We furthermore consider $A$    in $\KK^{C_{p^{k}}}$. 
 \begin{prop}[\cref{gijeoggerwfreffvvsdfvsfdvsfv}]\label{gijeoggerwfreffvvsdfvsfdvsfv1}
 If $p$ acts as an equivalence in $A$, then the evaluation  
 $$S_{\xi}(\K^{C_{p^{k}}}_{A}(*))\to \K^{C_{p^{k}}}_{A}(*)\to  \K^{C_{p^{k}}}_{A}(*)[\xi^{-1}]$$ of the sequence \eqref{sfvsdfvsfdcssc} at $A$ splits naturally.
 \end{prop} \begin{proof} By \cref{wtjigowefrfwref}, we have $\xi^{2}=p\xi$. 
     We let $\KK^{C_{p^{k}}}[p^{-1}]$ be the localizing subcategory of $ \KK^{C_{p^{k}}}$ on objects on which $p$ acts as an equivalence. 
    On this category, we can define the idempotent $\pi:=\frac{\xi}{p}$. Since   $A\in \KK^{C_{p^{k}}}[p^{-1}]$, we obtain a decomposition 
 $A\simeq \pi A\oplus (1-\pi)A$.

   Now $\xi$ acts as an equivalence on $\pi A$ (with inverse $\frac{1}{p}$) and annihilates $(1-\pi)A$.
 The cofibre sequence
 $$S_{\xi}(\K^{C_{p^{k}}}_{A}(*))\to \K^{C_{p^{k}}}_{A}(*)\to  \K^{C_{p^{k}}}_{A}(*)[\xi^{-1}]$$  splits as a sum of the two copies of the cofibre sequences from \eqref{sfvsdfvsfdcssc} applied to $\pi A$ and $(1-\pi )A$. 
 We now have $S_{\xi}(\pi A)\simeq 0$  and $(1-\pi)A[\xi^{-1}]\simeq 0$. Adding up what remains, we obtain the desired decomposition
   $$ \K^{C_{p^{k}}}_{A}(*)\simeq S_{\xi}(\K^{G}_{(1-\pi) A}(*))\oplus  \K^{C_{p^{k}}}_{\pi A}(*)[\xi^{-1}]\simeq  S_{\xi}(\K^{C_{p^{k}}}_{A}(*))\oplus  \K^{C_{p^{k}}}_{ A}(*)[\xi^{-1}]\ .
   $$
   \end{proof}

\begin{theorem}[\cref{erijgowergrefwrefr}]\label{erijgowergrefwrefr1} For any $A$ in $\KK^{C_{p^{k}}}$, 
 the assembly map
\begin{equation}\label{wervwoerijvewoivwev111} \colim_{BC_{p^{k}}} \K(\widehat \Res^{G}(A))\to   \K(A\rtimes_{r}C_{p^{k}})
\end{equation}   
 becomes    an   equivalence after $p$-completion.
\end{theorem}
\begin{proof} 
%
%
We perform a similar induction as in the proof of   \cref{geropfrfwrefrefrwfwrf}  over the family of all subgroups of $C_{p^{k}}$.  All non-trivial subgroups of $C_{p^{k}}$ are groups of the form $C_{p^{k'}}$ for some $k'$ between $1$ and $k$. For every group of this form we know that
$$\colim_{{{C_{p^{k'}}}_{\prp}\Orb}} \K_{A^{'}}^{C_{p^{k'}}} \to \K_{A^{'}}^{C_{p^{k'}}}(*)$$ (with $A':=\Res^{C_{p^{k}}}_{C_{p^{k'}}}(A)$)
becomes an equivalence after $p$-completion by \cref{thijwoptrgrgwergrgwrf} and \cref{wtjigowefrfwref}.
If we repeat the argument from the proof of \cref{geropfrfwrefrefrwfwrf} this observation is used in the last line of the argument instead of  \cref{gjkwegokwerpferfwerf}.  
The induction ends at the family of the trivial group where we get the statement that the map \eqref{wervwoerijvewoivwev111} becomes an equivalence after $p$-completion.
%
%
%
%
%
%
\end{proof}

\begin{lem}\label{wtgopwrfefewfsf}
If $A\simeq \Res_{C_{p}}(A')$ for some $A'$ in $\KK$ such that $K_{*}(A')$ is annihilated by $p^{n}$ for some $n$ in $\nat$,  then
\begin{equation}\label{vdfsvsdfvsdefeffvdfvsd} \colim_{BC_{p} } \K( A )\xrightarrow{\simeq} \K(A\rtimes_{r}C_{p})\ .
\end{equation} 
\end{lem}
\begin{proof}
The values of $V^{C_{p}}(S)\rtimes C_{p}$  are equivalent to  finite sums of copies of $\beins$ and therefore  belong to  the bootstrap class $\UCT$.   We can thus apply   \cref{qrejigoqrgregwergwerg} and obtain an equivalence 
  $$\K^{C_{p}}_{A}\simeq \K(A')\otimes_{KU} \K^{C_{p}}_{\beins_{\KK^{C_{p}}}}\ .$$
 By \cref{wekopgwerfrewfw}, the map  \eqref{vdfsvsdfvsdefeffvdfvsd} is then equivalent to the first map in the fibre sequence \begin{equation}\label{fpoqkwopefewdqwedq}
 \K(A')\otimes_{KU} S_{\xi}( \K^{C_{p}}_{\beins_{\KK^{C_{p}}}}(*)) \to   \K(A')\otimes_{KU}  \K^{C_{p}}_{\beins_{\KK^{C_{p}}}}(*)\to  \K(A')\otimes_{KU}  \K^{C_{p}}_{\beins_{\KK^{C_{p}}}}(*)[\xi^{-1}]\ .\end{equation}
Now on the one hand,  $p^{n}$  annihilates $K_{*}(A')$. On the other hand, it   acts as an equivalence on $\K^{C_{p}}_{\beins_{\KK^{C_{p}}}}(*)[\xi^{-1}]$  since $\xi^{2}$ is divisible by $p$ by \cref{wtjigowefrfwref}.  It follows  that  the right-most term in the sequence \eqref{fpoqkwopefewdqwedq} is trivial. Consequently  its first map is an equivalence.
\end{proof}

Recall the definition of the localizing subcategory $\bG(p)$ of  $\KK^{C_{p}}$ given before the statement of \cref{wetghkopwergwerfw}.
\begin{prop}[\cref{wetghkopwergwerfw}]\label{wetghkopwergwerfw1}
If $A$ is in $\bG(p)$, then
the assembly map is an equivalence
\begin{equation}\label{efdewdewdewdq1111} \colim_{BC_{p}}  \K(A)  \stackrel{\simeq}{\to}   \K( A \rtimes_{r}C_{p}) 
\end{equation} 
\end{prop}
\begin{proof} 
Since both sides of   the assembly map \eqref{frewfwrewoifjowrfwrefwerf} preserve colimits it suffices to show that it is an equivalence for the generators of $\bG(p)$.  

If $A=\Ind^{C_{p}}(A'')$ for some $A''$ in $\KK$, then the assertion  is an instance of \cref{werjoigpwegerfwref}.
In the second case where  $A=\Res_{C_{p}}(A')$ for some $A'$ in $\KK$ such that $\K_{*}(A')$  is annihilated by $p^{n}$ for some $n$ in $\nat$, the assertion is given  by  \cref{wtgopwrfefewfsf}.
\end{proof}

\subsection{Some calculations}\label{gwegire0fwrefwrefw}

In the following, $p$ is a prime, $k$ is in $\nat\setminus \{0\}$, $H$ is a subgroup of $C_{p^{k}}$,  $B:=\beins_{\KK^{H}}/\ell$ in $\KK^{H }$ for an integer $\ell$ with $p\not|\ell$, and $A:=B^{\otimes C_{p^{k}}/H}$.
 
\begin{lem}[Lemma \ref{oejgoerpwgrefdssfvfsfsv}]\label{oejgoerpwgrefdssfvfsfsv1}
We have $$\K(A\rtimes_{r}C_{p^{k}})[\xi^{-1}]\simeq R(C_{p^{k}})[\xi^{-1}]/\ell\simeq (KU/\ell)^{\oplus (p-1)p^{k-1}}\ .$$ 
\end{lem}
\begin{proof} In view of \cref{wekopgwerfrewfw},  we have an equivalence
$$K(A\rtimes_{r}C_{p^{k}})[\xi^{-1}]\simeq \Cof^{C_{p^{k}}}(A)\  .$$
%
%
 We now   employ  \cref{ewgojwoepgfrefwfwerfwrf},   stating that  the composition $$\Cof^{C_{p^{k}}}\circ (-)^{\otimes C_{p^{k}}/H}:\KK^{H}\to \Mod_{R(C_{p^{k}})}(\Sp)$$ preserves all colimits, for the marked equivalence in 
      \begin{align*} \Cof^{C_{p^{k}}}(A)&\simeq \Cof^{C_{p^{k}}}((\beins_{\KK^{H}}/\ell)^{\otimes C_{p^{k}}/H})\stackrel{!}{\simeq} \Cof^{C_{p^{k}}}(\beins_{\KK^{H}}^{\otimes C_{p^{k}}/H})/\ell\\&\simeq 
  \Cof^{C_{p^{k}}}(\beins_{\KK^{C_{p^{k}}}})/\ell\simeq R(C_{p^{k}})[\xi^{-1}]/\ell\ .\end{align*}
We have an inclusion of rings  $\pi_{0}R(C_{p^{k}})\hookrightarrow \C^{C_{p^{k}}}$ that  sends
$\rho$ to the function $C_{p^{k}}\ni c\mapsto \tr\rho(c)\in \C$. Using the result from \cref{werogkpwergerwfwf}, under this inclusion the multiplication by 
$\xi$ goes to the multiplication by the function $$c\mapsto \left\{\begin{array}{cc} p& \mbox{$c$ generates $C^{p^{k}}$}\\0 &\mbox{else}  \end{array} \right.\ .$$
We conclude that
$$\pi_{0}R(C_{p^{k}})[\xi^{-1}]\simeq \Z[p^{-1}]^{\oplus (p-1) p^{k-1} }\ ,$$ where the rank of this $\Z[p^{-1}]$-module is the number of generators of $C_{p^{k}}$. Since $p$ is invertible in $\Z/\ell$, we have
$$\pi_{0}R(C_{p^{k}})[\xi^{-1}]/\ell\simeq (\Z/\ell \Z)^{ \oplus (p-1) p^{k-1} }\ .$$ Since $\pi_{1}R(C_{p^{k}})\cong 0$, 
we obtain the eqivalence   $$R(C_{p^{k}})[\xi^{-1}]/\ell\simeq   (KU/\ell)^{\oplus (p-1)p^{k-1}}$$ in $\Mod_{KU}(\Sp)$.     \end{proof}

We now restrict to the case of $k=1$ and the trivial subgroup $H$. Then $A\simeq (\beins_{\KK}/\ell)^{\otimes C_{p}}$.
 \begin{lem}[\cref{wergojowpergrefgrwefwref}]\label{wergojowpergrefgrwefwref1}
  We have $$\colim_{BC_{p}} \K(\widehat \Res^{C_{p}}((\beins_{\KK}/\ell)^{\otimes C_{p}}))\simeq (KU/\ell)^{n_{0}(p)}\oplus (\Sigma KU/\ell)^{n_{1}(p)}$$ with
  $$n_{i}(p):=\left\{\begin{array}{cc} \frac{2^{p-2} +  \frac{(p+1)(p-1)}{2}  }{p}&i=0\\ \frac{2^{p-2}  - \frac{(p-1)^{2}}{2}   }{p}& i=1 \end{array} \right.$$ for $p>2$ and
  $$n_{i}(2):=\left\{\begin{array}{cc} 1&i=0\\ 0&i=1  \end{array} \right.$$
  \end{lem}
\begin{proof}
 We have an equivalence
$$ \K(\widehat \Res^{C_{p}}((\beins_{\KK}/\ell)^{\otimes C_{p}}))\simeq (KU/\ell)^{\otimes C_{p}}$$
in $\Fun(BC_{p},\Mod_{KU}(\Sp))$. If we tensorize the fibre sequence $\beins_{\KK}\xrightarrow{\ell}\beins_{\KK}\to \beins_{\KK}/\ell$ with $\beins_{\KK/\ell}$, then we obtain the fibre sequence $\beins_{\KK}/\ell\xrightarrow{0}\beins_{\KK}/\ell\to \beins_{\KK}/\ell\otimes \beins_{\KK}/\ell$. It   gives an equivalence 
 $KU/\ell\otimes KU/\ell\simeq KU/\ell\oplus \Sigma KU/\ell$ in $\Mod_{KU}(\Sp)$. Inductively,  we obtain  \begin{equation}\label{qfewdqewdqfqf} (KU/\ell)^{\otimes C_{p}}\simeq \bigoplus_{i=0}^{p-1}( \Sigma^{i}KU/\ell)^{\binom{p-1}{i}}
\end{equation} for any positive integer $p$.
 This calculation does not reflect the $C_{p}$-action. But we can see that $\pi_{*}(  (KU/\ell)^{\otimes C_{p}})$ is a sum of copies of $\Z/\ell \Z$. Furthermore, since $p\not| \ell$, we know that
 $p$ acts invertibly on the homotopy groups   so that
$ (KU/\ell)^{\otimes C_{p}}$ in $\Fun(BC_{p},\Mod_{KU}(\Sp))$ is projective by \cref{ekorgpefreqwdqedqed}.
Therefore, by \eqref{vsadoivjadsoivcadscsadc1}, we have an isomorphism
$$\pi_{*} (\colim_{BC_{p}}  (KU/\ell)^{\otimes C_{p}})\cong \colim_{BC_{p}} \pi_{*} ((KU/\ell)^{\otimes C_{p}})\ .$$

%
%
 
We have a cofibre sequence
$$\beins_{\KK}\to \beins_{\KK}/\ell\to \Sigma \beins_{\KK}$$
whose boundary map $ \Sigma \beins_{\KK}\to  \Sigma \beins_{\KK}$ is multiplication by $\ell$.
We use the spectral sequence from \cref{krwfopqwewfewdffqwf} associated to this cofibre sequence in place of \eqref{fqewifhqiofsffeadf} in order to calculate the associated graded of an invariant filtration of the  abelian groups 
$  \pi_{*} ((KU/\ell)^{\otimes C_{p}})$ with $C_{p}$-action.
We use that $\beins_{\KK}^{\otimes X}\simeq \beins_{\KK}$ and $(\Sigma \beins_{\KK})^{\otimes X}\simeq \Sigma^{|X|}\beins_{\KK}$ for any finite set $X$. Specializing \eqref{wrgeoigowerpfwerfwerf}, we obtain
  $$E^{C_{p},1}_{i,j}\cong \bigoplus_{F\in \cP_{p-j}(C_{p})} KU_{i-j+p}\ .$$
  In order to calculate the  matrix coefficients of the differential $d^{C_{p},1}_{i,j}$, we fix an order of the elements of  $C_{p}$.
 If $F'$ is a subset  of $C_{p}$ and $F:=F'\setminus\{f'\}$ for some $f'$ in $F'$, then  
 we define the sign $s(F,F')\in \{\pm 1\}$ as the negative of the parity of the position of $f'$. 
 By \cref{rkogpwergwerfwefwerfwref}, the matrix coefficient $d^{C_{p},1}_{i,j,F,F'}$ of the differential
  $$d^{C_{p},1}_{i,j}: \bigoplus_{F'\in \cP_{p-j}(C_{p})} KU_{i-j+p}\to  \bigoplus_{F\in \cP_{p-j-1}(C_{p})} KU_{i-j+p}$$
  is then given by 
  $$d^{C_{p},1}_{i,j,F,F'}= \left\{\begin{array}{cc} s(F,F')\ell &F\subseteq F'\\0 & else  \end{array} \right.\ .$$
  Note that in the formula for the differential in \cref{rkogpwergwerfwefwerfwref}, we use the identification
  $K(B)^{\otimes F'}\simeq K(B)^{\otimes \{f\}}\otimes K(B)^{\otimes F}$.  Since the homotopy groups of $B$ live in odd degrees, we get the sign from the Koszul sign rules.
  
We now consider the cochain complex of $\Z/2\Z$-graded abelian groups \begin{equation}\label{fwerfwevfsdfvf}
E^{C_{p},1}_{i,0}\xrightarrow{d^{C_{p},1}_{i,0}}E^{C_{p},1}_{i-1,1}  \xrightarrow{d^{C_{p},1}_{i-1,1}}\dots 
 \xrightarrow{d^{C_{p},1}_{i-p+1,p-1}}
E^{C_{p},1}_{i-p,p}
 \end{equation}
where $E^{C_{p},1}_{i-p,p}$ is in degree $0$.
It vanishes if $i-p$ is odd. If $i-p$ is even, then it is $C_{p}$-equivariantly isomorphic (see \cref{ertkohpperthgetr} below for justification)  to the cochain complex \begin{equation}\label{fqwefweoijdowqedqwedqwefcq}(\Z\stackrel{\ell}{\to} \Z)^{\otimes C_{p}}\ ,
\end{equation}  where $\Z\stackrel{\ell}{\to} \Z$ lives in degree $-1$ and $0$.

 \begin{rem} \label{ertkohpperthgetr}
Let $$Z_{-p}\xrightarrow{d_{-p}} Z_{1-p} \xrightarrow{d_{1-p}} \dots \xrightarrow{d_{-2}} Z_{-1}\xrightarrow{d_{-1}} Z_{0}$$ be the cochain complex \eqref{fqwefweoijdowqedqwedqwefcq}. Then  
$$Z_{-k}\cong \bigoplus_{F'\in \cP_{k}(C_{p})} \Z^{\otimes C_{p}\setminus F'}\otimes \Z^{\otimes  F'} \ ,$$
where of course every summand is isomorphic to $\Z$. The matrix coefficient $d_{-k,F,F'}$ of the differential $d_{-k} $   for $F'=F\sqcup \{f\}$ is 
 $$ \hspace{-0.5cm} \Z^{\otimes C_{p}\setminus F'}\otimes \Z ^{\otimes  F'}\xrightarrow{s(F,F')} \Z^{\otimes C_{p}\setminus F'}\otimes \Z^{\{f\}} \otimes  \Z^{\otimes  F}\xrightarrow{\id_{ \Z^{\otimes C_{p}\setminus F'}}\otimes \ell \otimes \id_{ \Z^{\otimes  F}} }   \Z^{\otimes C_{p}\setminus F'}\otimes \Z^{\{f\}} \otimes \Z^{\otimes  F} \cong \Z^{C_{p}\setminus F}\otimes \Z^{\otimes  F}\ . $$ This is precisely the description of the cochain complex \eqref{fwerfwevfsdfvf}  in the case that $i-p$ is even. \hB
\end{rem}

 By an inductive calculation in $D(\Z)$, we obtain the analogue of \eqref{qfewdqewdqfqf}
\begin{equation}\label{fewdwedwqdqwefqfeqfwefeqwf}(\Z/\ell\Z)^{\otimes C_{p}} \simeq \bigoplus_{i=0}^{p-1}(  \Z/\ell[i])^{\binom{p-1}{i}}
\end{equation} 
which describes the cohomology  of the cochain complex \eqref{fqwefweoijdowqedqwedqwefcq} non-equivariantly.
This cohomology  (two-periodically repeated) constitutes the   second page $E^{C_{p},2}_{*,*}$ of the spectral sequence.
 One checks by an explicit counting, using  \eqref{fewdwedwqdqwefqfeqfwefeqwf} and  \eqref{qfewdqewdqfqf},
 that  $$|\bigoplus_{j=0}^{p} E^{C_{p},2}_{i,j}|= |\pi_{i}((KU/\ell)^{\otimes C_{p}})|\ .$$
 This implies that the spectral sequence degenerates at the second page.
 
 We conclude that $\pi_{*}((KU/\ell)^{\otimes C_{p}})$ has a $C_{p}$-invariant filtration
 whose associated graded is $C_{p}$-equivariantly isomorphic to the cohomology of $(\Z \xrightarrow{\ell} \Z )^{\otimes C_{p}}$ (viewed two-periodically). We are interested in the $C_{p}$-coinvariants. Since $p$ acts invertibly on the cohomology and $\pi_{*}((KU/\ell)^{\otimes C_{p}}) $,
 taking coinvariants is an exact functor.  We see that the group of 
   coinvariants  $\colim_{BC_{p}}\pi_{*}((KU/\ell)^{\otimes C_{p}})$ again  has a filtration whose associated graded is isomorphic to the $C_{p}$-coinvariants in the cohomology of $(\Z \xrightarrow{\ell} \Z )^{\otimes C_{p}}$.
   We will see below that the coinvariants in the cohomology  of this cochain complex are again sums of copies of $\Z/\ell \Z$.  Since  
  $ \pi_{*}((KU/\ell)^{\otimes C_{p}})$ 
itself is a sum of copies of $\Z/\ell \Z$, we can conclude  (using, e.g., the classification of finite abelian groups) that
   $ \colim_{BC_{p}}\pi_{*}((KU/\ell)^{\otimes C_{p}})$  is also a sum of copies of $\Z/\ell \Z$. We must determine the number of such summands.

 We now note that 
  $\Z \xrightarrow{\ell} \Z$ is equivalent to $\Z[p^{-1}]\stackrel{\ell}{\to} \Z[p^{-1}]$.  Since now $p$ is invertible on the chain groups, we can interchange the order of taking homology and $C_{p}$-coinvariants on $(\Z[p^{-1}]\stackrel{\ell}{\to} \Z[p^{-1}])^{\otimes C_{p}}$.

  Let us start with then case $p=2$. Then the resulting cochain complex is (starting in degree $-2$)
  $$\Z[2^{-1}]\to \Z[2^{-1}]\oplus \Z[2^{-1}]\to \Z[2^{-1}]\ ,$$ where
  $C_{2}$ acts trivially on chain group in degree $0$, by the flip of summands on the group  in degree $-1$, and by the sign representation on the summand in degree $-2$. After taking coinvariants, we get a chain complex
 $\Z[2^{-1}]\stackrel{d_{-1}}{\to}\Z[2^{-1}]$ whose differential is multiplication by $\ell$.
We conclude that the odd cohomology vanishes and the even cohomology is given by $\Z/\ell$.
This shows that
$$ \colim_{BC_{2}} \K(\widehat \Res^{C_{2}}(\beins_{\KK}^{\otimes C_{2}}))\simeq  KU/\ell$$
 as asserted.
  
  We now consider the case $p>2$. 
  We write  $$(\Z[p^{-1}]\to \Z[p^{-1}])^{\otimes C_{p}}\cong (Z_{-p}\stackrel{d_{-p}}{\to}\dots Z_{1-p} \stackrel{d_{1-p}}{\to} \dots\to  Z_{-1} \stackrel{d_{-1}}{\to} Z_{0})\ ,$$
  where $Z_{-k}$ are free $\Z[p^{-1}]$-modules with a $C_{p}$-action. In detail  $$Z_{-k}\cong \bigoplus_{F'\in \cP_{k}(C_{p})} \Z[p^{-1}]^{\otimes C_{p}\setminus F'}\otimes \Z[p^{-1}]^{\otimes  F'}  \cong \left\{\begin{array}{cc} \triv^{C_{p}}(\Z[p^{-1}] )&k=0,p\\  \Ind^{C_{p}}(\Z[p^{-1}]^{\oplus \frac{1}{p} \binom{p}{k}})&k\in \{1,\dots,p-1\}  \end{array} \right.\ ,$$ where $\frac{1}{p} \binom{p}{k}$ is the number of $C_{p}$-orbits in $\cP_{k}(C_{p})$ for   $k$ in$ \{1,\dots,p-1\} $.
 The $C_{p}$-coinvariants  $Z_{-k,C_{p}}:=\colim_{BC_{p}} Z_{-k}$ of the chain groups  are therefore given by 
\begin{equation}\label{dasfoiahjiosdffqf}Z_{-k,C_{p}}\cong  \left\{\begin{array}{cc} \Z[p^{-1}] &k=0,p\\   \Z[p^{-1}]^{\oplus \frac{1}{p} \binom{p}{k}}&k\in \{1,\dots,p-1\}  \end{array} \right.\  .
\end{equation}    If we divide the differentials of the resulting chain complex 
 \begin{equation}\label{ewqdqwedewdqewd}0\to Z_{-p,C_{p}}\stackrel{d_{-p}}{\to} Z_{1-p,C_{p}}\stackrel{d_{1-p}}{\to} \dots\stackrel{d_{-1}}{\to} Z_{0,C_{p}}\to 0
\end{equation} 
  of coinvariants by $\ell$,
 then we get the analogous complex of coinvariants for $(\Z[p^{-1}] \xrightarrow{1} \Z[p^{-1}] )^{\otimes C_{p}}$ which is clearly exact.  Since the chain groups are $\ell$-torsion free,    we can conclude that the homology of 
 the complex \eqref{ewqdqwedewdqewd} is given by sums of copies of $\Z/\ell \Z$.

 It remains to determine the ranks. 
 To this end, we calculate the ranks $r_{i}$ of the differentials $d_{-i}$. Since the complex  \eqref{ewqdqwedewdqewd}  is exact after inverting $\ell$, we see that 
 $$r_{j}=\sum_{i=0}^{j-1} (-1)^{i-j+1} \rk Z_{-i,C_{p}}\ .$$
 Inserting \eqref{dasfoiahjiosdffqf}, we get  for $j<p$
 $$r_{j}=  \frac{(-1)^{j+1} }{p} \sum_{i=1}^{j-1} (-1)^{i} \binom{p}{i}+(-1)^{j+1}= \frac{\binom{p-1}{j-1}  +(-1)^{j}(1-p) }{p}\ .$$
We conclude that 
  $$\rk H_{-j}(\eqref{ewqdqwedewdqewd})= \left\{\begin{array}{cc}   \frac{\binom{p-1}{j}  +(-1)^{j-1}(1-p)  }{p}&j\in \{0,\dots,p-1\}\\ 0& j=p \end{array} \right.$$
 We now must add up the even and odd parts. 
 We  use 
 $$\sum_{j=0}^{\frac{p-1}{2}} \binom{p-1}{2j}=\sum_{j=0}^{\frac{p-3}{2}} \binom{p-1}{2j+1}=2^{p-2}$$
and get the asserted numbers \begin{eqnarray*}
n_{0}(p) &=&\frac{2^{p-2} +  \frac{(p+1)(p-1)}{2}  }{p}\\
 n_{1}(p)  &=& \frac{2^{p-2}  - \frac{(p-1)^{2}}{2}   }{p}\ .
\end{eqnarray*}
  \end{proof}
     \begin{prop}[\cref{ewfoijqwofdqewdqewdqewdq}]\label{ewfoijqwofdqewdqewdqewdq1}
     For a prime $p$ and $n$ in  $\nat$, we have an equivalence
$$\colim_{BC_{p}} \K( \widehat \Res^{C_{p}}((\beins_{\KK}/p^{n})^{\otimes C_{p}}))\simeq \K((\beins_{\KK}/p^{n})^{\otimes C_{p}}\rtimes_{r}C_{p})\ .$$
\end{prop}
\begin{proof}
In view of \eqref{vwvssfdvsfvs} it suffices to  show that  
\begin{equation}\label{ihjoifadfef}\Cof^{C_{p}}((\beins_{\KK}/p^{n})^{\otimes C_{p}} )\simeq 0\ .
\end{equation}
As a consequence of \cref{ewgojwoepgfrefwfwerfwrf},  we have an equivalence 
$$ \Cof^{C_{p}}((\beins_{\KK}/p^{n})^{\otimes C_{p}} )\simeq  \Cof^{C_{p}}(\beins_{\KK}^{\otimes C_{p}}\ )/p^{n}\simeq \Cof^{C_{p}}(\beins_{\KK^{C_{p}}})/p^{n}\ .$$
By \cref{wekopgwerfrewfw}, we have 
$$\Cof^{C_{p}}(\beins_{\KK^{C_{p}}})/p^{n}\simeq R(G)/p^{n}[\xi^{-1}]\ .$$
By \cref{wtjigowefrfwref}, we know that  $\xi^{2}$ is divisible by $p$ and consequently
$$R(G)/p^{n}[\xi^{-1}]\simeq 0\ .$$
Combining the last three equivalences we obtain \eqref{ihjoifadfef}. 
%
%
%
%
\end{proof}

 \bibliographystyle{alpha}
\bibliography{forschung2021}

\end{document}